\documentclass[format=acmsmall, review=false]{acmart}
\usepackage{acm-ec-26}
\usepackage{booktabs} % For formal tables
\usepackage[ruled]{algorithm2e} % For algorithms

\SetAlFnt{\small}
\SetAlCapFnt{\small}
\SetAlCapNameFnt{\small}
\SetAlCapHSkip{0pt}
\IncMargin{-\parindent}

\makeatletter
\AtBeginDocument{%
  \fancypagestyle{firstpagestyle}{%
    \fancyhf{}%
    \fancyfoot[L]{}
  }%
}
\makeatother

\usepackage[utf8]{inputenc} % allow utf-8 input
\usepackage[T1]{fontenc}    % use 8-bit T1 fonts
\usepackage{hyperref}       % hyperlinks
\usepackage{booktabs}       % professional-quality tables
\usepackage{amsfonts}       % blackboard math symbols
\usepackage{nicefrac}       % compact symbols for 1/2, etc.
\usepackage{microtype}      % microtypography
\usepackage{xcolor}         % colors
\usepackage{xr-hyper}
\usepackage{amsthm}
\usepackage{amsmath}
\usepackage{float}
\usepackage{mathrsfs}
\usepackage{amsfonts}
\usepackage{enumitem}
\usepackage{tikz}
\usetikzlibrary{shapes.geometric, arrows}
\usetikzlibrary{arrows.meta}
\usepackage{color}
\usepackage{cleveref}
\crefname{section}{§}{§§}
\Crefname{section}{§}{§§}

\usepackage{multirow}
\usepackage{array} 

\usepackage{threeparttable}
\usepackage{rotating}
\usepackage{hhline}
\allowdisplaybreaks

\def\bx{\boldsymbol x}

\def\bv{\boldsymbol \nu}

\def\bp{\boldsymbol p}

\def\bG{\boldsymbol G}
\def\bV{\boldsymbol V}
\def\bE{\boldsymbol E}

\def\blam{\boldsymbol \lambda}

\def\blambda{\boldsymbol \lambda}
\def\bdelta{\boldsymbol \delta}

\def\eps{\epsilon}

\def\cP{\mathcal{P}}
\def\cK{\mathcal{K}}
\def\cN{\mathcal{N}}

\def\cI{\mathcal{I}}

\def\ps-RAM{\mathcal{P}_{\textsc{s-RAM}}}

\newtheorem{assumption}{Assumption}

\setcitestyle{authoryear}

\title[Learning in Combinatorial Choice]{Going from a Representative Agent to Counterfactuals in Combinatorial Choice}
\author{Yanqiu Ruan}
\authornote{Institute of Operations Research and Analytics, National University of Singapore. \texttt{yanqiu.ruan@nus.edu.sg}}

\author{Karthyek Murthy}
\authornote{Industrial and Systems Engineering, Viterbi School of Engineering, University of Southern California. \texttt{karthyek@usc.edu}}

\author{Karthik Natarajan}
\authornote{Engineering Systems and Design, Singapore University of Technology and Design. \texttt{karthik\_natarajan@sutd.edu.sg}}

\begin{abstract}%
  We study decision-making problems where data comprises historical observations from a collection of binary polytopes, capturing aggregate information stemming from various combinatorial selection environments.  We propose a nonparametric approach for counterfactual inference in this setting based on a representative agent model, where the available data is viewed as arising from maximizing separable concave utility functions over the respective binary polytopes. Our first contribution is to precisely characterize the selection probabilities representable under this model and show that verifying the consistency of any given aggregated selection dataset reduces to solving a polynomial-sized linear program. Building on this characterization, we develop a nonparametric method for counterfactual prediction. When data is inconsistent with the model, finding a best-fitting approximation for prediction reduces to solving a compact mixed-integer convex program. Numerical experiments based on synthetic data demonstrate the method’s flexibility, predictive accuracy, and strong representational power even under model misspecification.
\end{abstract}

\begin{document}

% Title page for title and abstract only.
\begin{titlepage}

\maketitle

\makeatletter
% 清空首页底部由 acmart 输出的地址与声明区域
\gdef\@authorsaddresses{}
\gdef\@copyrightpermission{}
\gdef\@ACM@copyrighttext{}
\makeatother

% Optionally include a table of contents
\vspace{-0.1em}
\setcounter{tocdepth}{2} % adjust to 1 if desired
\tableofcontents

\end{titlepage}

% Paper body
\section{Introduction}\label{sec:intro}
Many decision-making problems involve agents selecting from complex combinatorial feasible sets, such as paths in networks, assignments in matching markets, or constrained item bundles. Given data on such choices, a central challenge is to infer models that both explain observed behavior and enable prediction in counterfactual environments not previously encountered. For instance, using historical congestion data from an evolving city transportation network, can we infer how the introduction of a new link would reshape congestion patterns? Such inferential capabilities are essential for understanding system behavior broadly and for informing downstream network and system design. 
%In such settings, observed decisions are often solutions to implicit optimization problems under uncertainty. Modeling the  behavior of such observed decisions is critical for  inference in counterfactual environments unobserved in the past and further downstream network or  system design questions. 

This paper builds on the following premise: By modeling observed selection data as outcomes of an implicit optimization problem, the induced structure on agent decisions can serve as a vehicle for predicting behavior in new environments.

An elementary, well-studied example of such modeling arises in capturing consumer choice over different product assortments. In particular, a flexible and well-established framework for modeling discrete choice behavior, grounded in microeconomic theory, is the representative agent model (RAM) or equivalently, the perturbed utility model~\citep{fudenberg2015stochastic, mcfadden2012theory, allen2019identification}, in which a representative agent chooses among $n$ alternatives on behalf of the entire population. Specifically, a representative agent chooses a consumption (or) choice probability vector $\bp$ %from the simplex $\mathcal{X} = \{ \bx \in [0,1]^n: \mathbf{1}^\intercal \bx = 1\},$ 
that solves the utility maximization $$ \bp = \arg\max \{ \bv^\top \bx - c(\bx): \bx \in [0,1]^n, \mathbf{1}^\intercal \bx = 1 \};$$ here the utility function is a linear function perturbed by subtracting a convex function that rewards the agent for diversifying away from the extremal $0-1$ vertices capturing deterministic choices. %In turn, this diversification is necessary for capturing datasets exhibiting  probabilistic selection.  
While typical applications of  RAM   focus on this  discrete choice setting, %where $\mathcal{X} = \{\bx \in \{0,1\}^n: 1^\intercal \bx = 1\}$ and the resulting $\text{conv}(\mathcal{X})$ is a simplex, 
we study its utility %and  significance 
in more challenging  combinatorial environments.
%, enabling its use in structured decision problems such as networks, assignments, and constrained subset selection. 

In particular, we consider settings where the observed dataset comprises tuples of the form $\{(\mathcal{X}^{(k)}, \bp^{(k)}): k \in [K]\}$, where $K$ is the total number of historical observations, $\mathcal{X}^{(k)} \subseteq \{0,1\}^n$ is a  combinatorial selection set, and   %(which can be significantly more complex than those arising in discrete choice)  %contained in $\mathcal{X} = \{ \bx \in \{0,1\}^n \mid \boldsymbol{A}\bx = \bb\},$ 
and  the respective $\bp^{(k)}$ lies in the convex hull of $\mathcal{X}^{(k)}.$ While the set $\mathcal{X}^{(k)}$ may be understood to be capturing the feasible set in a given network (or) a combinatorial selection environment, the respective $\bp^{(k)}$ can be understood as the resulting selection data aggregated/averaged across individual agents. For instance, the set $\mathcal{X}^{(k)}$ may correspond to the 
set of directed paths available in a transportation network configuration (captured by the flow conservation constraints) and the respective $\bp^{(k)}$ may capture the empirical aggregated flow of agents traversing these paths.  %in a market for matching agents from $A$ to users in $B$, the set $\mathcal{X}^{(k)}$ may denote the set of feasible assignments available at a given instance and the resulting $\bp^{(k)} = (p_{ij}^{k}: (i,j) \in A \times B)$ may capture the likelihood of a randomly selected agent $i \in A$ being matched with user $j \in B.$ 
The feasible set instances, denoted by  $\mathcal{X}^{(1)},\ldots,\mathcal{X}^{(K)},$ may differ among themselves due to different collections of variables being forced to be zero at different instances, say, due to link dysfunctions/unavailability in networks (or) unavailable pairings at a given point of time.  
%points from a collection of $0-1$ polytopes with some variables forced to be 0, such as feasible sets in networks with certain links disabled, or assignment problems where some pairs are disallowed. 

How might we model such a dataset $\{(\mathcal{X}^{(k)}, \bp^{(k)}): k \in [K]\},$ capturing aggregated selection statistics in combinatorial settings, so that we can perform counterfactual predictions for selections in a new variant, $\mathcal{X}^{\text{new}},$ differing from all the selection sets  $\mathcal{X}^{(1)},\ldots,\mathcal{X}^{(K)}$ observed in the past? 
The significance of this research question spans diverse settings, as indicated in Table 1 below.  %We illustrate more examples in the table below.
%\textcolor{red}{add a very concrete example here, then explain with this table}
\begin{table}[h!]
\label{tab:eg_app}
      \caption{Example application environments} 
\scalebox{0.8}{\begin{tabular}{|c|cc|cc|}
\hline
\multirow{2}{*}{\begin{tabular}[c]{@{}c@{}}Application \\ (Problem)\end{tabular}}                                                                                                               & \multicolumn{2}{c|}{\begin{tabular}[c]{@{}c@{}}Historial Observed Data \\ $\{(\mathcal{X}^{(k)} ,\bp^{(k)}): k\in [K]  \}$\end{tabular}}                                                                                                                                                   & \multicolumn{2}{c|}{Counterfactual Prediction}                                                                                                                                                                                     \\ \cline{2-5} 
                                                                                                                                                     & \multicolumn{1}{c|}{$\mathcal{X}^{(k)}$}                                                               & $\bp^{(k)}\in \text{conv}(\mathcal{X}^{(k)} )$                                                                                                                                    & \multicolumn{1}{c|}{Unseen Scenario $\mathcal{X}^{\text{new}}$}                                                                    & Prediction Task                                                                                \\ \hline
\begin{tabular}[c]{@{}c@{}}prediction of traffic \\ (shortest path)\end{tabular}                                                            & \multicolumn{1}{c|}{\begin{tabular}[c]{@{}c@{}}discrete feasible set \\ of $s-t$ paths\end{tabular}} & \begin{tabular}[c]{@{}c@{}}historial flows of the given\\  traffic network $\mathcal{X}^{(k)}$\end{tabular}                                                                       & \multicolumn{1}{c|}{\begin{tabular}[c]{@{}c@{}}modified network with \\  new constructed links\end{tabular}}                           & \begin{tabular}[c]{@{}c@{}}how travel patterns \\ would shift\end{tabular}                     \\ \hline
\begin{tabular}[c]{@{}c@{}}matching\\ (linear assignment )\end{tabular}                                                                              & \multicolumn{1}{c|}{\begin{tabular}[c]{@{}c@{}}discrete feasible set \\ of assignments\end{tabular}}   & \begin{tabular}[c]{@{}c@{}}historical pairwise matching  \\ rates in the given\\   assignment set $\mathcal{X}^{(k)}$\end{tabular}                                                          & \multicolumn{1}{c|}{\begin{tabular}[c]{@{}c@{}}a new assignment set\\ with some pairs  being\\  restricted/permitted\end{tabular}} & \begin{tabular}[c]{@{}c@{}}how would the \\ expected matching \\ structure shift\end{tabular} \\ \hline
\begin{tabular}[c]{@{}c@{}}project evaluation\\ (longest path )\end{tabular}                                                                         & \multicolumn{1}{c|}{\begin{tabular}[c]{@{}c@{}}discrete feasible set \\ of $s-t$ paths \end{tabular}} & \begin{tabular}[c]{@{}c@{}}the likelihoods of the activities \\ being a critical bottleneck \\ in meeting the milestones \\ of the given project $\mathcal{X}^{(k)}$\end{tabular} & \multicolumn{1}{c|}{a new project network}                                                                                         & \begin{tabular}[c]{@{}c@{}}how likely will \\ an activity \\ be critical\end{tabular}         \\ \hline
\begin{tabular}[c]{@{}c@{}}feature/asset selection \\ with cardinality \\ or matroid constraints\\ (constrained \\ subset selection)\end{tabular} 
& 
\multicolumn{1}{c|}{uniform matriod}                                                                   & \begin{tabular}[c]{@{}c@{}}the inclusion probabilities\\  of features/assets of  \\ the given constrained set $\mathcal{X}^{(k)}$\end{tabular}                                    & \multicolumn{1}{c|}{a new constrained set}                                                                                         & \begin{tabular}[c]{@{}c@{}}how do the \\ inclusion \\ probabilities \\ change\end{tabular}    \\ \hline
\end{tabular}}
\end{table}
%\textcolor{red}{for uniform matroid constraint, like $\sum_i x_i <= 5$, $x_i$ can be greater than 1. should we change the wording - inclusion probabilities here?}

% \noindent
% \textbf{Prediction of traffic via shortest path formulations:} Given aggregate traffic flows in historical road networks, can we predict how travel pattern would shift with the construction of new road links? 
% \noindent
% \textbf{Matching via linear assignment formulations:} From historical matchings in two-sided markets, %(e.g., job candidates to firms), 
% how would the expected matching structure shift if some pairings are restricted/permitted? \\
% \noindent
% \textbf{Project evaluation via longest path formulations:} Given various project networks comprising different milestones,  activities,  and the likelihoods of the activities being a critical bottleneck in meeting the milestones, how likely will an activity be critical when presented with a new project?\\ 
% \noindent
% \textbf{Constrained Subset Selection:} In feature or asset selection problems under cardinality or matroid constraints, how do inclusion probabilities change when the constraint set is altered?

In this work, we propose to respond to such challenges by  modeling the selection data $\{(\mathcal{X}^{(k)}, \bp^{(k)}): k \in [K]\}$ via the aforementioned representative agent model $\bp^{(k)} = \arg\max_{\bx \in \text{conv}(\mathcal{X}^{(k)})} \{ \bv^\top \bx - c(\bx)\},$ for $k \in [K].$ In particular, we focus on the case where $c(\bx)$ is in a separable additive form and we call the resulting model as the Separable Representative Agent Model (abbreviated as s-RAM). Formally, in the s-RAM, the representative agent selects a vector %$\bp \in \text{conv}(\mathcal{X})$ that solves $\arg\max_{\bx \in \text{conv}(\mathcal{X})} \{ \bv^\top \bx - \sum_{j \in [n]} c_j(x_j)\},$
\begin{align}
    \bx^* = \arg\max_{\bx \in \text{conv}(\mathcal{X})} \left\{ \bv^\top \bx - \sum_{j \in [n]} c_j(x_j) \right\}, \tag{s-RAM} \label{model:sram}
\end{align}
where $\bv \in \mathbb{R}^n$ is a vector of deterministic utilities and $\{c_j(x_j): j \in [n]\}$ are  convex perturbation functions. Besides proposing this model, one of our  key contributions is to precisely characterize the collection of $(\mathcal{X}^{(k)}, \bp^{(k)})$ tuples which can be modeled by s-RAM. %In other words, given 
Such exact characterizations have been critical in  discrete  choice modeling, both in leading to crisp understanding of the representative abilities of choice models and in performing counterfactual inference over product assortments unseen in the datasets; see, \cite{block1959random, farias2013nonparametric, jagabathula2019limit, fudenberg2015stochastic,ruan2022nonparametric}. Similarly, our modeling of combinatorial selection data via s-RAM leads towards  the following contributions: 

\begin{enumerate}
    \item \textbf{New modeling framework.} We propose a new nonparametric model for making counterfactual predictions in combinatorial selection environments, building on the s-RAM and the results developed below (see Figure \ref{fig:flowchat} for a flow chart). The approach is  tractable and applicable to a broad class of combinatorial problems, enabling downstream inference for any quantity of interest on new  instances.

     \item  \textbf{Polynomial-time verifiable s-RAM consistency.} We provide necessary and sufficient conditions for observed data in a collection of  $0 - 1$ polytopes to be representable by s-RAM. This reduces to solving a polynomial-sized linear program, enabling efficient model verification.

\item \textbf{Nonparametric counterfactual prediction.} When the data is consistent with s-RAM, we develop a nonparametric data-driven approach that builds on the exact characterization of s-RAM to support efficient prediction over new feasible sets using arbitrary functions of interest. This approach mitigates the risk of misspecification and enables end-to-end learning for s-RAM.
When the data is inconsistent with s-RAM, we show that a best fitting approximation can be computed via a compact mixed integer convex program, which can be further used for prediction.

\item \textbf{Empirical performance.} Through experiments on synthetic data from network problems including shortest path, longest path, assignment, and constrained subset selection, we demonstrate the practical flexibility, predictive accuracy, and strong representational power of our framework, even under model misspecification.
s

\end{enumerate}
%We propose a new nonparametric framework for counterfactual prediction in choice problems with combinatorial structures based on the s-RAM described above and the results below (see Figure \ref{fig:flowchat} for a flow chart). Our approach is general, tractable, and applicable to a broad class of combinatorial problems and enables downstream inference for any quantity of interest over new combinatorial structures.

\begin{figure}[h!] 
\vspace{-10pt}
\tikzstyle{startstop} = [rectangle, rounded corners, 
minimum width=3cm, 
minimum height=1cm,
text centered, 
draw=black, 
fill=white!30]

\tikzstyle{io} = [trapezium, 
trapezium stretches=true, % A later addition
trapezium left angle=70, 
trapezium right angle=110, 
minimum width=4.5cm, 
minimum height=1cm, 
text width=4.5cm, 
%text centered, 
draw=black, fill=white!30]

\tikzstyle{process} = [rectangle, 
minimum width=5cm, 
minimum height=1cm, 
%text centered, 
text width=5cm, 
draw=black, 
fill=white!30]

\tikzstyle{process2} = [rectangle, 
minimum width=4cm, 
minimum height=1cm, 
%text centered, 
text width=4cm, 
draw=black, 
fill=white!30]

\tikzstyle{decision} = [diamond, 
minimum width=3cm, 
minimum height=1cm, 
text width=3cm, 
text centered, 
draw=black, 
fill=white!30]

\tikzstyle{arrow} = [thick,->,>=stealth]

\tikzset{global scale/.style={
    scale=#1,
    every node/.append style={scale=#1}
  }
}

\begin{center}
\scalebox{0.8}{\begin{tikzpicture}[node distance=2cm]
\node (dec1)[decision]{ \textbf{Step 1:} Is data $\bp^\cK$  represented by s-RAM? };
%\node (in1) [io, above of = dec1, yshift = 1.5cm] {Inputs: historical choice data $\bp_\cK$, new assortment $A,$price vector $\boldsymbol{r}_A$};
\node (pro1)[process, right of = dec1, xshift = 3.5cm, yshift = 1.2cm ] {Set s-RAM  probabilities $\bx^\cK_* = \bp^\cK$}; 
\node (pro2) [process, right of = dec1, xshift = 3.5cm, yshift = -1.2 cm ]{\textbf{Step $1^\prime$: }Find the best fitting s-RAM probabilities $\bx^\cK_*$ to $\bp^\cK$};
\node (pro3) [process, right of = dec1, xshift = 9.5 cm ]{\textbf{Step 2:} Nonparametric prediction for new polytope $\mathcal{X}^{\text{new}}$ with $\bx^\cK_*$};
%\draw [arrow] (in1) -- (dec1);
\draw (2.25,0) -- (2.5,0);
\draw (2.5,0) -- (2.5,-1);
\draw (2.5,0) -- (2.5,1);
\draw [arrow] (2.5,-1) |- node[anchor=east] {NO} (pro2);
\draw [arrow] (2.5,1) |- node[anchor=east] {YES} (pro1);
\draw (8.15,1) -- (8.5,1);
\draw (8.15,-1) -- (8.5,-1);
\draw (8.5,1) -- (8.5,0);
\draw (8.5,-1) -- (8.5,0);
\draw [arrow] (8.5,0) -- (pro3);
\node[text width=4cm] at (0, - 2.8) {Polynomial-sized LP};
\node[text width=4cm] at (0.8, - 3.2) {(Theorem  \ref{thm:mdm-general-feascon}) };
%\node[text width=4cm] at (0.8, - 3.7) {(Theorem  \ref{thm:feascon-gmdm}) };
\node[text width=4cm] at (6.2, - 2.2)  {Compact MICP};
\node[text width=4cm] at (6.3, - 2.6) {(Proposition  \ref{prop:micp}) };
%\node[text width=4cm] at (6.3, - 3.1) {(Proposition  \ref{prop:gmdm_micp}) };
\node [text width=4cm] at (12.2,  -1)  {Compact IP};
\node[text width=4cm] at (12.3, - 1.4) {(Proposition \ref{prop:wc-obj-milp})};
%\node[text width=4cm] at (12.3, - 1.9) {(Proposition \ref{model:micp_gmdm_revenue_prediction})};
\end{tikzpicture}}
\vspace{-17pt}
\end{center}
\caption{ A new framework for counterfactual prediction with a summary of main contributions }
\label{fig:flowchat}
\end{figure}

The remainder of the paper is organized as follows. Section \ref{sec:prelim} introduces the problem, the model, and the related literature. Section \ref{sec:representation} provides a complete characterization of the data representable by the s-RAM model through polynomial-time verifiable necessary and sufficient conditions. In Sections \ref{sec:prediction} and \ref{sec:estimation}, we illustrate how these characterizations can be used to generate predictions in counterfactual environments and to estimate the best-fitting model, respectively.  Section \ref{sec:experiment} discusses results from numerical experiments on representability and counterfactual predictions.

\section{Model Description and Related Literature}\label{sec:prelim}
%We first formally present the problem setup and then review the related literature. %to the research question of this paper.

\subsection{Problem Description}
\label{sec:prob-desc}
%\textcolor{red}{(Notation of the combinatorial problem and the LP, adding script z for the integer discrete set)}
Consider $\mathcal{X} := \{\bx \in \{0,1\}^n \mid \boldsymbol{A}\bx \leq \boldsymbol{b}\},$ with $\boldsymbol{A} \in \mathbb{R}^{m \times n}$ and $\boldsymbol{b}\in \mathbb{R}^m$ characterizing the combinatorial problem of interest. Let $K$ be the number of historical observations and $\cK = \{1,\ldots,K\}$. Suppose for any $k\in \cK$, we observe $\bp^{(k)} \in \text{conv} (\mathcal{X}^{(k)})$ where $\mathcal{X}^{(k)} \subseteq \mathcal{X}$ and
 and $\boldsymbol{b} \in \mathbb{R}^m$. Specifically, 
\[
\mathcal{X}^{(k)} : = \{\bx \in \{0,1\}^n \mid \boldsymbol{A}\bx \leq \boldsymbol{b}, x_j = 0, \forall j \notin \cI^{(k)} \},
\]
where $\emptyset \subset \cI^{(k)} \subseteq [n]$. Thus, the set $\{ \text{conv}( \mathcal{X}^{(k)} ): k \in \cK  \}$ is the collection of historical $0-1$ polytopes. Let $\cP^\cK = \{\bp^{(k)}: k \in \cK \}$ and $\cI^\cK = \{ \cI^{(k)}: k\in \cK \}$ capture the collection of historical observations.   

Given any new structured polytope $\text{conv}(\mathcal{X}^{\text{new}}),$ where $\mathcal{X}^{\text{new}} := \{\bx \in \{0,1\}^n \mid \boldsymbol{A}\bx \leq \boldsymbol{b}, x_j = 0, \forall j \notin \cI^\text{new} \}$,  with no prior data, i.e., $ \cI^{\text{new}} \notin \cI^\cK,$ our goal is to:
\begin{enumerate}
    \item  verify whether $\cP^\cK$ is consistent with the s-RAM model hypothesis, and
\item  develop a nonparametric data-driven approach for making predictions for the polytope $\text{conv}(\mathcal{X}^{\text{new}})$ that is unseen from the historical data.
\end{enumerate}

\subsection{The s-RAM model}
\begin{assumption}\label{asp:convex_perturbation}
   The functions $c_j:[0,1]\to \mathbb{R}$ for $j \in [n]$ are univariate, strictly convex and continuously differentiable perturbation functions. 
\end{assumption}
{\color{black}Let $\boldsymbol{s} = \boldsymbol{b} - \boldsymbol{A}\bx $ be slack variables for constraints of $\mathcal{X}$. Then $\mathcal{X} = \{\bx \in \{0,1\}^n \mid \boldsymbol{A}\bx + \boldsymbol{s} =  \boldsymbol{b}, \boldsymbol{s} \geq \boldsymbol{0}\} .$ We assume that the convex hull of $\mathcal{X}$ is a $0-1$ polytope, as stated in Assumption \ref{asp:convexhull} below.
\begin{assumption}\label{asp:convexhull}
   The convex hull of $\mathcal{X}$ 
$\text{conv}(\mathcal{X}) = \{\bx \in \mathbb{R}^n \mid \boldsymbol{A}\bx + \boldsymbol{s} =  \boldsymbol{b}, \boldsymbol{0}\leq \bx \leq \boldsymbol{e}, \boldsymbol{s} \geq \boldsymbol{0} \}.$
\end{assumption}
%\textcolor{red}{[write downn explictly the assumption convex hull of $\mathcal{X}$ then it leads to the convex hull of $\mathcal{X}^{(k)}$  then goes to the equation (1) ]}.
Then, for $k\in \cK$,  the convex hull of $\mathcal{X}^{(k)} $ can be written as 
\[
\text{conv}(\mathcal{X}^{(k)}) = \{\bx \in \mathcal{R}^n \mid \boldsymbol{A}\bx + \boldsymbol{s} =  \boldsymbol{b}, \boldsymbol{0}\leq \bx \leq \boldsymbol{e},  x_j = 0, \forall j \notin \cI^{(k)}, \boldsymbol{s} \geq \boldsymbol{0} \}.
\] 
Thus, the set $\{ \text{conv}( \mathcal{X}^{(k)} ): k \in \cK  \}$ is the collection of historical $0-1$ polytopes. Note that all the examples in the table in  Section \ref{sec:intro} satisfy Assumption \ref{asp:convexhull}. }
 For $k\in \cK,$ the s-RAM solves the  convex program $\max_{\bx \in \text{conv}(\mathcal{X}^{(k)})} \big\{ \bv^\top \bx - \sum_{j\in \cI^{(k)}} c_j(x_j) \big\}$ to maximize the additively perturbed expected utility, which in turn is equivalent to 
\begin{align}
\label{s-RAM_k}
    \max_{\bx} \left\{ \bv^\top \bx - \sum_{j\in \cI^{(k)}} c_j(x_j) \middle| \ \
    \begin{aligned}
    &\sum_{j\in \cI^{(k)}} a_{ij} x_j + s_i =  b_i \quad  \forall  i \in [m], \\    
    & \quad 0 \leq x_j \leq 1  \quad \forall  j \in \cI^{(k)},\\
    & \quad s_i \geq 0 \quad\forall i \in [m]
    \end{aligned}
       \right\}. 
\end{align}
Any $\bx^*$ that maximizes \eqref{s-RAM_k} if and only if the following optimality conditions for \eqref{s-RAM_k} are satisfied:
\begin{align}
\begin{aligned}
     \nu_j  - c_j^\prime(x_j^*) + \sum_{i\in [m]} \lambda_i^k a_{ij} + \mu_j^k - w_j^k &= 0,\quad \forall j\in \cI^{(k)},\\
     \mu_j^k x_j^* = 0, \ w_j^k (1-x_j^*) & = 0, \quad \forall  j \in \cI^{(k)},\\
      \sum_{j\in \cI^{(k)}} a_{ij} x_j^* + s_i^* - b_i   = 0, \ \alpha_i s_i^* &= 0, \quad \forall i \in [m],\\
     %0 \leq x_j^* & \leq 1, \quad \forall j \in \cI^{(k)},
    0 \leq x_j^*  \leq 1, \ \mu_j^k, w_j^k &  \geq 0,\quad \forall j \in \cI^{(k)},\\
    s_i^*\geq 0 , \alpha_i &\geq 0, \quad \forall i \in [m],  \label{s-RAM-optcon}
\end{aligned}
\end{align}
where %$(\lambda_i^k: i \in [m] )$, $(\mu_j^k: j\in \cI^{(k)})$ and $(w_j^k: j\in \cI^{(k)})$
$\lambda_i^k$ and $\alpha_i$ for all $i \in [m]$, $ \mu_j^k$ and $ w_j^k$ for all $j\in \cI^{(k)}$ 
are dual variables of \eqref{s-RAM_k}.  %associated with constraints of \eqref{s-RAM_k}. 

For each observation $k$, since we observe that $\bp^{(k)} \in \text{conv} (\mathcal{X}^{(k)})$, we can compute $\bar{s}_i = b_i - \sum_{j\in \cI^{(k)}} a_{ij} p_j^{(k)}$, for $i\in [m].$ Then the first goal described in Section \ref{sec:prob-desc} above translates to verifying if $(\bp^{(k)},\bar{\boldsymbol{s}}^{(k)})$ satisfies the conditions in \eqref{s-RAM-optcon} for all $k\in \cK$.

%Our work complements the existing literature by shifting focus from model-based estimation in specific applications to a general, data-driven framework for verifying consistency and predicting with s-RAM in arbitrary combinatorial settings. Unlike prior studies that often assume known perturbation functions in s-RAM (except for \cite{liu2022product}), our work provides an assumption-light framework that supports nonparametric counterfactual prediction.

% \textbf{Related Literature.} Our paper is closely connected to the machine learning literature on learning and inference in combinatorial optimization, as well as to the choice modeling literature on combinatorial choice problems.

\subsection{Related Literature}
Our paper is closely connected to the machine learning literature on learning and inference in combinatorial optimization, and  as well as to the choice modeling literature on combinatorial selection problems.

%(ii) learning and optimization on graphs. It also makes a contribution to research on choice modeling in general combinatorial settings.
%\textcolor{red}{[first discuss the perturbation method to ]}
\subsubsection{Learning and Inference in Combinatorial Optimization.} Effective structured learning and inference over discrete combinatorial models is challenging. A prominent line of work uses perturbation methods, including perturbed optimizers \citep{berthet2020learning}, perturbed structural prediction \citep{indelman2021learning}, and stochastic softmax tricks for combinatorial spaces \citep{paulus2020gradient}. These methods can be viewed as instances of RAMs with perturbation functions and are typically designed for general binary combinatorial optimization. Our work differs  in two aspects: (i) we focus on RAMs with separable perturbation functions in combinatorial problems, and (ii) we study 0–1 polytopes, the convex hulls of combinatorial problems, and develop a nonparametric, end-to-end s-RAM framework that remains amenable via standard integer programming solvers.  RAM with separable perturbation functions have been shown to yield good tractability in discrete choice while maintaining expressiveness \citep{ruan2022nonparametric}. 

Another stream combines neural networks with combinatorial optimization, using RAMs to make prediction differentiable \citep{djolonga2017differentiable,dalle2022learning,wilder2019melding}. 
While these methods rely on deep learning models to approximate prediction, our framework leverages the polyhedral structure of $0- 1$ polytopes and the separable structure of s-RAM for tractable, nonparametric counterfactual learning.

\subsubsection{Choice Modeling in Combinatorial Problems.}
%\textbf{Combinatorial Choice with s-RAM.}
A large body of work applies logit models in the random utility model family to network route choice, starting with the Multinomial Logit (MNL) model for traffic assignment \citep{dial1971probabilistic} and later route-based extensions that capture overlap and correlated utilities, including path-size, nested, and recursive logit formulations \citep{chu1989paired,cascetta1996modified,bekhor1999formulations,ben1999discrete,zhou2012c,mai2015nested,mai2017dynamic,mai2016method}. More recently, s-RAM–based models have been introduced to traffic network problems, with the Perturbed Utility Route Choice model and its extensions addressing equilibrium assignment \citep{fosgerau2022perturbed,yao2024perturbed,fosgerau2024substitution}. It has been applied beyond transportation, such as product and ancillary pricing in the airline industry \citep{liu2022product}.

While these studies are mainly on route-based models, our framework extends s-RAM to a link-based choice model within networks, while enabling a data-driven, nonparametric approach to counterfactual prediction in more general combinatorial settings without restrictive parametric assumptions.

\subsubsection{Random Utility Interpretation of s-RAM} The s-RAM is equivalent to the Marginal Distribution Model (MDM) proposed by \cite{natarajan2009persistency}, an important class of semi-parametric choice models. 
Let $\tilde{\epsilon}_1,\ldots,\tilde{\epsilon}_n$ be random variables with marginals $F_1,\ldots,F_n.$
%Consider random variables $\tilde{\epsilon}_1$ to $\tilde{\epsilon}_n$ and $F_j$ denote the marginal distribution of $\tilde{\epsilon}_j$ for all $j\in [n]$. 
Let $\Theta$ denote the set of joint distributions for $\tilde {\boldsymbol{\epsilon}}$ with the given marginals $F_1, F_2,\ldots, F_n$. %\textcolor{black}{Under the following assumption on the marginal distributions, the choice probability vector under MDM is unique.}
%\begin{assumption}
%Each random term $\tilde{\epsilon}_{j},$ where $j\in [n],$ is an absolutely continuous random variable with a strictly increasing marginal distribution $F_j(\cdot)$ on its support and $\mathbb{E}|\tilde{\epsilon}_{j}| < \infty$.\label{asp:general}
%\end{assumption}
\begin{lemma}[\citealt{natarajan2009persistency}] \label{lemma:mdm_cv}
Suppose $ \mathcal{X} \subseteq \{0,1\}^n,$ the random variables $\{\tilde{\epsilon}_{j}\}_{j\in [n]}$ are absolutely continuous random variables with strictly increasing marginal distribution $\{F_j(\cdot)\}_{j \in [n]}$ on their respective supports, and $\mathbb{E}|\tilde{\epsilon}_{j}| < \infty$.
%and Assumption \ref{asp:general} holds. 
Then 
\begin{align}
    \sup_{\theta \in \Theta} \mathbb{E}_\theta  \Big( \max \{ (\bv + \tilde{\boldsymbol{\epsilon}})^\top \bx : \boldsymbol{x}\in  \mathcal{X} \} \Big) = \sup \left\{ \sum_{j \in [n] } \nu_jx_j +
\sum_{j \in [n]} \int_{1-x_j}^{1} F_j^{-1}(t) \, dt \ \middle| \ \bx \in \text{conv}(\mathcal{X})
\right\}. \nonumber
\end{align}
\end{lemma}
Thus, by letting $c_j(x_j) = -\int_{1-x_j}^{1} F_j^{-1}(t)$ for all $j\in [n]$, the s-RAM and MDM can be seen as equivalent to each other. Consequently, the broad range of applications studied for MDM also apply to s-RAM. In particular, MDM has been shown to exhibit good empirical performance in various applications, such as traffic assignment \citep{ahipacsaouglu2019distributionally}, product pricing \citep{zhenzhen,liu2022product}, assortment planing \citep{sun2020unified}. The formulation of MDM has also become useful in deriving prophet inequalities for Bayesian online selection problems (see \citealt{feldman2021online}) and solving smoothed optimal transport formulations (see \citealt{tacskesen2022semi}).

\section{Step 1: Checking s-RAM Representability to The Data} \label{sec:representation}

In this section, we develop the necessary and sufficient conditions for the observed data given by a collection of historical 0-1 polytopes that are representable by s-RAM. This serves as the foundation for Step 1 of the the workflow in Figure ~\ref{fig:flowchat}.

\begin{theorem}[s-RAM Characterization with $0-1$ Polytopes]
\label{thm:mdm-general-feascon}
    Under Assumptions \ref{asp:convex_perturbation} and \ref{asp:convexhull}, there exists an s-RAM to represent $\cP^\cK$ if and only if there exists $(\lambda_i^k \in \mathbb{R}: i\in [m], k \in \cK)$,
    such that
\begin{align}
\begin{aligned}
     &\sum_{i\in [m]}\lambda_i^k a_{ij}  < \sum_{i\in [m]}\lambda_i^l a_{ij}  \quad \text{if} \quad p_{j}^k < p_{j}^l,  \\
     &\sum_{i\in [m]}\lambda_i^k a_{ij}   = \sum_{i\in [m]}\lambda_i^l a_{ij}  \quad \text{if}  \quad 0< p_{j}^k = p_{j}^l < 1.  \label{eq:mdm-feascon-order}
\end{aligned}
\end{align}
 for any distinct $k,l \in \cK$ and $j\in \cI^{(k)} \cap \cI^{(l)}$.
% and for any $k\in \cK$, for any constraint $i \in [m]$,
% \begin{align}
%     & \lambda_i^k = 0 \quad \text{if} \quad b_i - \sum_{j\in \cI^{(k)} } a_{ij} p_j^k  > 0.\label{eq:mdm-feascon-cs}
% \end{align}
As a result, checking whether the given data $\cP^\cK$ satisfies the s-RAM hypothesis can be accomplished by solving a linear program with $\mathcal{O}(mK)$ continuous variables and $\mathcal{O}(nK)$ constraints.  
\end{theorem}

For any collection $\cI^\cK,$ let $\ps-RAM(\cK)$ denote the collection of probabilities for $\cI^\cK$, which are representable by any s-RAM. Due to the characterization in Theorem \ref{thm:mdm-general-feascon}, we have the following succinct description for s-RAM: for any collection $\cI^\cK$, 
\[\ps-RAM(\cK) = \text{Proj}_{\bx} (\Pi_\cK) :=  \left\{ \bx: (\bx,\boldsymbol{\lambda}) \in 
\Pi_\cK \right\},\] where $\Pi_\cK$ is
defined as 
\begin{align}\label{eq:lifted-set}
\Pi _\cK
=
\left\{
(\bx,\boldsymbol{\lambda}) :
\begin{aligned}
    &\sum_{j\in \cI^{(k)}} a_{ij} x_j \leq b_i, \ \forall k \in \cK, \forall  i \in [m],  
    \ 0 \leq x_j \leq 1, \ \forall k \in \cK,  \forall  j \in \cI^{(k)},  \\ 
    &\sum_{i\in [m]}\lambda_i^k a_{ij}   < \sum_{i\in [m]}\lambda_i^l a_{ij}  
    \quad \text{if} \quad x_{j}^k < x_{j}^l, 
    \ \forall k,l \in \cK, \forall j\in \cI^{(k)} \cap \cI^{(l)}, \\  
    &\sum_{i\in [m]}\lambda_i^k a_{ij}   = \sum_{i\in [m]}\lambda_i^l a_{ij}  
    \quad \text{if} \quad 0< x_{j}^k = x_{j}^l <1, 
    \ \forall k,l \in \cK, \forall j\in \cI^{(k)} \cap \cI^{(l)}
\end{aligned}
\right\}.
\end{align}

% \begin{align}\label{eq:lifted-set}
%     \Pi _\cK = \left\{  (\bx,\boldsymbol{\lambda}) :  \,  
%     &\sum_{j\in \cI^{(k)}} a_{ij} x_j \leq b_i, \ \forall k \in \cK, \forall  i \in [m],  0 \leq x_j \leq 1, \ \forall k \in \cK,  \forall  j \in \cI^{(k)},  \nonumber \\ 
%     &\sum_{i\in [m]}\lambda_i^k a_{ij}   < \sum_{i\in [m]}\lambda_i^l a_{ij}  \quad \text{if} \quad x_{j}^k < x_{j}^l, \forall k,l \in \cK, \forall j\in \cI^{(k)} \cap \cI^{(l)}, \\  
%     & \sum_{i\in [m]}\lambda_i^k a_{ij}   = \sum_{i\in [m]}\lambda_i^l a_{ij}  \quad \text{if} \quad 0< x_{j}^k = x_{j}^l <1, \forall k,l \in \cK, \forall j\in \cI^{(k)} \cap \cI^{(l)} \nonumber
%     \right\}.  
% \end{align}
While Theorem \ref{thm:mdm-general-feascon} provides a general and abstract characterization of s-RAM representability, it may be difficult to interpret directly. We therefore present the following examples to show how this characterization simplifies in concrete combinatorial problems. Taking the longest path problem as an example, the abstract ordering conditions reduce to intuitive comparisons of edge probabilities to be consistent with the relative shadow prices across instances, leading to simple linear feasibility checks. 
This example provides intuition for the characterization, which in turn facilitates understanding of the prediction and estimation formulations introduced later. For the representability conditions of the other examples listed in Table 1, we refer the reader to Section \ref{app:example_theory} in the appendix.

\begin{example}[s-RAM Characterization for the Longest Path Problem]
Suppose that $\bG=(s\cup t \cup \bV, \bE)$ is the complete Directed Acyclic Graph (DAG) that describes all possible nodes and edges of the network, and we have flow data for a collection $K$ DAGs, denoted as $\bG^1,\bG^2, \cdots, \bG^K,$ with $\bG^k = (s\cup t \cup \bV^k, \bE^k)$, for every $k\in \cK$. For a graph $\bG^k$, let $p_{ij}^k \in [0,1]$ denote the flow of edge $(i,j)$ in  $\bG^k$. For $k\in \cK$, the convex hull of the longest path problem  is given as 
\begin{align}
\operatorname{conv}(\mathcal{X}^{(k)}) =
\left\{\bx:
\begin{aligned}
& \sum_{j:(i,j)\in \bE^k} x_{ij} - \sum_{j:(j,i)\in \bE^k} x_{ji}
= 
\begin{cases}
1, & i=s,\\
-1, & i=t,\\
0, & i\in \bV^k,
\end{cases}\\
& x_{ij} \ge 0,\quad \forall (i,j)\in \bE^k
\end{aligned}
\right\}.
\end{align}
% \begin{align}
%   \text{conv}(\mathcal{X}^{(k)}) = \left\{   &  \sum_{j:(i,j)\in \bE^k} x_{ij} - \sum_{j:(j,i)\in \bE^k} x_{ji} = \begin{cases}
%          1,\quad &i=s, \\
%          -1,\quad &i=t, \\
%          0,\quad &i\in \bV^k,
%      \end{cases} \\
% %& 0 \, \leq\,  \lambda_j^k- \lambda_i^k \leq\,  1, \quad \forall \, k,l \in \cK, \  \forall (i,j) \in \bE^k \cap \bE^l ,\\
% & x_{ij}^{k} \geq 0, \  \forall \,  k \in \cK,\  \forall (i,j) \in \bE^k, \right\}
% \end{align}
Given observed flow data $\bp^\cK = (p_{ij}^k: (i,j) \in \bE^k, k \in \cK)$, by Theorem \ref{thm:mdm-general-feascon}, the s-RAM characterization for longest path problems is given below.  

%\textbf{s-RAM characterization for Longest Path Problems.} 
Under Assumption \ref{asp:convex_perturbation} and \ref{asp:convexhull}, an observed flow collection $\bp^\cK$ is representable by an s-RAM  if and only if there exists $(\lambda_i^k \in \mathbb{R}: i\in s\cup t\cup \bV^k, k \in \mathcal{K})$,
    such that for any two distinct graph $k,l \in \cK$ containing a common edge $(i,j) \in \bE^k \cap \bE^l,$ 
\begin{align}
\begin{aligned}
     &\lambda_j^k - \lambda_i^k  > \lambda_j^l - \lambda_i^l  \quad \text{if} \quad p_{ij}^k > p_{ij}^l,\\
     &\lambda_j^k - \lambda_i^k  = \lambda_j^l - \lambda_i^l  \quad \text{if} \quad 0 < p_{ij}^k = p_{ij}^l <1.   \label{eq:mdm-feascon-longestpath}
\end{aligned}
\end{align}
As a result, checking whether the given flow data $\bp^\cK$ satisfies the s-RAM hypothesis can be accomplished by solving a linear program
with $\mathcal{O}(\vert \bV \vert K )$ continuous variables and $\mathcal{O}(\vert \bE \vert K) $ constraints, which is equivalent to checking whether the optimal objective value of \eqref{consistencyLP_longestpath} below is strictly positive:
\begin{align}
\begin{aligned}\label{consistencyLP_longestpath}
    \max_{\epsilon,\blam}\quad & \epsilon\\
  \text{s.t.}\quad &\lambda_j^k - \lambda_i^k  \geq \lambda_j^l - \lambda_i^l + \epsilon  \quad \text{if} \quad p_{ij}^k > p_{ij}^l ,\quad \forall k,l \in \cK, \forall (i,j) \in \bE^{k}\cap \bE^{l} ,\\
     & \lambda_j^k - \lambda_i^k  = \lambda_j^l - \lambda_i^l  \quad \text{if} \quad 0 < p_{ij}^k = p_{ij}^l < 1 , \quad \forall k,l \in \cK, \forall (i,j) \in \bE^{k}\cap \bE^{l}.  
\end{aligned}
\end{align}
\end{example}

\section{Step 2: Nonparametric Predictions with A New Polytope}\label{sec:prediction}

%\textcolor{black}{Parametric methods such as maximum likelihood estimation (MLE) provide a natural starting point for counterfactual prediction by first estimating parameters from data and then incorporating the estimators into a specified model for prediction. However, these methods typically require strong parametric assumptions and suffer from model misspecification risk. Alternatively, our framework offers a nonparametric, data-driven approach that directly learns from the observed choice data in an end-to-end manner without relying on parametric specifications.}

When the data $\cP^\cK$  is s-RAM-representable, we proceed to Step 2 with the workflow in Figure \ref{fig:flowchat}.
As an application of the exact characterization derived in Theorem \ref{thm:mdm-general-feascon}, we develop a nonparametric data-driven approach for making predictions for a new structured polytope 
$\text{conv}(\mathcal{X}^{\text{new}})$ where $\mathcal{X}^{\text{new}} := \{\bx \in \{0,1\}^n \mid \boldsymbol{A}\bx = \boldsymbol{b}, x_j = 0, \forall j \notin \cI^\text{new} \}$ and $ \cI^{\text{new}} \notin \cI^\cK.$ 
%$\mathcal{X}^{\text{new}} := \{\bx \in [0,1]^n \mid \boldsymbol{A}\bx = \boldsymbol{b}, x_j = 0, \forall j \notin \cI^\text{new} \}$ and $\mathcal{X}^{\text{new}} = \text{conv}( \mathcal{X}^{\text{new}} \cap \{0,1\}^n)$ with no prior data, i.e., $\cI^{\text{new}} \neq \cI^{(k)}$ for all $k\in \cK$. 

The key idea of the nonparametric approach is as follows. 
To predict the value of a function on $\text{conv}(\mathcal{X}^{\text{new}})$, we consider all s-RAMs that are consistent with the observed data, and use the worst-case (or best-case) value of this function across the consistent models as our estimate.
Thus, robust optimization forms the basis of our approach for allowing data to select a suitable s-RAM based on the prediction task at hand. 

\subsection{A Robust Optimization Formulation for Prediction} \label{sec:robust_prediction}

The collection of all s-RAM choice probability vectors $\bx = (x_{j}: j \in \cI^{\text{new}} )$ for the new polytope $\text{conv}(\mathcal{X}^{\text{new}})$ being consistent with the observed data $\bp^\cK$ is given by 
\[\mathcal{U}_{\text{new}}
   := \left\{ \bx: (\bp^\cK, \bx) \in \ps-RAM( \cK^\prime) \right\},\] where $\cI^{\cK^\prime} = \cI^\cK \cup \{ \cI^\text{new} \}$. 
   {\color{black}Let $f(\bx)$ denote any function of interest for the prediction problem e.g., the expected revenue of a new assortment, the flow of a chosen link in a new network, total welfare of a new assignment, or other performance measures derived from the choice probabilities. Define the worst case estimate of $f(\bx)$ with $\bx \in \mathcal{U}_{\text{new}}$ as $\underline{f}_\text{new}$.  Due to the exact characterization of s-RAM in \eqref{eq:lifted-set}, we obtain the following formulation for the prediction problem with $\text{conv}(\mathcal{X}^{\text{new}})$.}
\begin{proposition}
    \label{prop:wc-obj-reform}
    Suppose that the given data $\cP^\cK$ is representable by s-RAM and $\cI^\text{new} \subseteq [n]$ and $f(\bx)$ is a continuous function of $\bx.$ Then the worst-case estimate of $f(\bx)$ with $\text{conv}(\mathcal{X}^{\text{new}})$, denoted as $\underline{f}_\text{new}$,  equals:

\begin{subequations}
\begin{align}
     \min_{\bx,\blam} \quad & f(\bx) \nonumber\\
     \text{s.t.}\quad &
    x_{j}  \  \geq \   p_{j}^k \quad \text{if} \quad  \sum_{i\in [m]}\lambda_i^\text{new} a_{ij}  \geq \sum_{i\in [m]}\lambda_i^k a_{ij} , \quad \forall k \in \cK, \forall j\in \cI^\text{new} \cap \cI^{(k)}, \label{model:mdm_prediction-a}\\     
    & x_{j}  \  \leq \   p_{j}^k \quad \text{if} \quad  \sum_{i\in [m]}\lambda_i^\text{new} a_{ij}  \leq \sum_{i\in [m]}\lambda_i^k a_{ij}, \quad  \forall k \in \cK, \forall j\in \cI^\text{new} \cap \cI^{(k)}, \label{model:mdm_prediction-b}\\  
    & \sum_{j\in \cI^\text{new}} a_{ij} x_j \leq b_i, \qquad\qquad\qquad\qquad\quad\quad\quad\;\, \forall i \in [m],\nonumber\\
    & 0 \leq x_{j} \leq 1, \qquad\qquad\qquad\qquad\qquad\qquad\quad\quad \forall j \in \cI^\text{new}, \nonumber  \\
    &\sum_{i\in [m]}\lambda_i^k a_{ij}   < \sum_{i\in [m]}\lambda_i^l a_{ij}  \quad \text{if} \quad p_{j}^k < p_{j}^l, \qquad\qquad \forall k ,l \in \cK, \forall j\in \cI^{(k)} \cap \cI^{(l)}, \nonumber\\
     &\sum_{i\in [m]}\lambda_i^k a_{ij}   = \sum_{i\in [m]}\lambda_i^l a_{ij}  \quad \text{if} \quad 0<p_{j}^k = p_{j}^l <1,\quad \forall k,l \in \cK, \forall j\in \cI^{(k)} \cap \cI^{(l)}.\nonumber
    \end{align}
\end{subequations}
\end{proposition}
%One may obtain worst-case predictions for $x_i,$ with $i\in \cI^\text{new}$, by letting $r_j = 1$ and $r_j = 0$ for $j \neq i$ in the objective in the formulation in Proposition \ref{prop:wc-obj-reform}. Similarly, replacing the minimization in this formulation with a maximization yields a best-case (optimistic) estimate of the objective value $\bar{r}_\text{new}.$ 
Replacing the minimization in this formulation with a maximization yields a best-case (optimistic) estimate of the objective value $\bar{f}_\text{new}.$ 
Observe when the data is available for a richer collection, it leads to a less-conservative estimate for $\underline{f}_\text{new}$ and a narrower interval $[\underline{f}_\text{new},\bar{f}_\text{new}]$ as plausible values for the objective value estimates which are consistent with the given data and the s-RAM hypothesis.  This is because the number of constraints in the constraint collections \eqref{model:mdm_prediction-a}-\eqref{model:mdm_prediction-b} is larger when the collection $\cK$ for which the data is available is made richer. %In other words, $\cI_{\cK_1} \subseteq \cI_{\cK_2}$ when $\cK_1 \subseteq \cK_2$ and therefore the resulting $\ps-RAM(\cK_2 \cup \{A\})$ is nested within $\ps-RAM(\cK_1 \cup \{A\}).$ 
\subsection{A Compact-sized Integer Formulation for Prediction}
\vspace{-5pt} 
As a generally applicable approach for solving the prediction problem in Proposition \ref{prop:wc-obj-reform}, one may model the ``if'' conditions in \eqref{model:mdm_prediction-a} - \eqref{model:mdm_prediction-b} via additional binary variables $(\delta_{\text{new},k}, \delta_{k,\text{new}}: k \in \cK)$ to obtain the mixed-integer linear reformulation with $\mathcal{O}(nK)$ binary variables, $\mathcal{O}(n+mK)$ continuous variables, and $\mathcal{O}(nK)$ constraints as follows.
\begin{proposition}
\label{prop:wc-obj-milp}
Suppose that the assumptions in Proposition \ref{prop:wc-obj-reform} are satisfied. Then 
for any $0 < \epsilon < 1 /(2nK)$, $\underline{f}_\text{new}$ in Proposition \ref{prop:wc-obj-reform} equals the value of the following integer program:   
\begin{subequations}
\begin{align}
\min_{\bx,\blambda,\bdelta,\boldsymbol{z}}  & \quad f(\bx) \nonumber\\
 \text{s.t.}\  -&\delta_j^{\text{new},k}  \leq     \sum_{i\in [m]}\lambda_i^\text{new} a_{ij}  -  \sum_{i\in [m]}\lambda_i^k a_{ij}  \leq   1 -  (1+\epsilon)\delta_j^{\text{new},k}, \, \forall k \in \cK, \forall j \in \cI^\text{new} \cap \cI^{(k)}, \label{mdm:milp-a}\\
 -&\delta_j^{k,\text{new}}   \leq    \sum_{i\in [m]}\lambda_i^k a_{ij}  -  \sum_{i\in [m]}\lambda_i^\text{new} a_{ij}  \leq   1 -  (1+\epsilon)\delta_j^{k,\text{new}}, \; \forall k \in \cK, \forall j\in \cI^\text{new} \cap \cI^{(k)}, \label{mdm:milp-b}\\
&  \delta_j^{\text{new},k} - 1 \ \leq \  p_{j}^k - x_{j}   \ \leq \ 1 - \delta_j^{k,\text{new}}, \qquad\qquad\qquad\quad\;\;   \forall k \in \cK, \forall j\in \cI^\text{new} \cap \cI^{(k)},\label{mdm:milp-c}\\
 -&(\delta_j^{\text{new},k} +  \delta_j^{k,\text{new}}) \ \leq \  p_{j}^k - x_{j} \ \leq \ \delta_j^{\text{new},k} + \delta_j^{k,\text{new}},\ \qquad\quad    \forall k \in \cK, \forall j\in \cI^\text{new} \cap \cI^{(k)},\label{mdm:milp-d}\\ 
& \sum_{i\in [m]}\lambda_i^k a_{ij} - \sum_{i\in [m]}\lambda_i^l a_{ij} \ \geq \ \epsilon, \qquad \forall k ,l \in \cK, \forall j\in \cI^{(k)} \cap \cI^{(l)} \ \textnormal{ s.t. }\  p_{j}^k > p_{j}^{l} \nonumber,\\
&\sum_{i\in [m]}\lambda_i^k a_{ij} - \sum_{i\in [m]}\lambda_i^l a_{ij} \ = \ 0, \qquad \forall k ,l \in \cK, \forall j\in \cI^{(k)} \cap \cI^{(l)} \ \textnormal{ s.t. }\ 0< p_{j}^k = p_{j}^{l} <1 \nonumber,\\
& \sum_{j\in \cI^\text{new}} a_{ij} x_j \leq b_i, \qquad\qquad\qquad\qquad\quad \forall i \in [m], \nonumber\\
& 0 \leq x_{j} \leq 1,  \qquad\qquad\qquad\qquad\quad\quad\quad\;\; \forall j \in \cI^\text{new},\nonumber\\ 
%-&1 \, \leq\, \sum_{i\in [m]}\lambda_i^\text{new} a_{ij} -  \sum_{i\in [m]}\lambda_i^k a_{ij} \leq\,  1, \qquad \forall k \in \cK, \forall j \in \cI^\text{new} \cap \cI^{(k)},\nonumber \\
&\delta_j^{\text{new},k}, \delta_j^{k,\text{new}} \, \in \, \{0,1\},  \qquad\qquad\qquad\quad\; \forall \, k \in \cK, \forall j \in \cI^\text{new} \cap \cI^{(k)}. \nonumber 
\end{align}
\end{subequations}
Likewise, the optimistic estimate for the objective value $\bar{f}_\text{new}$  equals the optimal value obtained by maximizing over the constraints in the above integer program.  
\end{proposition}

\section{Step $1^\prime$: Nonparametric Estimation to Find the Best Fit }\label{sec:estimation}

When the observed data does not satisfy the s-RAM hypothesis perfectly. We turn to Step $1^\prime$ in Figure \ref{fig:flowchat}. Considering data instances that are not s-RAM-representable, we next seek to quantify the limit or the cost of approximating given data with s-RAM and a procedure for identifying  s-RAM-representable probabilities offering the best fit. Then, we use the best-fitting approximated s-RAM probabilities for the downstream prediction tasks. 
\vspace{-5pt} 
\subsection{ A Nonparametric Estimation Formulation for s-RAM}
\vspace{-5pt} 
Given observed data $\bp^\cK$, suppose that a loss function $\bx^\cK \mapsto  \text{loss}(\bp^\cK,\bx^\cK)$ measures the degree of inconsistency in approximating the data $\bp^\cK$ with an s-RAM-consistent probability assignment $\bx^\cK \in \ps-RAM(\cK).$ 
We assume the loss function is non-negative and convex, and satisfies the property that 
$\text{loss}(\bp^\cK,\bx^\cK )=0$ if and only if $\bx^\cK = \bp^\cK$. Suppose that $(n_k: k \in \cK)$ is a vector of non-negative weights over polytopes in the historical dataset. Then common choices include a norm-based loss such as $\sum_{k \in \cK} n_k \Vert \bp_k - \bx_k \Vert$ or a Kullback-Liebler divergence based loss such as $-\sum_{k \in \cK} n_k \sum_{i \in \cI^{(k)}} p_i^k \log(x_i^k/p_i^k)$. 

We define the limit of the s-RAM , denoted $\mathcal{L}(\bp^\cK)$, is the minimum possible loss achievable by approximating $\bp^\cK$ using a model consistent with s-RAM: %, denoted by $\mathcal{L}(\bp^\cK),$ as the smallest value of $\text{loss}(\bp^\cK,\bx^\cK )$ attainable by fitting the observed data $\bp_{\cK}$ with an s-RAM choice model:
\begin{align}
    \mathcal{L}(\bp^\cK) = \inf \left\{ \text{loss}(\bp^\cK,\bx^\cK)\,:\, \bx_{\cK} \in \ps-RAM(\cK) \right\}.
    \label{eq:limit-of-mdm}
\end{align}
% From the definition above,  evaluating the limit  $\mathcal{L}(\bp^\cK)$ can be viewed as identifying a probability assignment $\bx^\cK_\ast$ which is consistent with the s-RAM hypothesis and is about as close to any s-RAM can be to the observed data $\bp^\cK.$ Thus, 
From the above definition, evaluating the limit $\mathcal{L}(\bp^\cK)$ amounts to finding a probability assignment $\bx^\cK_\ast$ in the s-RAM class that best fits the observed data $\bp^\cK$. %Specifically, any minimizer $\bx^\cK_\ast$ of \eqref{eq:limit-of-mdm} best fits the observed data within the s-RAM family. 
Specifically, suppose we take $\text{loss}(\bp^\cK,\bx^\cK ) = - \sum_{k \in \cK} n_k \sum_{j \in \cI^{(k)}} p_j^k \log(x_j^k/p_j^k)$ with the weight $n_k$ for $k \in \cK$ to be equal to the number of observations of a polytope $\text{conv}(\mathcal{X}^{(k)})$ in the dataset. A solution to \eqref{eq:limit-of-mdm} can be viewed as obtained from maximum likelihood estimation without any parametric assumption on s-RAM.

Due to the characterization of s-RAM in Theorem \ref{thm:mdm-general-feascon}, we have the equivalent formulation for $\mathcal{L}(\bp^\cK)$ below. 
\begin{proposition}\label{prop:limit-mdm}
Under Assumption \ref{asp:convex_perturbation},
 the limit $\mathcal{L}(\bp^\cK)$ equals 
\begin{align}
     \min_{\bx^{\cK},\blambda }\quad &  \sum_{k \in \cK}  \textnormal{loss}(\bp^k,\bx^k) \nonumber  \\ 
      \text{s.t.}\quad 
      %& x_{j}^k  \  \leq \   x_{j}^l \quad \text{if} \quad  \sum_{i\in [m]}\lambda_i^l a_{ij}  \leq \sum_{i\in [m]}\lambda_i^k a_{ij} , \qquad \forall k,l \in \cK, \forall j\in \cI^{(k)} \cap \cI^{(l)}, \\     
    & x_{j}^k  \  \leq \   x_{j}^l \quad \text{if} \quad  \sum_{i\in [m]}\lambda_i^k a_{ij}  \leq \sum_{i\in [m]}\lambda_i^l a_{ij} , \quad\;\; \forall k,l \in \cK, \forall j\in \cI^{(k)} \cap \cI^{(l)},\nonumber  \\ 
    & \sum_{j\in \cI^{(k)}}  a_{ij} x_j^k \leq b_i, \qquad\qquad\qquad\qquad\qquad\quad \forall k \in \cK, \forall i \in [m], \label{model:lom} \\  
    & 0 \leq x_j^k \leq 1, \qquad\qquad\qquad\qquad\qquad\qquad\quad\;\; \forall k \in \cK, \forall j \in \cI^{(k)}. \nonumber 
\end{align}
\end{proposition}

\subsection{A Compact-sized Mixed Integer Convex Reformulation for Estimation of s-RAM}

The Proposition below provides a generally applicable mixed-integer convex reformulation for \eqref{model:lom}. 
\begin{proposition} 
\label{prop:micp} 
Suppose that Assumption \ref{asp:convex_perturbation} holds. Then for any $0 < \epsilon < 1 /(2nK)$, 
the limit  $\mathcal{L}(\bp^\cK)$ equals the value of the following mixed integer convex program: 
\begin{align}
 \min_{\bx,\blambda,\bdelta}\quad &  \sum_{k \in \cK}  \textnormal{loss}(\bp^k,\bx^k) \nonumber \\
 \text{s.t.}\quad    
 -&\delta_j^{k,l} \ \leq\    \sum_{i\in [m]}\lambda_i^k a_{ij}  -  \sum_{i\in [m]}\lambda_i^l a_{ij} \ \leq \  1 -  (1+\epsilon)\delta_j^{k,l}, \quad \forall k,l \in \cK, \forall j\in \cI^{(l)} \cap \cI^{(k)}, \nonumber \\
&  \delta_j^{k,l} - 1 \ \leq \ x_j^l - x_j^k \ \leq \ 1 - \delta_j^{l,k}, \qquad\qquad\qquad\qquad\;\;\,  \forall k,l \in \cK, \forall j\in \cI^{(l)} \cap \cI^{(k)}, \nonumber \\
 -&(\delta_j^{k,l} +  \delta_j^{l,k}) \ \leq \ x_j^l - x_j^k \ \leq \ \delta_j^{k,l} + \delta_j^{l,k},\ \qquad\qquad\qquad\,   \forall k,l \in \cK, \forall j\in \cI^{(l)} \cap \cI^{(k)}, \label{model:mdm_micp} \\ 
& \sum_{j\in \cI^{(k)}} a_{ij} x_j \leq b_i, \qquad\qquad\qquad\qquad\qquad\qquad\qquad\qquad\!\! \forall k  \in  \cK, \, \forall i  \in  [m], \nonumber \\
& 0 \leq x_{j}^k \leq 1, \    \qquad\qquad\qquad\qquad\qquad\qquad\qquad\qquad\qquad\! \forall k\in \cK, \forall j \in \cI^{(k)}, \nonumber \\ 
& \delta_j^{k,l} \ \in \ \{0,1\},  \qquad\qquad\qquad\qquad\qquad\qquad\qquad\qquad\quad\; \forall \, k,l \in \cK,  \forall j\in \cI^{(l)} \cap \cI^{(k)} . \nonumber      
\end{align} 
\end{proposition}
Suppose that  $(\bx^\cK_\ast, \blambda_\ast)$ attains the minimum in \eqref{model:lom} or equivalently in  \eqref{model:mdm_micp}. Since  $\ps-RAM(\cK)$ is not a closed set and the constraints in \eqref{model:lom} allow ${x_j^k}_\ast = {x_j^l}_\ast $ even when $\sum_{i\in [m]} {\lambda_i^k}_\ast a_{ij} \neq \sum_{i\in [m]} {\lambda_i^l}_\ast a_{ij},$ the solution $\bx^\cK_\ast$ can only be guaranteed to be arbitrarily close to the s-RAM-representable collection $\ps-RAM(\cK).$ Then, $\bx^\cK_\ast$ can not be directly used for prediction since it does not satisfy the last two conditions of the formulation in   Proposition \ref{prop:wc-obj-reform} or equivalently the corresponding constraints in Proposition \ref{prop:wc-obj-milp}.
Therefore, to proceed with the prediction task, if one wishes to obtain a $\delta$-optimal s-RAM-representable choice probability assignment, for some $\delta > 0$, they may do so as follows: Equipped with the optimal value $\mathcal{L}(\bp_\cK) = \text{loss}(\bp^\cK,\bx^\cK_\ast)$ and an  optimal $\blambda_\ast,$ a $\delta$-optimal s-RAM-representable probability assignment $\bx_\cK$ can be obtained  by solving the following convex program: 
\begin{align} 
&\!\!\!\!\!\!\!\!\!\!\!\!\!\!\!\!\!\!\!\max_{\bx^\cK,\epsilon}  \quad  \epsilon  \nonumber \\   
 \text{s.t.}\quad &
\textnormal{loss}(\bp^\cK, \bx^\cK) \leq \mathcal{L}(\bp^\cK)(1 + \delta), \nonumber \\
& x_{j}^k  \  \geq \   x_{j}^l+\epsilon \quad \text{if} \quad  \sum_{i\in [m]} {\lambda_i^k}_\ast a_{ij}  > \sum_{i\in [m]}{\lambda_i^l}_\ast a_{ij} , \quad \forall k,l \in \cK, \forall j\in \cI^{(k)} \cap \cI^{(l)}, \nonumber\\
& x_{j}^k  \  = \   x_{j}^l \quad \text{if} \quad  \sum_{i\in [m]} {\lambda_i^k}_\ast a_{ij}  = \sum_{i\in [m]}{\lambda_i^l}_\ast a_{ij} , \qquad\quad\!\! \forall k,l \in \cK, \forall j\in \cI^{(k)} \cap \cI^{(l)}, \label{mdm:perturb}\\
 & \sum_{j\in \cI^{(k)}}  a_{ij} x_j^k \leq b_i, \qquad\qquad\qquad\qquad\qquad\quad \forall k \in \cK, \forall i \in [m],\nonumber \\
 & 0 \leq x_j^k \leq 1, \qquad\qquad\qquad\qquad\qquad\qquad\qquad\! \forall k \in \cK, \forall j \in \cI^{(k)}.  \nonumber 
\end{align}
Given the constraints in \eqref{mdm:perturb} and the characterization in Theorem \ref{thm:mdm-general-feascon}, any choice probability collection $\bx^\cK$ obtained by solving \eqref{mdm:perturb} is s-RAM-representable. Moreover, it cannot be improved to achieve a fit better than within a magnitude of $\delta$ for any arbitrary $\delta > 0$, due to the constraint $\text{loss}(\bp^{\cK}, \bx^{\cK}) \leq \mathcal{L}(\bp^\cK)(1 + \delta)$. 
Note that when $\text{loss}(\cdot,\cdot)$ is defined in terms of the $L_1$-norm, the formulation \eqref{mdm:perturb} is a linear program with $\mathcal{O}(nK)$ continuous variables and $\mathcal{O}(nK^2+nm)$ constraints and \eqref{model:mdm_micp} is a mixed-integer linear program  with $\mathcal{O}(nK^2)$ binary variables, $\mathcal{O}(nK)$ continuous variables and $\mathcal{O}(nK^2+nm)$ constraints.

\section{Numerical Experiments}\label{sec:experiment}

In this section, we generate synthetic datasets to evaluate (i) the explanatory ability of s-RAM as captured by the representability in Step 1 and the best-fitting cost of s-RAM data in Step $1^\prime$, and (ii) the efficacy of predictions obtained in Step 2 of the workflow in Figure \ref{fig:flowchat} for a wide range of combinatorial problems including longest path, shortest path, linear assignment, and constrained subset selection problems. Missing details on experiment implementation in this section are provided in Section \ref{app:exp_examples}. We use a MacBook Pro with a 12-core Apple M4 Pro chip and 24 GB RAM for all experiments. 

\subsection{Longest Path Problem} \label{sec:experiment_longestpath}
\subsubsection{Explanatory ability of s-RAM in the longest path problem} We generate data from models that lie outside the s-RAM framework to assess the explanatory ability of s-RAM under model misspecification.%\par\noindent

\textit{Generation of historical $\mathcal{X}^{(k)}$ instances.} We take the number of historically observed graph instances to vary across the range \( K \in \{5, 10, 15, 20\} \). We generate directed acyclic graphs (DAGs) with 10 nodes. Edges are included with probability drawn uniformly from \([0.4, 0.6]\) while ensuring source-to-sink connectivity. On average, there are around 25 edges in the DAGs.

\textit{Generation of historical $\bp^{(k)}$ instances.}
We generate the respective synthetic $\bp^{(k)}$ data values outside the s-RAM framework by assigning edge costs in the form \(v_e + \epsilon_e\), where \(v_e\) is a deterministic component drawn uniformly from \([1, 7]\), and \(\epsilon_e\) is a Gaussian noise term. We consider three noise structures: independent, negatively correlated (correlation coefficients drawn from \([-0.9, -0.5]\)), and positively correlated (from \([0.5, 0.9]\)), with variances sampled uniformly from \([9, 25]\). For each graph, we solve the deterministic version of the problem 100 times under different noise realizations and average the solutions to construct the observed data.

% \textit{Representability test and estimation with $\bp^\cK$.} For each instance, we verify the representability of s-RAM for it via \eqref{eq:mdm-feascon-order} and, if violated, we compute the average absolute deviation loss. In this case, \eqref{model:mdm_micp} for the longest path problem with the chosen loss function can be written as 
% \textbf{Estimation for Longest Path Problems.}
% In the case of longest path problems, by taking $0 < \epsilon < 1 /(2|\bV|K),$ formulation \eqref{model:mdm_micp} can be written as: 

\textit{Representability test and estimation with $\bp^\cK$.} 
For each instance, we first verify whether the observed data are s-RAM representable by solving \eqref{consistencyLP_longestpath}. 
If an instance is not s-RAM representable, we quantify the degree of violation by computing the average absolute deviation loss. For the longest path problem, the resulting estimation problem corresponding to \eqref{model:mdm_micp} with the chosen loss function can be written as
\begin{align}
\begin{aligned}
    \min_{\bx,\blambda,\bdelta\!\!\,}\quad&  \quad   \frac{1}{K}\sum_{k=1}^K\sum_{(i,j)\in \bE^k}p_{ij}^k |x_{ij}^k - p_{ij}^k| \\
    \text{s.t.}\quad \ -&\delta^{k,l}_{i,j} \ \leq\    \lambda_j^k- \lambda_i^k  -  \lambda_j^l + \lambda_i^l \ \leq \  1 -  (1+\epsilon)\delta^{k,l}_{i,j}, \quad\  \ \forall \,   k,l \in \cK,\ \forall (i,j) \in \bE^{k}\cap \bE^{l}, \\
&  \delta^{k,l}_{i,j} - 1 \ \leq \  x_{ij}^{l}  - x_{ij}^{k}  \ \leq \ 1 - \delta^{l,k}_{i,j}, \quad\qquad\quad   \forall \,   k,l \in \cK,\ \forall (i,j) \in \bE^{k}\cap \bE^{l},\\
 -&(\delta^{k,l}_{i,j} +  \delta^{l,k}_{i,j} ) \ \leq \ x_{ij}^{l}  - x_{ij}^{k}   \ \leq \ \delta^{k,l}_{i,j} + \delta^{l,k}_{i,j} ,\ \quad     \forall \,   k,l \in \cK,\ \forall (i,j) \in \bE^{k}\cap \bE^{l},\\ 
&  \sum_{j:(i,j)\in \bE^k} x_{ij} - \sum_{j:(j,i)\in \bE^k} x_{ji} = \begin{cases}
         1,\quad &i=s, \\
         -1,\quad &i=t, \\
         0,\quad &i\in \bV^k,
     \end{cases} \forall k \in \cK,\\
%& 0 \, \leq\,  \lambda_j^k- \lambda_i^k \leq\,  1, \quad \forall \, k,l \in \cK, \  \forall (i,j) \in \bE^k \cap \bE^l ,\\
& x_{ij}^{k} \geq 0, \  \forall \,  k \in \cK,\  \forall (i,j) \in \bE^k,\\
&\delta^{k,l}_{i,j} \, \in \, \{0,1\}, \  \forall \,  k,l \in \cK,\  \forall (i,j) \in \bE^k \cap \bE^l, 
\label{model:lp_micp}
\end{aligned}
\end{align}
where $0 < \epsilon < 1 /(2|\bV|K)$ and $\bV$ is the collection of all possible nodes. For each value of $K$ and each correlation structure, we report in Table~\ref{tab:explanatory_longestpath} the fraction of s-RAM-representable instances and the corresponding average absolute deviation loss, computed over 20 independently generated instances. 
Because the data are generated randomly, enforcing exact representability can be overly restrictive: instances with vanishingly small loss may be classified as non-representable due to numerical or sampling noise. 
Accordingly, we treat an instance as s-RAM representable if its average absolute deviation loss is below $10^{-6}$, and report the representable fraction under this criterion.

% where $0 < \epsilon < 1 /(2|\bV|K).$
% %$$ \min_{\bx} \frac{1}{K}\sum_{k=1}^K\sum_{e\in \bE^k}p_{e}^k |x_{e}^k - p_{e}^k|$$using \eqref{model:mdm_micp}.
% Then for each \(K\) and correlation structure, we report the representable proportion and the average absolute deviation loss over 20 instances in Table \ref{tab:explanatory_longestpath} below. Due the randomness of data generation, if we require exact representability of the instance then instances with extremely small loss (below $10^{-6}$) can be viewed as nonrepresentable instances (which might be not a good presentation of the representable results), here we report the fraction of s-RAM representable instances as fraction of instances with loss below $10^{-6}$.

\textit{Representation and estimation results in the longest path problem.} 
Table~\ref{tab:explanatory_longestpath} shows that s-RAM performs well across all tested scenarios despite model misspecification, achieving high representability for $K \leq 15$. 
For instance, when $K=5$, more than $85\%$ of the instances are s-RAM representable across all correlation structures, and the representable fraction remains above $60\%$ even when $K=15$. 
Moreover, the associated estimation error remains small: the average absolute deviation loss is consistently on the order of $10^{-3}$ or smaller, even in settings where exact representability is violated. 
Naturally, the representability decreases with $K$ due to the increasing number of feasibility constraints in \eqref{eq:mdm-feascon-longestpath} induced by additional historical observations. 
Overall, we observe the framework to be robust across different correlation structures.

%\textit{Representation and estimation results in the longest path problem.} We observe from Table \ref{tab:explanatory_longestpath} that  s-RAM performs well across all tested scenarios despite model misspecification, achieving high representability for \(K \leq 15\) \textcolor{red}{(request: explain a bit like for example K=5 ... )} and low estimation error \textcolor{red}{(request: explain a bit like even at the scale of 10-3... )}  . Naturally, the representability decreases with $K$ due to increasing number of constraints in \eqref{eq:mdm-feascon-longestpath} introduced by additional data availability. We observe the  framework to be robust across correlation structures.%, with better performance under independent and negatively correlated noise.

\begin{table}[htb]
\centering
\caption{Representation and estimation in the longest path problem with Gaussian instances}
\label{tab:explanatory_longestpath}
\scalebox{0.85}{
\begin{tabular}{|c|cc|cc|cc|}
\hline
\multirow{2}{*}{K} 
& \multicolumn{2}{c|}{independent gaussian} 
& \multicolumn{2}{c|}{gaussian with -ve correlation} 
& \multicolumn{2}{c|}{gaussian with +ve correlation} \\ \cline{2-7} 
\multicolumn{1}{|c|}{}                   & \multicolumn{1}{c|}{\begin{tabular}[c]{@{}c@{}}fraction \\ of s-RAM\\ representable\\ instances\end{tabular}} & \multicolumn{1}{c|}{\begin{tabular}[c]{@{}c@{}}average absolute\\ deviation loss (se)\\ $(10^{-3})$\end{tabular}} & \multicolumn{1}{c|}{\begin{tabular}[c]{@{}c@{}}fraction \\ of s-RAM\\ representable\\ instances\end{tabular}} & \multicolumn{1}{c|}{\begin{tabular}[c]{@{}c@{}}average absolute\\ deviation loss (se)\\ $(10^{-3})$ \end{tabular}} & \multicolumn{1}{c|}{\begin{tabular}[c]{@{}c@{}}fraction \\ of s-RAM\\ representable\\ instances\end{tabular}} & \multicolumn{1}{c|}{\begin{tabular}[c]{@{}c@{}}average absolute\\ deviation loss (se)\\ $(10^{-3})$ \end{tabular}} \\ \hline
5  & \multicolumn{1}{c|}{0.95} & 0.02 (0.02) 
   & \multicolumn{1}{c|}{0.85} & 1.35 (0.92) 
   & \multicolumn{1}{c|}{1.00} & 0.00 (0.00) \\ \hline
10 & \multicolumn{1}{c|}{0.85} & 0.78 (0.57) 
   & \multicolumn{1}{c|}{0.70} & 0.62 (0.23) 
   & \multicolumn{1}{c|}{0.45} & 1.31 (0.58) \\ \hline
15 & \multicolumn{1}{c|}{0.70} & 0.65 (0.28) 
   & \multicolumn{1}{c|}{0.60} & 1.52 (0.48) 
   & \multicolumn{1}{c|}{0.55} & 0.88 (0.37) \\ \hline
20 & \multicolumn{1}{c|}{0.45} & 1.92 (0.72) 
   & \multicolumn{1}{c|}{0.45} & 1.36 (0.56) 
   & \multicolumn{1}{c|}{0.30} & 1.74 (0.52) \\ \hline
\end{tabular}}
\end{table}

\subsubsection{Predictive ability of s-RAM in the longest path problem}
We evaluate the predictive accuracy of nonparametric s-RAM on counterfactual tasks using data generated from parametric s-RAMs. For the longest path, we test whether s-RAM can predict the probability of an edge being included in the critical path when the graph changes to a new structure.  

\textit{Generation of training instances $\bp^{(k)}$ and ground truth for the longest path problem.}
Since s-RAM and MDM are equivalent, to retain interpretability in utility theory, we generate data from parametric MDMs. In the longest path, marginals of noise follow nonidentical exponential distributions with rates sampled from $[1,2]$. To ensure that the edge of interest is influential in the optimal solution, we set its deterministic cost to $8$. All other edge costs are sampled independently and uniformly from $[1,7]$. We generate directed acyclic graphs following the same procedure as in the explanatory ability test above. Then, for each instance $k$, we use formulation \eqref{data:longestpath} to generate training instances $\bp^{(k)}.$ We evaluate the prediction performance of the nonparametric s-RAM across collection sizes $K \in \{5,10,20\}$. For each value of $K$, we generate $20$ training instances. Given a new graph $\bG^{\text{new}} = (s \cup t \cup \bV^{\text{new}}, \bE^{\text{new}})$. Let $\bx^{\text{true}}$ denote the resulting vector of edge inclusion probabilities for $\bG^{\text{new}}.$ Then we use the same parameter values as in the training instances together formulation \eqref{data:longestpath} to generate $\bx^{\text{true}}$. 
Let $(i^*,j^*)$ be the edge of interest. The corresponding ground-truth quantity is then $x^{\text{true}}_{i^* j^*}$. 

\textit{Prediction with nonparametric s-RAM in the longest path.} 
For the longest path problem, we compute the s-RAM worst- and best-case estimates for the probability that the selected edge $(i^*,j^*)$ is included in the critical path.  Since the observed data $\bp^\cK$ are representable by s-RAM, given a new graph $\bG^{\text{new}} = (s \cup t \cup \bV^{\text{new}}, \bE^{\text{new}})$, Proposition \ref{prop:wc-obj-milp} implies that the worst-case estimate $\underline{f}_{\text{new}}$ can be obtained by solving the following mixed-integer linear program, with $0 < \epsilon < 1/(2|\bV|K)$:
%For longest path problem, we compute the s-RAM worst- and best-case estimates for the probability being included in the critical path for the selected edge. Let the selected edge denote as $(i^*,j^*).$ Note that $\bp^\cK$ is representable by s-RAM, given a new graph $\bG^{\text{new}} = (s\cup t \cup \bV^{\text{new}}, \bE^{\text{new}})$, by Proposition \ref{prop:wc-obj-milp}, the prediction formulation to compute the worst case estimate $\underline{f}_\text{new}$ by taking $0 < \epsilon < 1 /(2|\bV|K),$ is as follows.   
\begin{align}
\min_{\bx,\blambda,\bdelta\!\!\,}\quad&  \quad   x_{i^* j^*} \label{model:lp_milp} \\
    \text{s.t.}\quad \ -&\delta^{\text{new},k}_{i,j} \ \leq\    \lambda_j^{\text{new}} - \lambda_i^{\text{new}} -  \lambda_j^k+ \lambda_i^k \ \leq \  1 -  (1+\epsilon)\delta^{\text{new},k}_{i,j}, \   \forall \,  (i,j) \in \bE^\text{new}, \forall k\in \cK, (i,j) \in \bE^k, \nonumber \\
 -&\delta^{k, \text{new}}_{i,j}  \ \leq\    \lambda_j^k- \lambda_i^k  -  \lambda_j^{\text{new}} + \lambda_i^{\text{new}}\ \leq \  1 -  (1+\epsilon)\delta^{k, \text{new}}_{i,j}, \ \forall \,  (i,j) \in \bE^\text{new}, \forall k\in \cK, (i,j) \in \bE^k, \nonumber \\
&  \delta^{\text{new},k}_{i,j} - 1 \ \leq \  p_{ij}^k - x_{ij}   \ \leq \ 1 - \delta^{k, \text{new}}_{i,j}, \quad\qquad\quad\quad\;\;   \forall \,  (i,j) \in \bE^\text{new}, \, \forall k\in \cK, (i,j) \in \bE^k,\nonumber \\
 -&(\delta^{\text{new},k}_{i,j} +  \delta^{k, \text{new}}_{i,j} ) \ \leq \ p_{ij}^k - x_{ij} \ \leq \ \delta^{\text{new},k}_{i,j} + \delta^{k, \text{new}}_{i,j} ,\ \quad     \forall \,  (i,j) \in \bE^\text{new}, \, \forall k\in \cK, (i,j) \in \bE^k,\nonumber \\ 
& \lambda_j^k- \lambda_i^k - \lambda_j^l + \lambda_i^l \ \geq \ \epsilon, \qquad\qquad\! \forall \, k,l \in \cK, \  \forall (i,j) \in \bE^k \cap \bE^l \ \textnormal{ s.t. }\  p_{ij}^k > p_{ij}^l,\nonumber \\
&\lambda_j^k- \lambda_i^k - \lambda_j^l + \lambda_i^l \ = \ 0, \qquad\quad\ \ \, \forall \, k,l \in \cK, \  \forall (i,j) \in \bE^k \cap \bE^l \ \textnormal{ s.t. }\  0< p_{ij}^k = p_{ij}^l <1 ,\nonumber \\
&  \sum_{j:(i,j)\in \bE^\text{new}} x_{ij} - \sum_{j:(j,i)\in \bE^\text{new}} x_{ji} = \begin{cases}
         1,\quad &i=s, \\
         -1,\quad &i=t, \\
         0,\quad &i\in \bV^\text{new},
     \end{cases}\nonumber \\
& x_{ij}\geq 0, \  \forall k\in \cK, (i,j) \in \bE^k,\nonumber \\
&\delta^{\text{new},k}_{i,j}, \delta^{k, \text{new}}_{i,j} \, \in \, \{0,1\}, \  \forall \,  (i,j) \in \bE^\text{new}, \, \forall k\in \cK, (i,j) \in \bE^k. \nonumber
\end{align}
%Likewise, the best-case estimate $\bar{f}_{\text{new}}$ is obtained by maximizing the same objective subject to the identical set of constraints.
Likewise, the best-case estimate for the objective value $\bar{f}_\text{new} $ equals the optimal value obtained by maximizing the same objective over the constraints in the above mixed integer linear program. % In our experiment, we take $f(\bx) = x_{\hat{e}}$ where $\hat{e}$ is the pre-specified edge. 

\begin{figure}[htb]
\centering
\begin{minipage}[t]{0.3\textwidth}
\centerline{longest path with $K = 5$ }
\centerline{\includegraphics[scale=0.071]{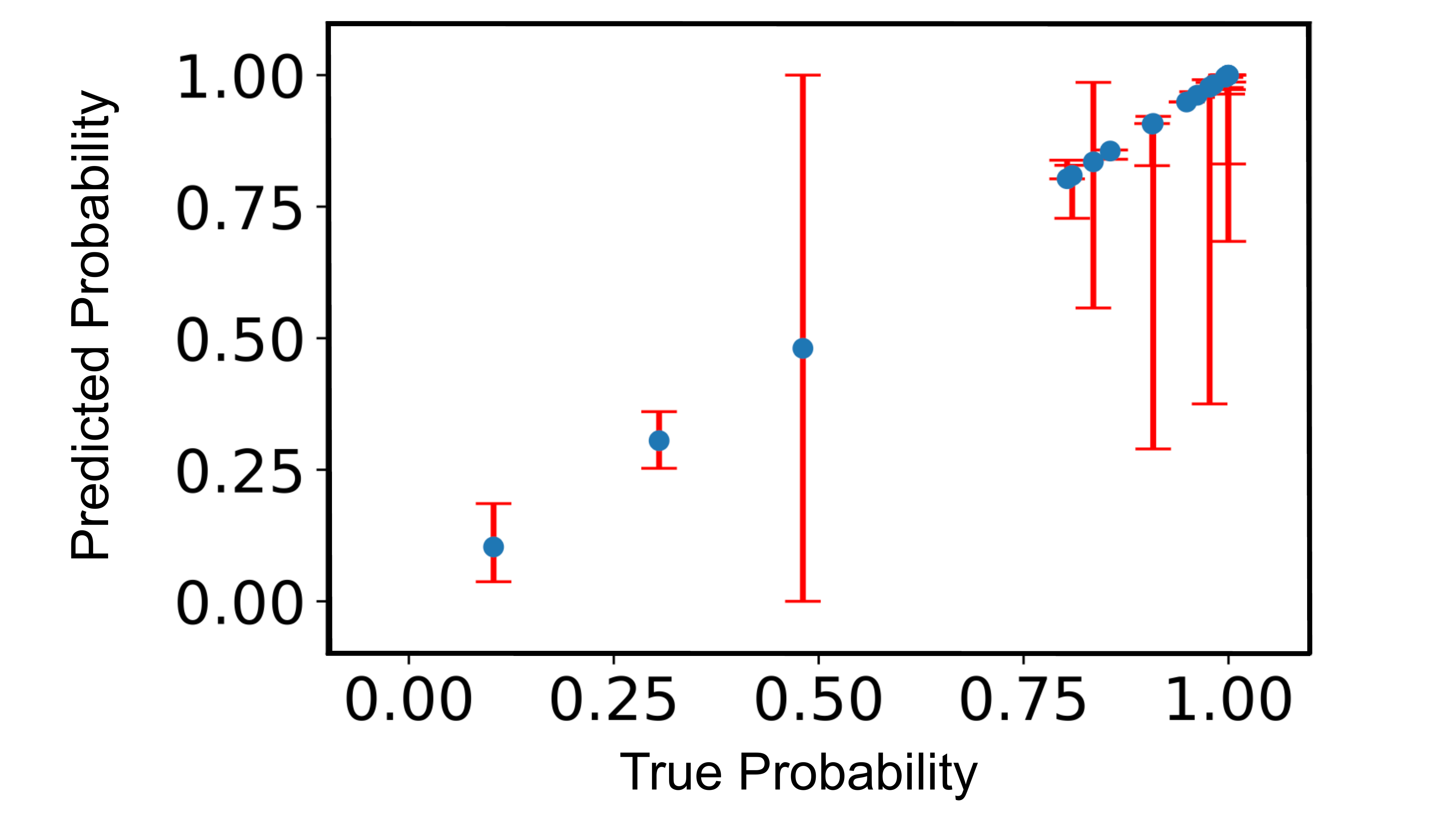}}
\end{minipage}
\;
\begin{minipage}[t]{0.3\textwidth}
\centerline{longest path with $K = 10$}
\centerline{\includegraphics[scale=0.071]{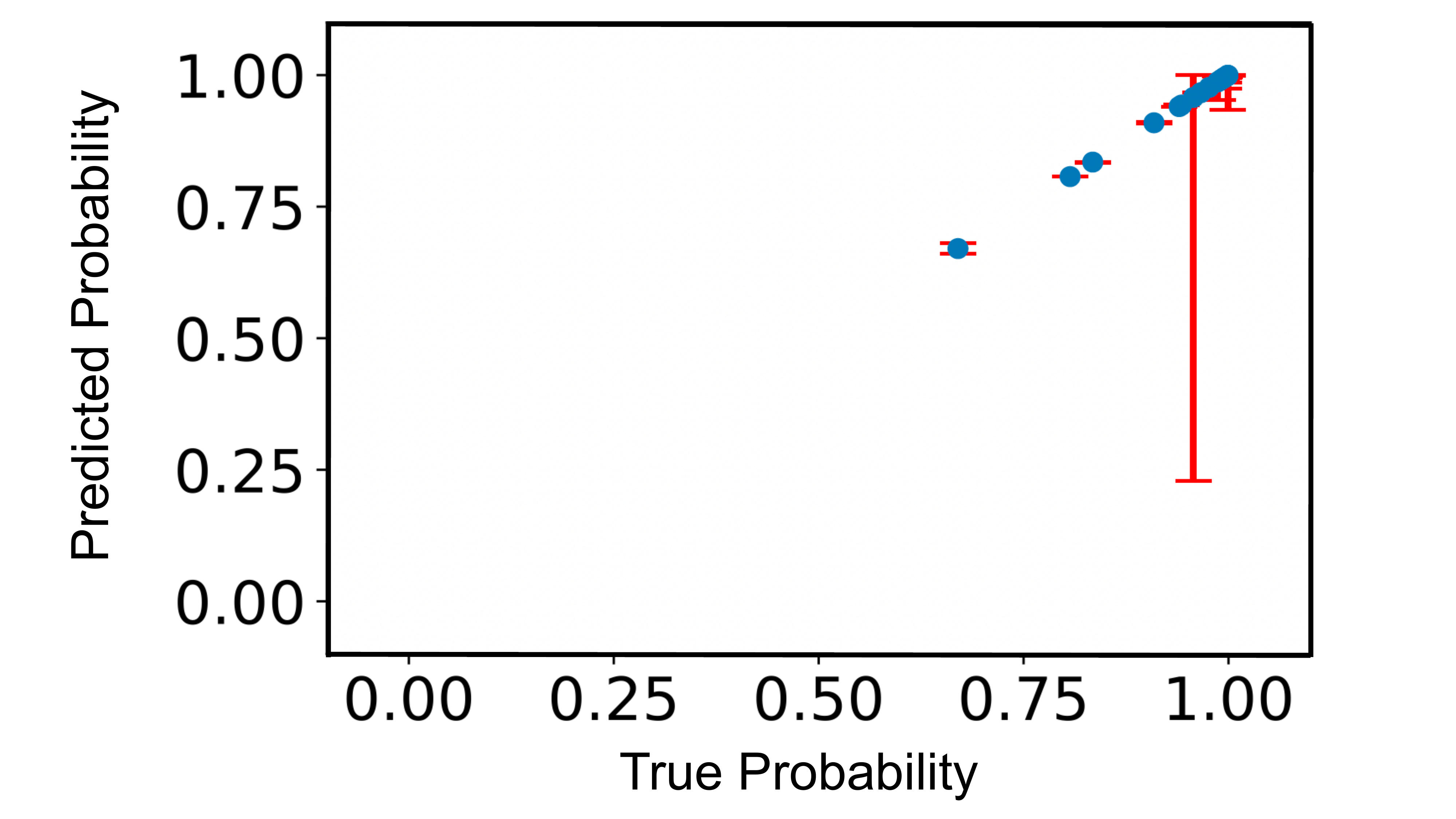}}
\end{minipage}
\;
\begin{minipage}[t]{0.3\textwidth}
\centerline{longest path with $K = 20$ }
\centerline{\includegraphics[scale=0.071]{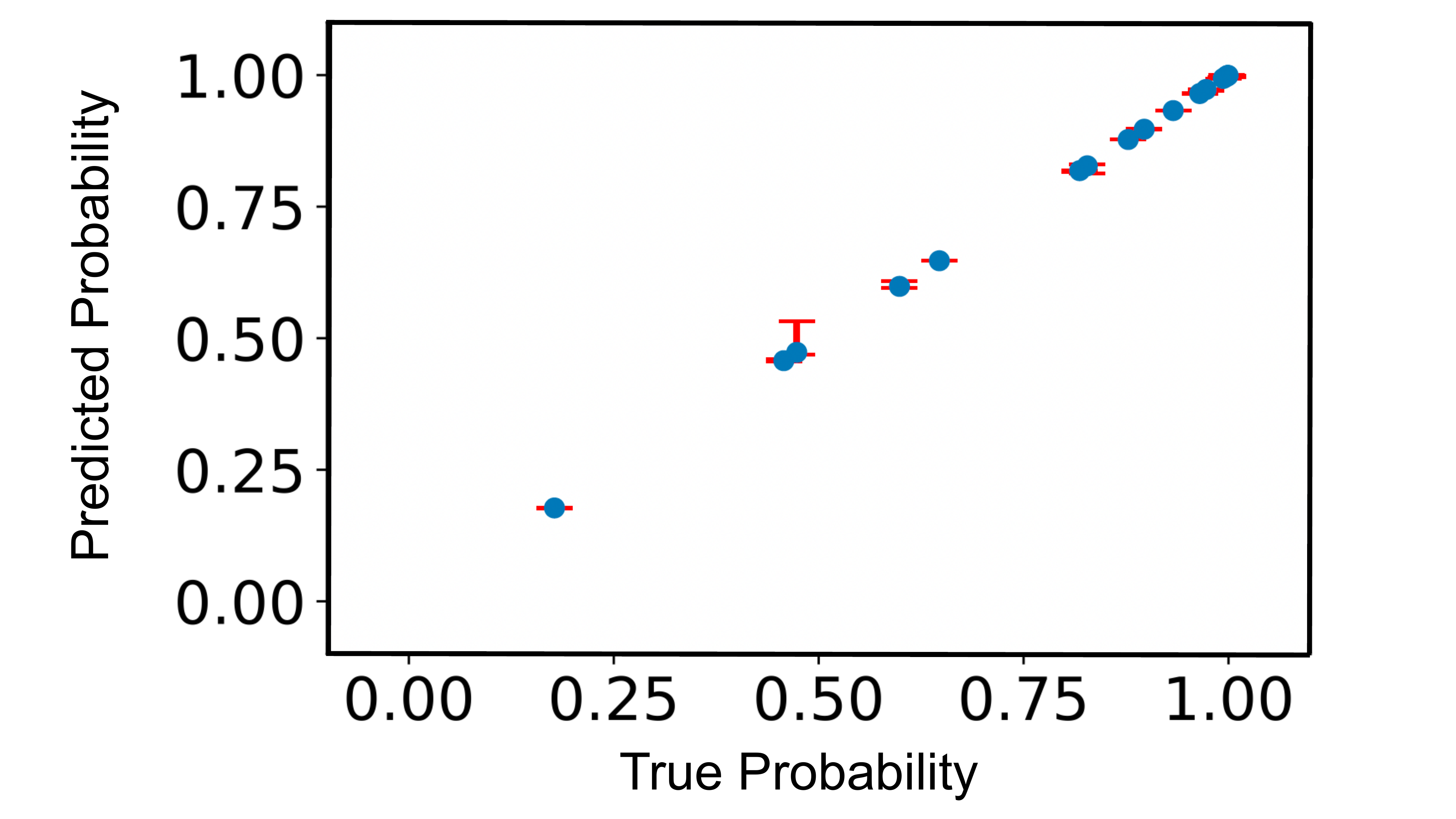}}
\end{minipage}
\caption{Prediction accuracy of nonparametric s-RAM on data generated from an underlying parametric s-RAM. Blue dots denote ground-truth values. Red ranges represent prediction intervals for the predicted probabilities by using nonparametric s-RAM.}
\label{fig:prediction_longestpath}
\end{figure}

The prediction results reported in Figure \ref{fig:prediction_longestpath} show that the proposed nonparametric approach accurately recovers the ground truth. The estimated probability of the selected edge concentrates around the true probability, and the resulting prediction intervals become progressively tighter as the number of historical graphs increases. These findings indicate that the nonparametric s-RAM effectively captures the underlying structure governing the longest path outcomes.

%The prediction results reported in Figure \ref{fig:prediction_longestpath} shows that the proposed nonparametric approach accurately recovers the ground truth. The estimated probability of the selected edge concentrate around the true probability, with prediction intervals narrowing as the number of historical graphs increases. These results demonstrate that the nonparametric s-RAM effectively captures the underlying ground truth.

\subsection{Shortest Path Problem}
Next, we evaluate the performance of s-RAM for the shortest path problem. Due to space limitations, we refer the reader to Section \ref{app:shortestpath} for details on the representability check, estimation, and prediction formulations, and to Section \ref{app:exp_data} for a complete description of the data generation.

\subsubsection{Explanatory ability of s-RAM in the shortest path problem}\par\noindent

\textit{Generation of historical $\mathcal{X}^{(k)}$ and $\bp^{(k)}$  instances.} We follow the same experimental setup as in the longest path problem described in Section \ref{sec:experiment_longestpath} to generate instances for evaluating the explanatory ability of s-RAM in the shortest path setting.  

\textit{Representability test and estimation with $\bp^\cK$.} 
The results of representability and estimation loss for the shortest path problem are reported in Table \ref{tab:explanatory_shortestpath}. Similar to the longest path case, Table \ref{tab:explanatory_shortestpath} shows that s-RAM performs well across all tested scenarios despite model misspecification, and exhibits even stronger explanatory ability for the shortest path problem. In particular, s-RAM achieves high representability for $K \leq 15$. For example, when $K=5$, at least $90\%$ of the instances are s-RAM representable across all correlation structures, and the representable fraction remains above $70\%$ even when $K=15$. Moreover, the associated estimation error is substantially smaller: the average absolute deviation loss is consistently on the order of $10^{-4}$ or smaller, even in settings where exact representability is violated. This represents an order-of-magnitude improvement compared to the longest path problem. As in the longest path, representability decreases with $K$ due to the increasing number of feasibility constraints in \eqref{sp-feascon} induced by additional historical observations. Overall, the framework remains robust across different correlation structures, with particularly stronger performance under independent and negatively correlated noise than positively correlated noise.

\begin{table}[htb]
\centering
\caption{Representation and estimation in the shortest path problem with Gaussian instances}
\label{tab:explanatory_shortestpath}
\scalebox{0.85}{
\begin{tabular}{|c|cc|cc|cc|}
\hline
\multirow{2}{*}{K} 
& \multicolumn{2}{c|}{independent gaussian} 
& \multicolumn{2}{c|}{gaussian with -ve correlation} 
& \multicolumn{2}{c|}{gaussian with +ve correlation} \\ \cline{2-7} 
\multicolumn{1}{|c|}{}                   
& \multicolumn{1}{c|}{\begin{tabular}[c]{@{}c@{}}fraction \\ of s-RAM\\ representable\\ instances\end{tabular}} 
& \multicolumn{1}{c|}{\begin{tabular}[c]{@{}c@{}}average absolute\\ deviation loss (se)\\ $(10^{-4})$\end{tabular}} 
& \multicolumn{1}{c|}{\begin{tabular}[c]{@{}c@{}}fraction \\ of s-RAM\\ representable\\ instances\end{tabular}} 
& \multicolumn{1}{c|}{\begin{tabular}[c]{@{}c@{}}average absolute\\ deviation loss (se)\\ $(10^{-4})$\end{tabular}} 
& \multicolumn{1}{c|}{\begin{tabular}[c]{@{}c@{}}fraction \\ of s-RAM\\ representable\\ instances\end{tabular}} 
& \multicolumn{1}{c|}{\begin{tabular}[c]{@{}c@{}}average absolute\\ deviation loss (se)\\ $(10^{-4})$\end{tabular}} \\ \hline

5  & \multicolumn{1}{c|}{0.90} & 0.62 (0.53) 
   & \multicolumn{1}{c|}{0.90} & 0.32 (0.23) 
   & \multicolumn{1}{c|}{0.95} & 0.04 (0.04) \\ \hline
10 & \multicolumn{1}{c|}{0.70} & 3.98 (2.25) 
   & \multicolumn{1}{c|}{0.95} & 0.85 (0.83) 
   & \multicolumn{1}{c|}{0.70} & 2.66 (1.38) \\ \hline
15 & \multicolumn{1}{c|}{0.70} & 3.61 (2.12) 
   & \multicolumn{1}{c|}{0.75} & 2.64 (2.09) 
   & \multicolumn{1}{c|}{0.75} & 3.60 (1.93) \\ \hline
20 & \multicolumn{1}{c|}{0.60} & 3.08 (1.32) 
   & \multicolumn{1}{c|}{0.65} & 1.02 (0.44) 
   & \multicolumn{1}{c|}{0.15} & 8.68 (3.40) \\ \hline
\end{tabular}}
\end{table}

\subsubsection{Predictive ability of s-RAM in the shortest path problem}
We again evaluate the predictive accuracy of the nonparametric s-RAM on counterfactual tasks using data generated from parametric s-RAMs. For the shortest path problem, we focus on predicting the bottleneck congestion (the maximum flow on any link) on a new graph $\bG^{\text{new}}$. 

\textit{Generation of training instance $\bp^{(k)}$ and ground truth in the shortest path.} We adopt a data generation procedure that is largely identical to that used for the longest path problem, except that we vary the parametric noise distributions in order to test the robustness of the nonparametric s-RAM to different distributional assumptions. Specifically, for the shortest path problem, the marginals follow non-identical uniform distributions, where the lower bounds are sampled uniformly from $(1,3)$, and the upper bounds are sampled uniformly from $(5,10)$. All other edge costs are sampled independently and uniformly from $[3,5]$. Directed acyclic graphs are generated using the same procedure as in the longest path experiments. Then, we use formulation \eqref{data:shortestpath} to generate training instance $\bp^{(k)}$. We evaluate the prediction performance of the nonparametric s-RAM across collection sizes $K \in \{5,10,20\}$. For each value of $K$, we generate $20$ training instances. Given a new graph $\bG^{\text{new}} = (s \cup t \cup \bV^{\text{new}}, \bE^{\text{new}})$, let $\bx^{\text{true}}$ denote the ground truth flow for $\bG^{\text{new}}$.
To generate $\bx^{\text{true}}$,  we use the same parameter values as in the training instances together with formulation \eqref{data:shortestpath} to compute the true flow on the new graph.  Then, the corresponding counterfactual ground-truth bottleneck congestion is then given by $\max_{e\in \bE^{\text{new}}}x_{e}^{\text{true}}$.

\begin{figure}[htb]
\centering
\begin{minipage}[t]{0.3\textwidth}
\centerline{shortest path with $K = 5$ }
\centerline{\includegraphics[scale=0.25]{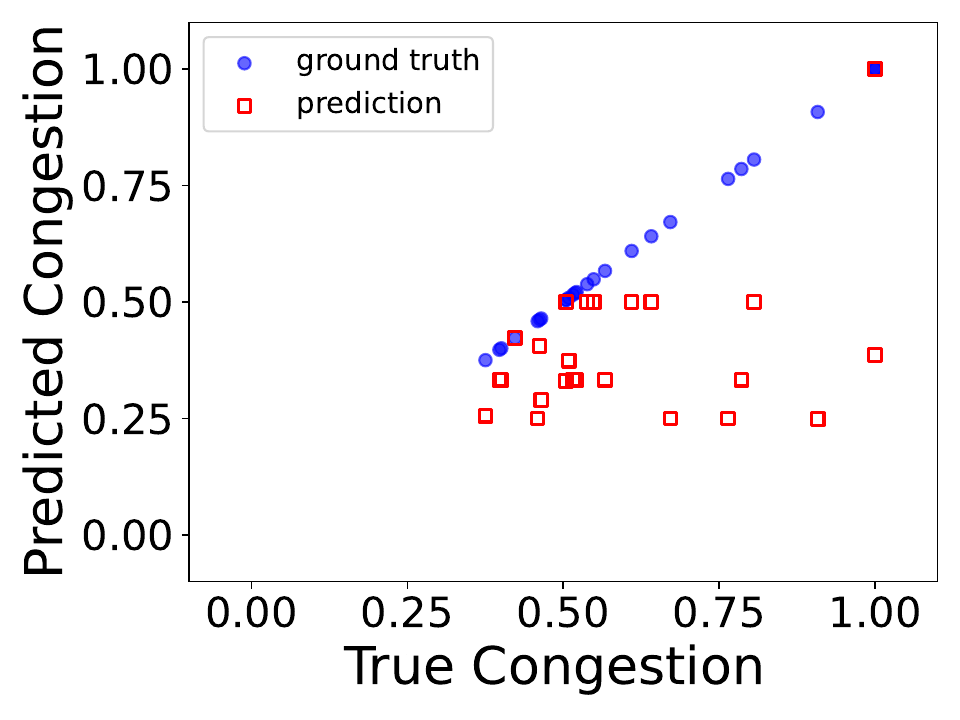}}
\end{minipage}
\;
\begin{minipage}[t]{0.3\textwidth} %\centerline{collection size = 40}
\centerline{shortest path with $K = 10$}
\centerline{\includegraphics[scale=0.25]{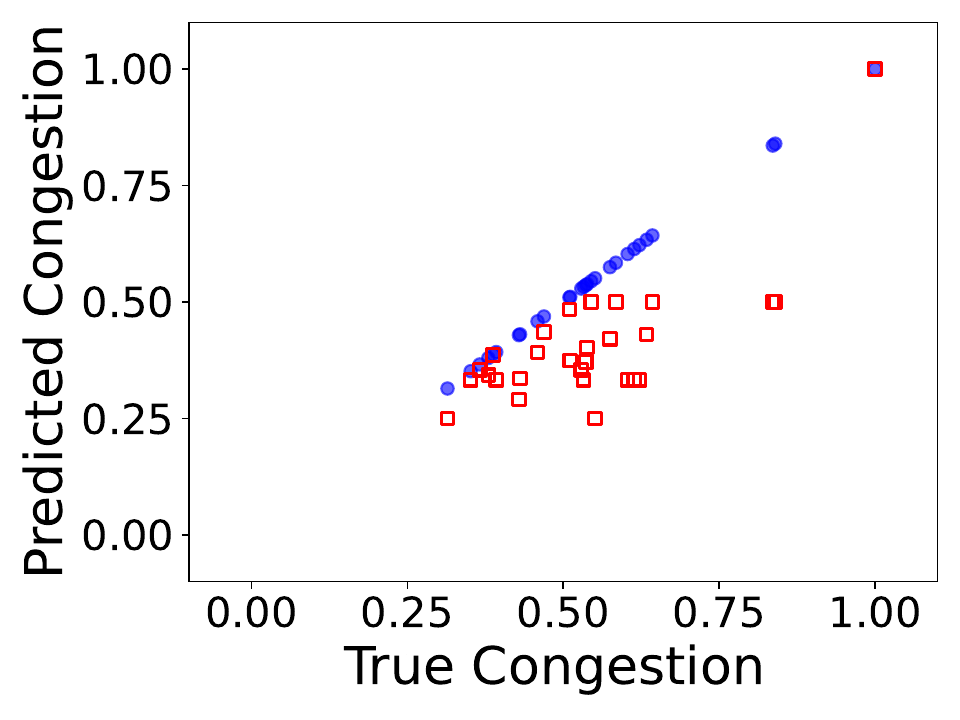}}
\end{minipage}
\;
\begin{minipage}[t]{0.3\textwidth} %\centerline{collection size =80}
\centerline{shortest path with $K = 20$ }
\centerline{\includegraphics[scale=0.25]{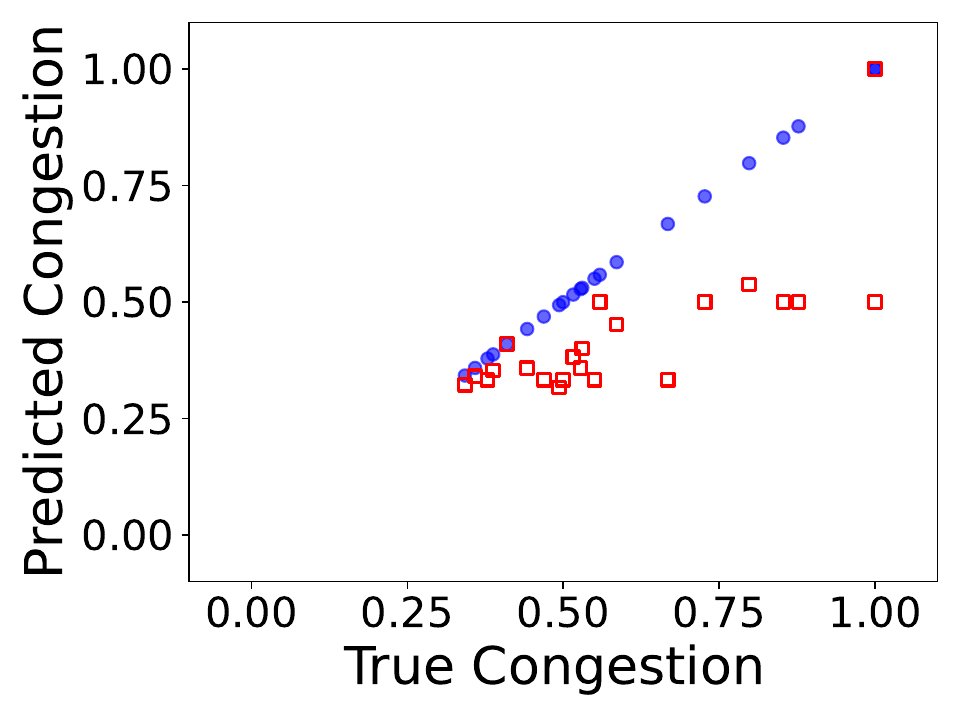}}
\end{minipage}
\caption{Prediction accuracy of nonparametric s-RAM on data generated from an underlying parametric s-RAM. Blue dots denote ground-truth values. Red dots indicate optimistic estimates of bottleneck congestion in the shortest path problem by using nonparametric s-RAM.}
\label{fig:prediction_shortestpath}
\end{figure}

\textit{Prediction with nonparametric s-RAM in the shortest path.}
For the shortest path problem, we use the s-RAM best-case estimate of the bottleneck congestion, defined as the minimum achievable maximum edge flow over all s-RAMs consistent with the observed data on the new graph. Formally, this quantity is given by $\min_{\bx \in \mathcal{U}{\text{new}}} \max_{e \in \bE^{\text{new}}} x_e$. Recall that $\mathcal{U}{\text{new}}$ denotes the feasible s-RAM nonparametric prediction region defined in Section \ref{sec:robust_prediction}.  As shown in Figure \ref{fig:prediction_shortestpath}, in the shortest path problem, the optimistic bottleneck congestion predicted by nonparametric s-RAM converges to the true bottleneck congestion as more graphs are observed.

\subsection{Assignment Problem}
Next, we evaluate the performance of s-RAM for the shortest path problem. Due to space limitations, we refer the reader to Section \ref{app:assignment} for details on the representability check, estimation, and prediction formulations, and to Section \ref{app:exp_data} for a complete description of the data generation.

\subsubsection{Explanatory ability of s-RAM in the assignment problem} \par\noindent

\textit{Generation of historical $\mathcal{X}^{(k)}$ instances.}
For the assignment problem, we consider bipartite graphs with $4$ nodes on the left and $6$ nodes on the right. The number of historical instances varies across $K \in \{5,10,15,20\}$. Edges are included independently with probabilities drawn uniformly from $[0.4,0.6]$, subject to the constraint that each node is incident to at least one edge, ensuring feasibility of the assignment problem. On average, the resulting graphs contain approximately $12$ edges.

\textit{Generation of historical $\bp^{(k)}$ instances.}
The generation of $\bp^{(k)}$ for the assignment problem follows the same procedure as in the longest path experiments, with edge costs given by $v_e + \epsilon_e$ under the same noise structures and parameter settings. The only difference is that the deterministic assignment problem is solved to obtain the observed data.

\textit{Representability test and estimation with $\bp^\cK$.} 
The results for the assignment problem are reported in Table \ref{tab:explanatory_assignment}. Consistent with the longest and shortest path experiments, s-RAM exhibits strong explanatory ability across all tested scenarios despite model misspecification. In particular, the fraction of s-RAM representable instances remains high throughout the tested range of $K$. For example, when $K=5$, at least $95\%$ of the instances are s-RAM representable under all noise structures, and representability remains above $60\%$ even when $K=20$. Moreover, the associated estimation error is extremely small: the average absolute deviation loss is consistently on the order of $10^{-5}$ or smaller across all settings, including cases where exact representability is not achieved. As in the path problems, representability decreases with $K$ due to the growing number of feasibility constraints induced by additional historical observations. Nevertheless, the assignment problem appears particularly well-suited to the s-RAM framework, with robust explanatory performance and negligible estimation error across a wide range of correlation structures.

\begin{table}[htb]
\centering
\caption{Representation and estimation in the assignment problem with Gaussian instances}
\label{tab:explanatory_assignment}
\scalebox{0.85}{
\begin{tabular}{|c|cc|cc|cc|}
\hline
\multirow{2}{*}{K} 
& \multicolumn{2}{c|}{independent gaussian} 
& \multicolumn{2}{c|}{gaussian with -ve correlation} 
& \multicolumn{2}{c|}{gaussian with +ve correlation} \\ \cline{2-7} 
\multicolumn{1}{|c|}{}                   
& \multicolumn{1}{c|}{\begin{tabular}[c]{@{}c@{}}fraction \\ of s-RAM\\ representable\\ instances\end{tabular}} 
& \multicolumn{1}{c|}{\begin{tabular}[c]{@{}c@{}}average absolute\\ deviation loss (se)\\ $(10^{-5})$\end{tabular}} 
& \multicolumn{1}{c|}{\begin{tabular}[c]{@{}c@{}}fraction \\ of s-RAM\\ representable\\ instances\end{tabular}} 
& \multicolumn{1}{c|}{\begin{tabular}[c]{@{}c@{}}average absolute\\ deviation loss (se)\\ $(10^{-5})$\end{tabular}} 
& \multicolumn{1}{c|}{\begin{tabular}[c]{@{}c@{}}fraction \\ of s-RAM\\ representable\\ instances\end{tabular}} 
& \multicolumn{1}{c|}{\begin{tabular}[c]{@{}c@{}}average absolute\\ deviation loss (se)\\ $(10^{-5})$\end{tabular}} \\ \hline

5  & \multicolumn{1}{c|}{0.95} & 0.50 (0.49)
   & \multicolumn{1}{c|}{1.00} & 0.00 (0.00)
   & \multicolumn{1}{c|}{0.95} & 0.10 (0.10) \\ \hline
10 & \multicolumn{1}{c|}{0.90} & 1.55 (1.24)
   & \multicolumn{1}{c|}{0.95} & 1.40 (1.36)
   & \multicolumn{1}{c|}{0.85} & 3.15 (2.92) \\ \hline
15 & \multicolumn{1}{c|}{0.85} & 4.03 (2.45)
   & \multicolumn{1}{c|}{0.85} & 2.50 (1.82)
   & \multicolumn{1}{c|}{0.90} & 0.20 (0.14) \\ \hline
20 & \multicolumn{1}{c|}{0.75} & 2.38 (1.75)
   & \multicolumn{1}{c|}{0.60} & 8.53 (4.40)
   & \multicolumn{1}{c|}{0.60} & 1.97 (1.13) \\ \hline

\end{tabular}}
\end{table}

\subsubsection{Predictive ability of s-RAM in the assignment problem}

We evaluate the predictive accuracy of the nonparametric s-RAM on counterfactual tasks for the assignment problem using data generated from parametric s-RAMs. In this setting, we examine whether s-RAM can accurately predict the total system welfare on a new graph when some nodes become unavailable.

\textit{Generation of training instances $\bp^{(k)}$ and ground truth in the assignment problem.}
In assignment, we let the noise marginals be non-identical and uniformly distributed, with lower bounds drawn from $[1,2]$ and upper bounds drawn from $[5,10]$. Counterfactual instances are generated by randomly designating a subset of nodes as unavailable in the bipartite graph and removing all edges incident to those nodes. Then we use formulation \eqref{data:assignment} to generate training instances $\bp^{(k)}$. Specifically, let $\bv$ denote the deterministic weight matrix, and let $\mathbf{a}$ and $\mathbf{b}$ denote the matrices of lower and upper bounds of the uniform noise distributions, respectively. Given a new bipartite graph $\bG^{\text{new}} = (\bV_L^{\text{new}}, \bV_R^{\text{new}}, \bE^{\text{new}})$, let $\bx^{\text{true}}$ denote the ground-truth matching rates for $\bG^{\text{new}}$. We then compute the ground-truth system welfare for each counterfactual instance with formulation \eqref{data:assignment} under the same parameter values as in the training instances. %Specifically, let $\bv$ denote the weight, $\textbf{a}$ be the matrix of lower bounds and  $\textbf{b}$ be the matrix of upper bounds of the uniform distribution. Let $\bx^{\text{true}}$ denote the ground truth matching rate. Given a new bipartite graph $\bG^{\text{new}} = (\bV_L^{\text{new}}, \bV_R^{\text{new}}, \bE^{\text{new}})$, the ground truth system welfare is 
The corresponding ground-truth system welfare is given by
$$\sum_{i\in \bV_L^{\text{new}} }\sum_{j \in \bV_R^{\text{new}} } (\nu_{ij}+ b_{ij})x_{ij}^{\text{true}} - \sum_{i\in \bV_L^{\text{new}} }\sum_{j \in \bV_R^{\text{new}} }  \frac{b_{ij} - a_{ij}}{2} {(x_{ij}^{\text{true}})}^2.$$

\textit{Prediction with nonparametric s-RAM in the assignment problem.}
For the assignment problem, computing the worst-case welfare prediction leads to a nonconvex mixed-integer optimization problem, which is computationally challenging and may be unreliable when solved using standard solvers. We therefore focus on the best-case s-RAM estimate for the total system welfare, defined as the maximum achievable welfare over all s-RAMs consistent with the observed data and the counterfactual graph structure. Formally, this quantity is given by 
$$\max_{\bx \in \mathcal{U}{\text{new}}} \sum_{i\in \bV_L^{\text{new}} }\sum_{j \in \bV_R^{\text{new}} } (\nu_{ij}+ b_{ij})x_{ij} - \sum_{i\in \bV_L^{\text{new}} }\sum_{j \in \bV_R^{\text{new}} }  \frac{b_{ij} - a_{ij}}{2} x_{ij}^2 .$$ 
Recall that $\mathcal{U}{\text{new}}$ denotes the feasible s-RAM nonparametric prediction region defined in Section \ref{sec:robust_prediction}.

\begin{figure}[htb]
\centering
\begin{minipage}[t]{0.3\textwidth}
\centerline{assignment with $K = 5$ }
\centerline{\includegraphics[scale=0.25]{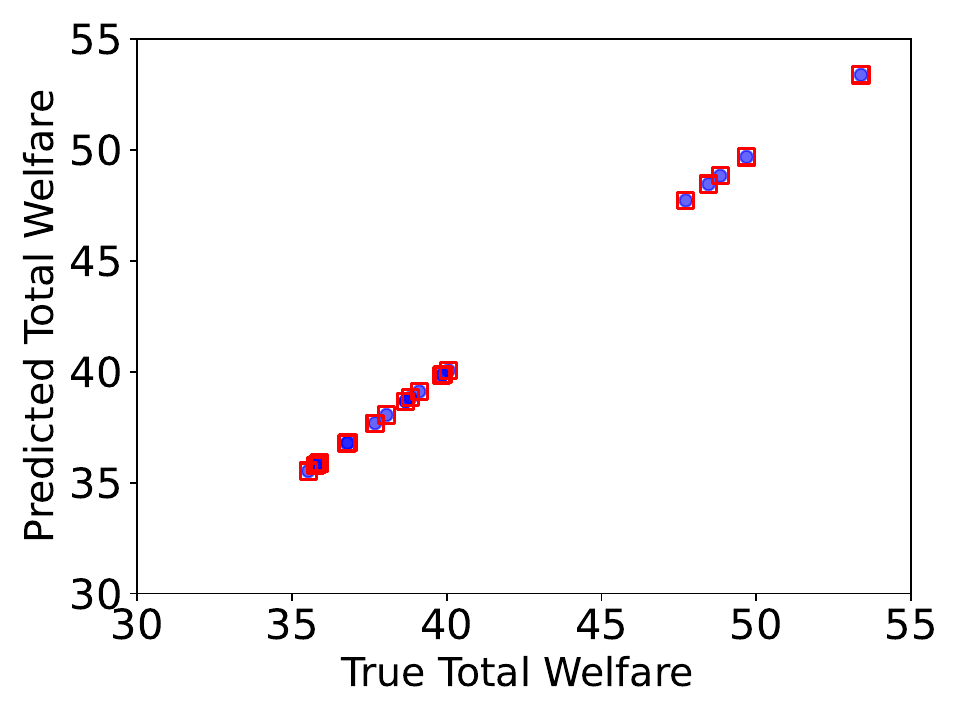}}
\end{minipage}
\;
\begin{minipage}[t]{0.3\textwidth} 
\centerline{assignment with $K = 10$ }
\centerline{\includegraphics[scale=0.25]{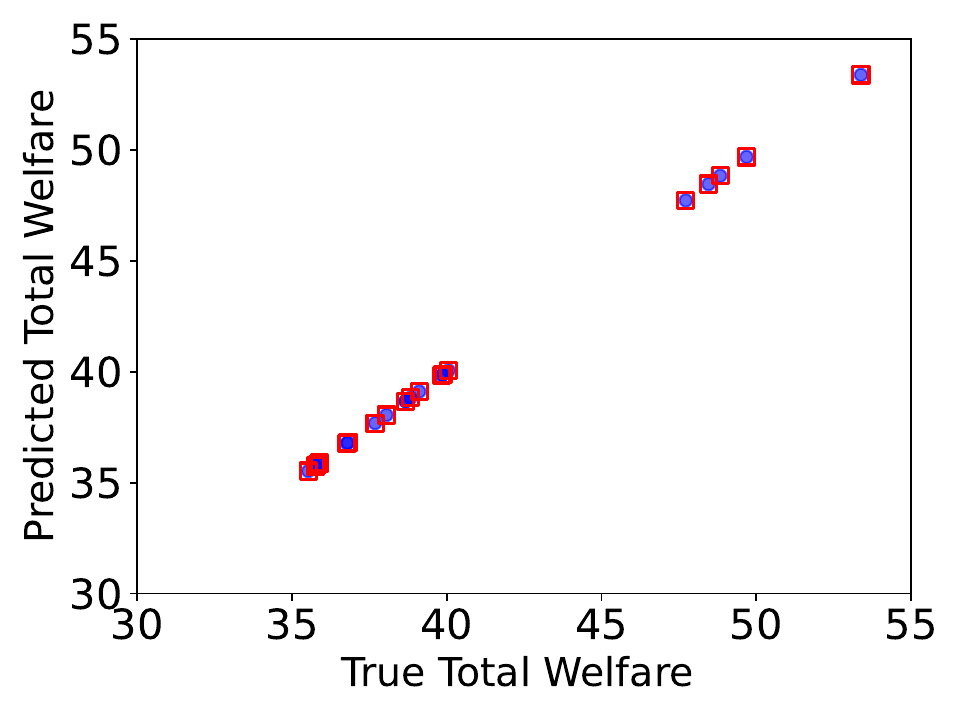}}
\end{minipage}
\;
\begin{minipage}[t]{0.3\textwidth} 
\centerline{assignment with $K = 20$ }
\centerline{\includegraphics[scale=0.25]{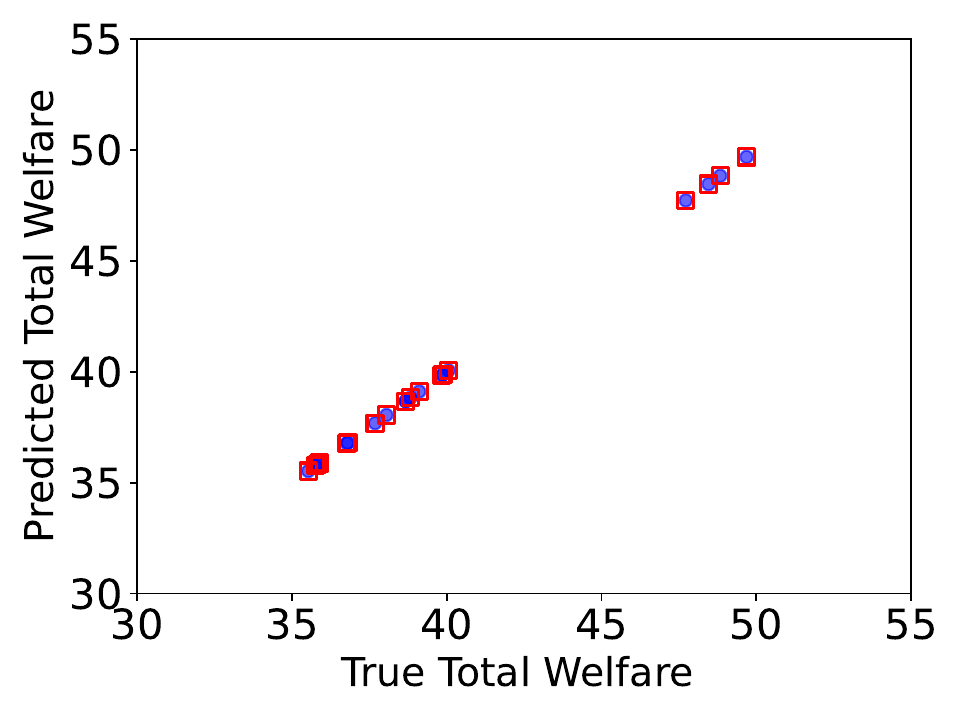}}
\end{minipage}
\caption{Prediction accuracy of nonparametric s-RAM on data generated from an underlying parametric s-RAM. Blue dots denote ground-truth values. Red squares indicate optimistic total welfare in the assignment problem by using nonparametric s-RAM.}
\label{fig:prediction_assignment}
\end{figure}

As shown in Figure \ref{fig:prediction_assignment}, across all test instances, the best-case welfare predicted by s-RAM closely matches the ground-truth system welfare, indicating strong predictive performance of the nonparametric approach in the assignment setting.

% \subsubsection{Predictive ability of s-RAM in the assignment problem}
% For the assignment problem, we examine whether s-RAM can predict the total system welfare on a new graph with some nodes being unavailable. 
% \textit{Generation of training instance $\bp^{(k)}$ and ground truth in the assignment problem.} 
% In assignment, marginals of noise follow nonidentical uniform distributions in the form of $\text{Uniform}(a, b)$ and $a$ being sampled uniformly from \([1, 2]\) and $b$ being sampled uniformly from \([5, 10]\). We generate the bipartite graphs by randomly assigning some nodes to be unavailable and removing all the associated edges. 
% \textit{Prediction with nonparametric s-RAM in the assignment problem.}
% For assignment, predicting the worst-case welfare leads to a nonconvex mixed integer program, which is computationally challenging and might not be accurate when solving by existing solvers. We use the best-case estimate for the total welfare.
% In the assignment problem, the predicted total welfare closely matches the ground truth across all test instances.

\subsection{Constrained Subset Selection Problem}
Next, we evaluate the performance of s-RAM for the shortest path problem. Due to space limitations, we refer the reader to Section \ref{app:subsetselection} for details on the representability check, estimation, and prediction formulations, and to Section \ref{app:exp_data} for a complete description of the data generation.

\subsubsection{Explanatory ability of s-RAM in the constrained subset selection problem}\par\noindent

\textit{Generation of historical $\mathcal{X}^{(k)}$ instances.}
For the constrained subset selection problem, we consider a ground set of available items $\cN$ with $|\cN| = 10$. We focus on a cardinality-constrained setting, where 
$\mathcal{X}^{(k)} = \left\{\bx: \sum_{i\in \cI^{(k)}} x_i \leq C, \  x_i\in \{0,1\}, \ \forall i \in \cI^{(k)}\right\}$ 
with $C = 2$ and $\cI^{(k)} \subseteq \cN$ for each historical instance $k$. The number of historical instances varies across $K \in \{5,10,15,20\}$. For each instance, we randomly generate the available item set $\cI^{(k)}$ with cardinality $|\cI^{(k)}| \in \{3,4\}$.

\textit{Generation of historical $\bp^{(k)}$ instances.}
The generation of $\bp^{(k)}$ for the constrained subset selection problem follows the same procedure as in the longest path experiments, with item costs given by $v_i + \epsilon_i$ under the same noise structures and parameter settings as edge costs. The only difference is that the deterministic constrained problem is solved to obtain the observed data.  

\textit{Representability test and estimation with $\bp^\cK$.} 
The results for the constrained subset selection problem are reported in Table \ref{tab:explanatory_subsetselection}. Compared with the longest path, shortest path, and assignment problems, s-RAM exhibits weaker explanatory performance in this setting, which is expected since the data are generated outside the s-RAM framework and thus constitute a stress test. In particular, s-RAM is able to represent a nontrivial fraction of instances for small $K$, but representability deteriorates more rapidly as $K$ increases, reflecting the increasing number of feasibility constraints induced by additional historical observations. Nevertheless, the estimation loss remains moderate: the average absolute deviation loss stays on the order of $10^{-2}$ across all settings. This larger scale is natural, as outcomes in the constrained subset selection problem are not normalized to sum to one, unlike the flow- and matching-based problems considered earlier. Overall, despite model misspecification and the absence of normalization, s-RAM continues to provide reasonable explanatory performance in this more challenging setting.

\begin{table}[htb]
\centering
\caption{Representation and estimation in the constrained subset selection problem with Gaussian instances}
\label{tab:explanatory_subsetselection}
\scalebox{0.85}{
\begin{tabular}{|c|cc|cc|cc|}
\hline
\multirow{2}{*}{K} 
& \multicolumn{2}{c|}{independent gaussian} 
& \multicolumn{2}{c|}{gaussian with -ve correlation} 
& \multicolumn{2}{c|}{gaussian with +ve correlation} \\ \cline{2-7} 
\multicolumn{1}{|c|}{}                   
& \multicolumn{1}{c|}{\begin{tabular}[c]{@{}c@{}}fraction \\ of s-RAM\\ representable\\ instances\end{tabular}} 
& \multicolumn{1}{c|}{\begin{tabular}[c]{@{}c@{}}average absolute\\ deviation loss (se)\\ $(10^{-2})$\end{tabular}} 
& \multicolumn{1}{c|}{\begin{tabular}[c]{@{}c@{}}fraction \\ of s-RAM\\ representable\\ instances\end{tabular}} 
& \multicolumn{1}{c|}{\begin{tabular}[c]{@{}c@{}}average absolute\\ deviation loss (se)\\ $(10^{-2})$\end{tabular}} 
& \multicolumn{1}{c|}{\begin{tabular}[c]{@{}c@{}}fraction \\ of s-RAM\\ representable\\ instances\end{tabular}} 
& \multicolumn{1}{c|}{\begin{tabular}[c]{@{}c@{}}average absolute\\ deviation loss (se)\\ $(10^{-2})$\end{tabular}} \\ \hline

5  & \multicolumn{1}{c|}{0.75} & 0.35 (0.20)
   & \multicolumn{1}{c|}{0.90} & 0.13 (0.09)
   & \multicolumn{1}{c|}{0.85} & 0.25 (0.21) \\ \hline
10 & \multicolumn{1}{c|}{0.50} & 0.57 (0.20)
   & \multicolumn{1}{c|}{0.45} & 0.73 (0.33)
   & \multicolumn{1}{c|}{0.45} & 0.60 (0.19) \\ \hline
15 & \multicolumn{1}{c|}{0.15} & 0.66 (0.13)
   & \multicolumn{1}{c|}{0.10} & 0.98 (0.17)
   & \multicolumn{1}{c|}{0.20} & 0.67 (0.17) \\ \hline
20 & \multicolumn{1}{c|}{0.05} & 1.96 (0.36)
   & \multicolumn{1}{c|}{0.10} & 1.29 (0.27)
   & \multicolumn{1}{c|}{0.00} & 1.31 (0.21) \\ \hline

\end{tabular}}
\end{table}

\subsubsection{Predictive ability of s-RAM in the constrained subset selection problem}

We evaluate the predictive accuracy of the nonparametric s-RAM on counterfactual tasks for the constrained subset selection problem under a cardinality constraint. In this setting, we examine whether s-RAM can predict the counterfactual marginal allocation of a selected item, defined as its optimal fractional weight in the solution to the constrained subset selection problem on a new instance.

\textit{Generation of training instances $\bp^{(k)}$ and ground truth for the constrained subset selection problem.}
In the constrained subset selection problem, all items are assigned identical deterministic weights, while the noise marginals follow non-identical exponential distributions with rates sampled uniformly from $(1,2)$. Then, we use formulation \eqref{data:subsetselection} to generate training instances $\bp^{(k)}$. We evaluate the prediction performance of the nonparametric s-RAM across collection sizes $K \in \{20,40,80\}$, generating $20$ training instances for each value of $K$. For a new instance, the ground-truth marginal allocation is obtained by solving the corresponding constrained subset selection problem under the same parameter values with formulation \eqref{data:subsetselection}. Let $\bx^{\text{true}}$ denote the resulting vector of marginal allocation of items in a new instance $\cI^{(new)}$, and let $i^*$ be the item of interest. The corresponding ground-truth quantity is then $x^{\text{true}}_{i^*}$.  

\textit{Prediction with nonparametric s-RAM in the constrained subset selection problem.}
For the constrained subset selection problem, we compute the s-RAM worst- and best-case estimates for the marginal allocation of the selected feature over all s-RAMs consistent with the observed data. Given a new set  $\cI^{(new)}$, the worst case estimate is given as $\min_{\bx \in \mathcal{U}{\text{new}}} x_{i^*}$ and the best case estimate is given as $\max_{\bx \in \mathcal{U}{\text{new}}} x_{i^*}$.  Recall that $\mathcal{U}{\text{new}}$ denotes the feasible s-RAM nonparametric prediction region defined in Section \ref{sec:robust_prediction}. 

\begin{figure}[htb]
\centering
\begin{minipage}[t]{0.3\textwidth}
\centerline{subset selection with $K = 20$ }
\centerline{\includegraphics[scale=0.25]{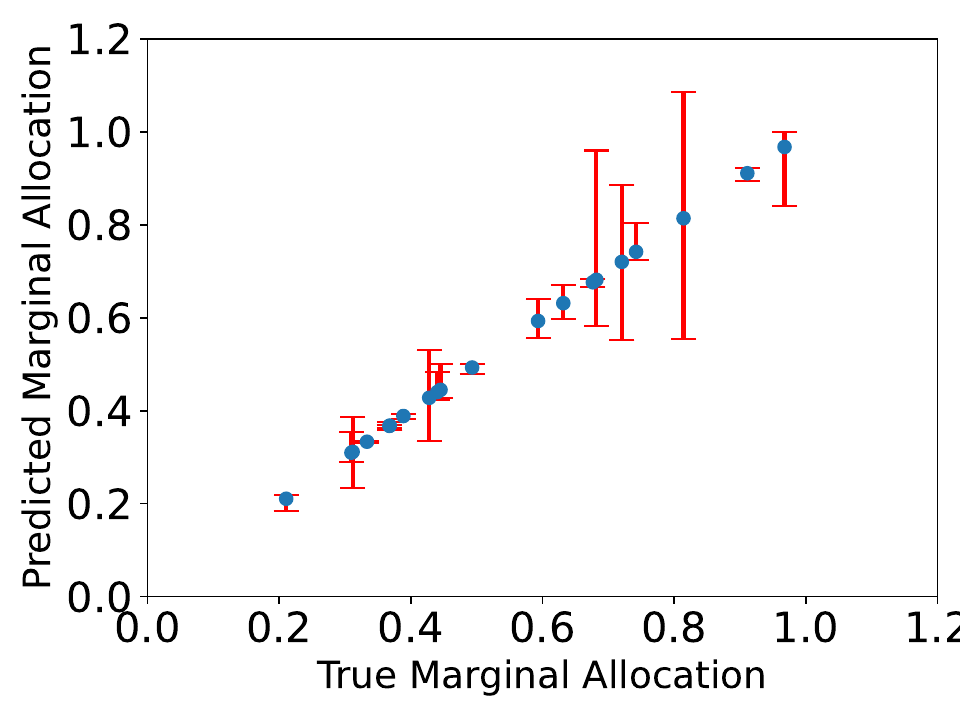}}
\end{minipage}
\;
\begin{minipage}[t]{0.3\textwidth} %\centerline{collection size = 40}
\centerline{subset selection  with $K = 40$ }
\centerline{\includegraphics[scale=0.25]{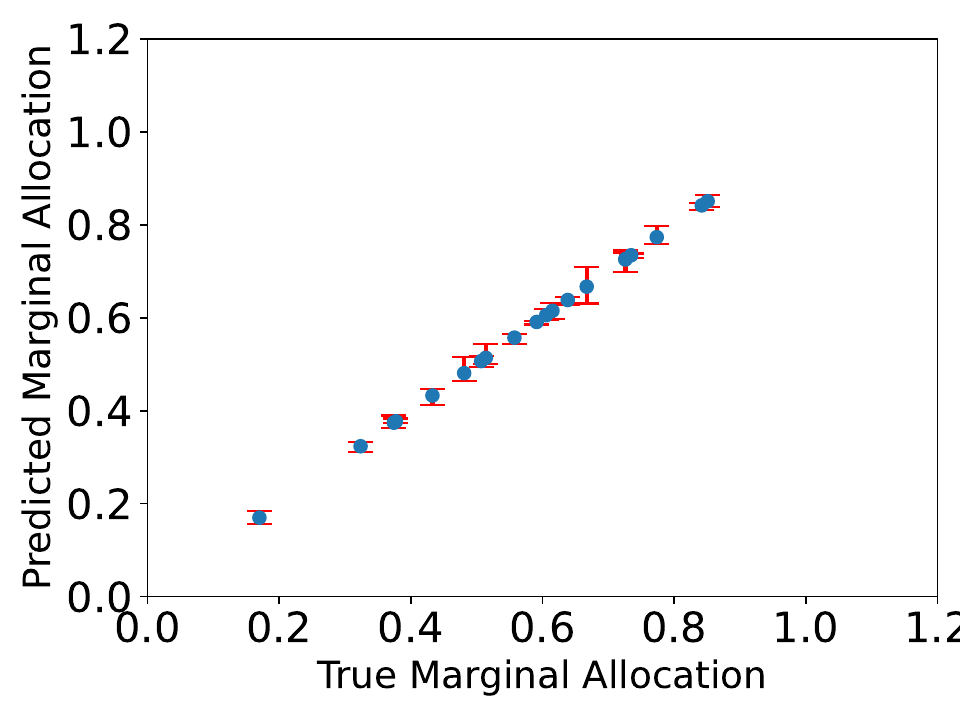}}
\end{minipage}
\;
\begin{minipage}[t]{0.3\textwidth} %\centerline{collection size =80}
\centerline{subset selection with $K = 80$ }
\centerline{\includegraphics[scale=0.25]{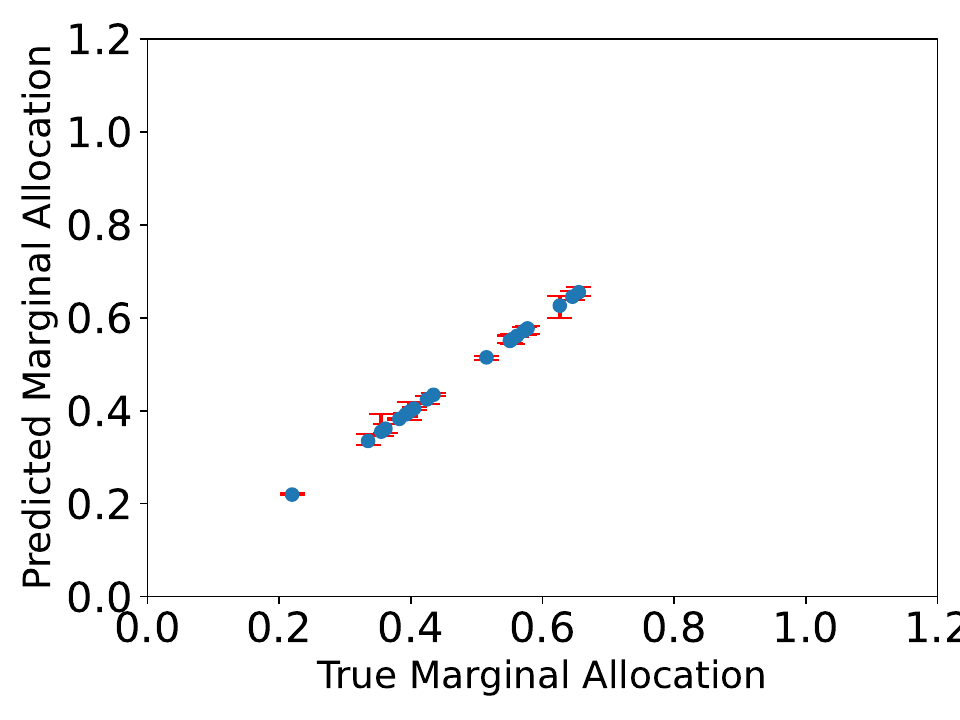}}
\end{minipage}
\caption{Prediction accuracy of nonparametric s-RAM on data generated from an underlying parametric s-RAM. Blue dots denote ground-truth values. Red ranges represent prediction intervals for the predicted marginal allocation in the constrained subset selection problem by using nonparametric s-RAM.}
\label{fig:prediction_subsetselection}
\end{figure}

Figure \ref{fig:prediction_subsetselection} shows that the estimated marginal allocations closely track the ground truth values across all tested settings. Moreover, the prediction intervals become progressively tighter as the number of historical instances increases, indicating that the nonparametric s-RAM effectively captures the underlying structure of the constrained subset selection problem.

%As shown in Figure \ref{fig:prediction_subsetselection}, the prediction results show that the estimated marginal allocations closely track the ground-truth values across all tested settings. Moreover, the prediction intervals become progressively tighter as the number of historical instances increases, indicating that the nonparametric s-RAM effectively captures the underlying structure of the constrained subset selection problem.

% for constrained subset selection:
% \textcolor{blue}{For the  constrained subset selection problem under the uniform matroid constraint, we test whether s-RAM can predict the counterfactual marginal allocation of a key feature, i.e., its optimal fractional weight in the solution to the  constrained subset selection problem on a new instance.}

% \textcolor{blue}{For constrained subset selection, we compute the s-RAM worst- and best- case estimates for the marginal allocation of the selected feature.}

% Finally, in the  constrained subset selection problem, the estimated marginal allocations recover the true marginals, with prediction intervals shrinking as the number of historical observations increases.

\section{Conclusion and Future Work}\label{sec:conclusion}

In this paper, we develop a nonparametric framework based on a representative agent model to infer how aggregate selection behavior may vary across new combinatorial choice environments unobserved in historical data. Our first main contribution is an exact characterization of the selection probabilities representable under the s-RAM model over general 0–1 polytopes, and a demonstration that consistency verification can be reduced to solving a polynomial-sized linear program. This characterization, in turn, enables a tractable nonparametric approach for generating predictions in counterfactual environments. For data that are inconsistent with the s-RAM model, we introduce a compact mixed-integer convex formulation to compute the best-fitting s-RAM approximation. We conduct extensive numerical experiments that demonstrate the strong empirical performance of the proposed methods. 

We conclude by highlighting  several challenges for future research: Developing advanced methods for handling large scale instances where the proposed integer programs formulations may fail to scale well is an important direction for future research. Another direction is to extend the framework to settings where the data is from a polytope for which a $0$–$1$ representation is not available, such as the Traveling Salesman Problem (TSP), or even beyond $0$–$1$ polytopes, as in the knapsack problem, thereby accommodating richer combinatorial structures, and broadening the framework’s applicability.

% Bibliography
\bibliographystyle{ACM-Reference-Format}
\bibliography{ref}

% Appendix
\appendix
\section{Missing Proofs of the Paper} \label{app:proof}
%Technical appendices with additional results, figures, graphs and proofs may be submitted with the paper submission before the full submission deadline (see above), or as a separate PDF in the ZIP file below before the supplementary material deadline. There is no page limit for the technical appendices.
% Text of your paper here
%\subsection{Proof of Theorem \ref{thm:mdm-general-feascon}}

\begin{proof}{Theorem \ref{thm:mdm-general-feascon}}
\textit{Necessity of \eqref{eq:mdm-feascon-order}}: 
Suppose $\cP^\cK$ is s-RAM-representable. 
Then there exist perturbation functions $\{c_j:j\in[n] \}$ and deterministic utilities $\{\nu_{j}: j\in[n] \}$ such that for any $k \in \cK,$  the given $(p_{j}^k: j\in\cI^{(k)}, k \in \cK)$, and the respective Lagrange multipliers $\{\lambda_{i}^k :i\in [m],k\in \cK\}$, are obtainable by solving the optimality conditions $\eqref{s-RAM-optcon}.$   That is, there exist $\{\lambda_{i}^k :i\in [m],k\in \cK\}$ for some fixed choice of $\{c_j:j\in[n] \}$ and  $\{\nu_{j}: j\in[n] \}$, for any $k\in \cK$, such that  
\begin{align}
     & \nu_j  - c_j^\prime(p_j^k) + \sum_{i\in [m]} \lambda_i^k a_{ij} + \mu_j^k - w_j^k= 0,\quad \forall j\in \cI^{(k)},\label{eq:general-optcon1}\\
    %& \lambda_i^k (b_i - \sum_{j\in \cI^{(k)} } a_{ij} p_j^k ) = 0, \quad \forall i \in [m], \label{eq:general-optcon2}\\
    & \mu_j^k x_j^* = 0, \ w_j^k (1-x_j^*) = 0, \quad \forall  j \in \cI^{(k)}. \label{eq:general-optcon3}
\end{align}
%The complementary slackness conditions \eqref{eq:general-optcon2} implies \eqref{eq:mdm-feascon-cs} because we have $b_i - \sum_{j\in \cI^{(k)}} a_{ij} p_j^{k} \leq 0,$ for all $i \in [m]$. Then, when $b_i - \sum_{j\in [n]} a_{ij} p_j^k  > 0$, we must have $\lambda_i^k = 0$ . In the following, we focus on showing \eqref{eq:mdm-feascon-order} is necessary. 
For each $j\in[n]$ and any two distinct observations  $k,l \in \cK$ containing $j$ as a common variable which is not forced to be zero, i.e., $j \in \cI^{(k)}\cap \cI^{(l)}$, we have
\begin{align*}
    & \sum_{i\in [m]} \lambda_i^k a_{ij} - \nu_j = c_j^\prime(p_j^k) + \mu_j^k - w_j^k \mbox{~ and ~}\\
    & \sum_{i\in [m]} \lambda_i^l a_{ij} - \nu_j =c_j^\prime(p_j^l) + \mu_j^l - w_j^l
\end{align*}
%$$\sum_{i\in [m]} \lambda_i^k a_{ij} - \nu_j = F_j^{-1}(1-p_j^k) + \mu_j^k - w_j^k  \mbox{~ and ~}  \sum_{i\in [m]} \lambda_i^k a_{ij} - \nu_j = F_j^{-1}(1-p_j^l) -\sum_{i\in [m]} \lambda_i^l a_{ij}  + \mu_j^l - w_j^l.$$
From Assumption \ref{asp:convex_perturbation}, since $c_j:[0,1]\to \mathbb{R}$ for $j \in [n]$ are univariate, strictly convex and continuously differentiable perturbation functions, we have $c^\prime_j:[0,1]\to \mathbb{R}$ for $j \in [n]$ are strictly increasing functions.  If $0\leq p_{j}^k < p_{j}^l \leq 1$, we have $c^\prime_j(p_{j}^k) < c^\prime_j(p_j^l)$ since $c^\prime_j(p)$ is a strictly increasing function over $p \in [0,1]$, and $\mu_j^k \geq 0$, $w_j^k =0$ and $\mu_j^l = 0$, $w_j^l \geq 0$ because of the complementary slackness conditions  \eqref{eq:general-optcon3}. Then we obtain 
\[
\sum_{i\in [m]} \lambda_i^k a_{ij}  - \nu_{j} \   = \    c^\prime_{j}(p_{j}^k)\  < \   c^{\prime}_{j}(p_{j}^l)  \ =\  \sum_{i\in [m]} \lambda_i^k a_{ij}  - \nu_{j}.
\]
%$\sum_{i\in [m]} \lambda_i^k a_{ij}  - \nu_{j} \   = \    c^\prime_{j}(p_{j}^k)\  < \   c^{\prime}_{j}(p_{j}^l)  \ =\  \sum_{i\in [m]} \lambda_i^k a_{ij}  - \nu_{j}.$
Adding $\nu_{j}$ on both sides, we obtain that the Lagrange multipliers should satisfy $\sum_{i\in [m]} \lambda_i^k a_{ij}  <  \sum_{i\in [m]} \lambda_i^k a_{ij}.$ If on the other hand, $0< p_{j}^k = p_{j}^l < 1$, we have $\mu_j^k = w_j^k = \mu_j^l = w_j^l = 0$, because of complementary slackness. Then, we have 
\[
\sum_{i\in [m]} \lambda_i^k a_{ij}  - \nu_{j} \   = \    c^{\prime}_{j}(p_{j}^k)\  = \   c^{\prime}_{j}(p_{j}^l)  \ =\  \sum_{i\in [m]} \lambda_i^k a_{ij}  - \nu_{j}.
\]
Again, adding $\nu_{j}$ on both sides, we obtain that $\{\lambda_{i}^k :i\in [m],k\in \cK\}$ should satisfy $\sum_{i\in [m]}\lambda_i^k a_{ij}   = \sum_{i\in [m]}\lambda_i^l a_{ij}.$  So, \eqref{eq:mdm-feascon-order} is also necessary for $(p_{j}^k: j\in\cI^{(k)}, k \in \cK)$ being s-RAM-representable.

\textit{Sufficiency of \eqref{eq:mdm-feascon-order}}: 
Given $(p_{j}^k: j\in\cI^{(k)}, k \in \cK)$ and  $(\lambda_i^k \in \mathbb{R}: i\in [m], k \in \cK)$ such that \eqref{eq:mdm-feascon-order} holds for all $j\in [n]$, $k\in \cK$ , we next exhibit a construction of perturbation functions $(c_{j}: j\in [n])$ and deterministic utilities $(\nu_j: j\in[n])$ for s-RAM. This construction will be such that it yields the given $(p_{j}^k: j\in\cI^{(k)}, k \in \cK)$ as the corresponding solutions from the optimality conditions in \eqref{s-RAM-optcon}, for any observation $k \in \cK.$ 

For any $j \in [n],$ let $\mathcal{K}_{j} = \{ k \in \cK: \exists \ i \in [m],\ x_j\geq 0 \}$ denote the collection of observations $k \in \cK$ which does not force $x_j$ to take value of zero and let $K_j = \vert \mathcal{K}_{j} \vert.$ 
Further, let $L_j$ denote the number of such observations containing $j$ for which $0 < p_j^k <1.$ Let $Z_j$ denote the number of such observations containing $j$ for which $ p_j^k  = 0.$ 
%Here $L_{j}=K_j$ when the observations of solutions $\{p_j^k: k \in \mathcal{K}_j\}$ are all non-zero. 
Equipped with this notation, we construct the derivatives of the perturbation functions $c_{j}^\prime (\cdot)$ for any $j\in [n]$ as follows:
\begin{enumerate}[label=(\alph*),leftmargin=*]
\item Consider any ordering $( 1, 2, \ldots,Z_j, Z_j+1, \ldots, Z_j+L_j, Z_j+ L_j + 1, \ldots, K_j)$  over the elements in $\mathcal{K}_j$ for which $\sum_{i\in [m]}\lambda_i^1 a_{ij} \leq \sum_{i\in [m]}\lambda_i^2 a_{ij}  \leq \ldots \leq \sum_{i\in [m]}\lambda_i^{Z_j} a_{ij} \leq \sum_{i\in [m]}\lambda_i^{{Z_j} +1 } a_{ij} \leq \ldots \leq \sum_{i\in [m]}\lambda_i^{Z_j + L_j} a_{ij} \leq \sum_{i\in [m]}\lambda_i^{Z_j + L_j + 1} a_{ij} \leq \ldots \leq \sum_{i\in [m]}\lambda_i^{K_j} a_{ij}.$
With  $Z_j$ defined as the size of collection $\mathcal{K}_j$ for which $p_j^k = 0,$ note that it is necessary to have $\sum_{i\in [m]}\lambda_i^{1} a_{ij} < \sum_{i\in [m]}\lambda_i^{{Z_j} +1 } a_{ij}$ whenever $Z_j >0$. Together with $L_j$ defined as the size of collection $\mathcal{K}_j$ for which $0 < p_j^k > 1,$ note that it is necessary to have $\sum_{i\in [m]}\lambda_i^{Z_j + L_j} a_{ij} < \sum_{i\in [m]}\lambda_i^{{K_j} } a_{ij}$ whenever $Z_j + L_j < K_j$. 
%With $L_j$ defined as the size of collection $\mathcal{K}_j$ for which $p_j^k > 0,$ note that it is necessary to have $\sum_{i\in [m]}\lambda_i^{L_j} a_{ij} < \sum_{i\in [m]}\lambda_i^{{L_j} +1 } a_{ij}$ whenever $L_j < K_j$.
%This follows from the observations that $\sum_{i\in [m]}\lambda_i^{L_j} a_{ij}^{L_j}$ and $\sum_{i\in [m]}\lambda_i^{{L_j} +1 } a_{ij}$ satisfies \eqref{eq:mdm-feascon-order} and $p_j^{L_j} > 0 = p_j^{L_j + 1}.$  
Further, due to the conditions in \eqref{eq:mdm-feascon-order}, the data $(p_j^k: k \in \mathcal{K}_j)$ must necessarily satisfy the ordering $p_j^1 = \ldots p_j^{Z_j} =0 <  p_j^{Z_j + 1} \leq p_j^{Z_j + 2} \leq \ldots \leq p_j^{Z_j + L_j } < 1 = p_j^{Z_j + L_j +1 } = \ldots = p_j^{K_j}$.% and $p_j^{L_j +1 } = p_j^{L_j +2 } = \ldots, p_j^{K_j} = 0.$

\item Construct the function $c_j^\prime(\cdot)$ by first setting $c_j^\prime (p_j^k) = \sum_{i\in [m]}\lambda_i^k a_{ij}$ for $k=Z_j + 1,\ldots, Z_j + L_j$. 
With this assignment, we complete the construction of the function $c_j^\prime$ in between these points by connecting them with line segments as follows: 
For any two consecutive observations $k$ and ${k+1}$ in the ordering satisfying $p_j^k < p_j^{k+1}$, connect the respective points $(p_j^k , \sum_{i\in [m]}\lambda_i^k a_{ij} )$ and $(p_j^{k+1}, \sum_{i\in [m]}\lambda_i^{k+1} a_{ij})$ with a line segment (see Figure \ref{fig:mdm-construction}). Note that if the consecutive observations $k$ and ${k+1}$ are such that $p_j^k = p_j^{k+1}$, then the corresponding points  $(p_j^k , \sum_{i\in [m]}\lambda_i^k a_{ij} )$ and $(p_j^{k+1}, \sum_{i\in [m]}\lambda_i^{k+1} a_{ij})$ coincide and there is no need to connect them. Further note that  $\sum_{i\in [m]}\lambda_i^k a_{ij} < \sum_{i\in [m]}\lambda_i^{k+1} a_{ij},$ when  $p_j^k < p_j^{k+1} $, and hence, $c_j^\prime$ is strictly increasing in the interval $(0,1)$.
% $\sum_{i\in [m]}\lambda_i^k a_{ij}  < \sum_{i\in [m]}\lambda_i^k a_{ij}^{k+1},$ connect the respective points $(\sum_{i\in [m]}\lambda_i^k a_{ij}, 1 - p_j^k)$ and $(\sum_{i\in [m]}\lambda_i^k a_{ij}, 1 - p_j^{k+1})$ with a line segment (see Figure \ref{fig:mdm-construction}). 
% For $k \leq L_j,$ note that if the consecutive observations $k$ and ${k+1}$ are such that $\sum_{i\in [m]}\lambda_i^k a_{ij}  = \sum_{i\in [m]}\lambda_i^k a_{ij},$ then the corresponding points $(\sum_{i\in [m]}\lambda_i^k a_{ij}, 1 - p_j^k)$ and $(\sum_{i\in [m]}\lambda_i^k a_{ij}^{k+1}, 1 - p_j^{k+1})$  coincide and there is no need to connect them. 
% Further note that $p_j^k > p_j^{k+1} $ when  $\sum_{i\in [m]}\lambda_i^k a_{ij} < p_j^{k+1},$ because of \eqref{eq:mdm-feascon-order}, and hence the cumulative distribution function $F_j$ is strictly increasing in the interval $[\sum_{i\in [m]}\lambda_i^1 a_{ij}, \sum_{i\in [m]}\lambda_i^{L_j} a_{ij}].$

\item Lastly we construct the remaining parts of $c_j^\prime$ as follows: if $p_j^1  = 0$ meaning $Z_j >0$, connect the points $(p_j^{Z_j}, \sum_{i\in [m]}\lambda_i^{Z_j} a_{ij})$ and $(p_j^{Z_j + 1}, \sum_{i\in [m]}\lambda_i^{Z_j + 1} a_{ij})$ with a line segment (see Figure \ref{fig:mdm-construction} (a)), if $p_j^1 >0$ meaning $Z_j = 0$, connect the points $(0, \sum_{i\in [m]}\lambda_i^1 a_{ij}) - \delta $ and $(p_j^1, \sum_{i\in [m]}\lambda_i^1 a_{ij})$ with a line segment by choosing any arbitrary $\delta > 0$ (see Figure \ref{fig:mdm-construction} (b)); if $p_j^{K_j} <1$, connect the points  $(p_j^{K_j}, \sum_{i\in [m]}\lambda_i^{K_j} a_{ij})$ and $(1, \sum_{i\in [m]}\lambda_i^{K_j} a_{ij} + \delta)$  with a line segment by choosing any arbitrary $\delta > 0$  (see Figure \ref{fig:mdm-construction} (a)), if $Z_j + L_j < K_j$, connect the points $(p_j^{Z_j + L_j}, \sum_{i\in [m]}\lambda_i^{Z_j + L_j} a_{ij}$ and $(p_j^{Z_j + L_j + 1}, \sum_{i\in [m]}\lambda_i^{Z_j + L_j +1} a_{ij})$  with a line segment (see Figure \ref{fig:mdm-construction} (b)). Hence we ensure that $c_j^\prime$ is strictly increasing in the interval $[0,1]$.
% For the right tail,  connect the points $(\sum_{i\in [m]}\lambda_i^{L_j} a_{ij}, 1 - p_j^{L_j})$ and $(\sum_{i\in [m]}\lambda_i^{L_j + 1} a_{ij},\, 1)$ with a line segment if $L_j < K_j.$ 
% We then have $F_j(x) = 1$ for every $x \geq \sum_{i\in [m]}\lambda_i^{L_j + 1} a_{ij}$ and therefore $F_j^{-1}(1) = \sum_{i\in [m]}\lambda_i^{L_j + 1} a_{ij}.$
% If $L_j = K_j,$ connect the points $(\sum_{i\in [m]}\lambda_i^{L_j} a_{ij}, 1 - p_j^{L_j})$and $(\sum_{i\in [m]}\lambda_i^{L_j} a_{ij} + \delta, \, 1)$ by choosing any arbitrary $\delta > 0$ (see Figure \ref{fig:mdm-construction}). 
% In this case, we will have $F_j(x) = 1$ for every $x \geq \sum_{i\in [m]}\lambda_i^{L_j} a_{ij} + \delta.$
% For the left tail, if $p_j^1 = 1$, then we have $F_j(x) = 0$ for every $x \leq \sum_{i\in [m]}\lambda_i^1 a_{ij}.$ 
% On the other hand,  if $p_j^1 < 1$, we use a line segment to connect 
% $(\sum_{i\in [m]}\lambda_i^1 a_{ij}, 1-p_j^1)$ and $(\sum_{i\in [m]}\lambda_i^1 a_{ij} -\delta, \,0)$ by choosing an arbitrary $\delta > 0$. 
% In this case,  $F_j(x) = 0$ for every $x \leq \sum_{i\in [m]}\lambda_i^1 a_{ij} -\delta.$
\end{enumerate}

\begin{figure}[!htbp]
\centering
\begin{tikzpicture}[scale=0.65]
  % Extended axes
  \draw[-stealth] (0,0) -- (0,12) node[right]{{$c_j^\prime(x)$}};
  \draw[-stealth] (0,0) -- (11,0) node[below]{{$x$}};

  % Origin label
  %\node[below left] at (0,0) {$0$};

  % Plot line (transposed)
  \draw[] (0,0.08) node[left]{$\underset{i\in [m]}{\sum}\lambda_i^{Z_j} a_{ij}$}
           -- (1.5,2) 
           -- (2,3.3) 
           -- (3.4,4.2) 
           -- (5.5,6.2) 
           -- (6.5,7.2) 
           -- (8.5,8.2) 
           -- (10,10.2);

  % Dashed horizontal lines
  \draw[dashed] (0,2) node[left]{$\underset{i\in [m]}{\sum}\lambda_i^{Z_j + 1} a_{ij}$} -- (1.5,2);  
  \draw[dashed] (0,4.2) node[left]{$\underset{i\in [m]}{\sum}\lambda_i^k a_{ij}$} -- (3.4,4.2);
  \draw[dashed] (0,6.2) node[left]{$\underset{i\in [m]}{\sum}\lambda_i^{k+1} a_{ij}$} -- (5.5,6.2);
  \draw[dashed] (0,8.2) node[left]{$\underset{i\in [m]}{\sum}\lambda_i^{K_j} a_{ij}$} -- (8.5,8.2);
  \draw[dashed] (0,10.2) node[left]{$\underset{i\in [m]}{\sum}\lambda_i^{K_j} a_{ij} + \delta$} -- (10,10.2);

  % Dashed vertical lines
  \draw[dashed] (-0.5,0) node[below]{$ p_j^{Z_j } = 0$};  % Keep this as-is
  \draw[dashed] (1.5,0) node[below]{$ p_j^{Z_j + 1}$} -- (1.5,2);
  \draw[dashed] (3.4,0) node[below]{$p_j^k$} -- (3.4,4.2);
  \draw[dashed] (5.5,0) node[below]{$p_j^{k+1}$} -- (5.5,6.2);
  \draw[dashed] (8.5,0) node[below]{$ p_j^{K_j}$} -- (8.5,8.2);
  \draw[dashed] (10,0) node[below]{$1$} -- (10,10.2);

  % Dots (transposed)
  \fill[black] (0,0.05) circle (2pt);
  \fill[black] (1.5,2) circle (2pt);
  \fill[black] (2,3.3) circle (2pt);
  \fill[black] (3.4,4.2) circle (2pt);
  \fill[black] (5.5,6.2) circle (2pt);
  \fill[black] (6.5,7.2) circle (2pt);
  \fill[black] (8.5,8.2) circle (2pt);
  \fill[black] (10,10.2) circle (2pt);
\end{tikzpicture}
\vspace{-0.2cm}
\begin{center}
    (a)
\end{center}  
\;
\begin{tikzpicture}[scale=0.56]
  % Axes (both 12 units)
  \draw[-stealth] (0,0) -- (0,12.5) node[right]{{$c_j^\prime(x)$}};
  \draw[-stealth] (0,0) -- (14.5,0) node[below]{{$x$}};
  \node[below left] at (0,0) {$0$}; % label at origin

  % Main line (rescaled)
  \draw[] (0.05,0) node[left]{$\underset{i\in [m]}{\sum}\lambda_i^1 a_{ij} - \delta$}
           -- (1,2)  
           -- (2,4) 
           -- (3,5) 
           -- (5,7.5) 
           -- (5.5,9) 
           -- (7.5,9.5) 
           -- (10,11.5);

  % Dots
  \fill[black] (0.05,0) circle (2pt);
  \fill[black] (1,2) circle (2pt);
  \fill[black] (2,4) circle (2pt);
  \fill[black] (3,5) circle (2pt);
  \fill[black] (5,7.5) circle (2pt);
  \fill[black] (5.5,9) circle (2pt);
  \fill[black] (7.5,9.5) circle (2pt);
  \fill[black] (10,11.5) circle (2pt);

  % Dashed vertical lines (x-axis tick values from y)
  \draw[dashed] (1,0) node[below]{$ p_j^1$} -- (1,2);
  \draw[dashed] (3,0) node[below]{$ p_j^k$} -- (3,5) ;
  \draw[dashed] (5,0) node[below]{$ p_j^{k+1}$} -- (5,7.5);
  \draw[dashed] (7.5,0) node[below]{$ p_j^{Z_j + L_j}$} -- (7.5,9.5);
  \draw[dashed] (10,0) -- (10,11.5);
  \draw[dashed] (11.8,0) node[below]{$ p_j^{Z_j + L_j + 1} = 1$}; 
  % Dashed horizontal lines (y-axis tick values from x)
  \draw[dashed] (0,2) -- (1,2);
  \draw[dashed] (0,5) -- (3,5);
  \draw[dashed] (0,7.5) -- (5,7.5);
  \draw[dashed] (0,9.5) -- (7.5,9.5);
  \draw[dashed] (0,11.5) -- (10,11.5);

  % Y-axis labels
  \node[left] at (0,2) {$\underset{i\in [m]}{\sum}\lambda_i^1 a_{ij}^1$};
  \node[left] at (0,5) {$\underset{i\in [m]}{\sum}\lambda_i^{k} a_{ij}$};
  \node[left] at (0,7.5) {$\underset{i\in [m]}{\sum}\lambda_i^{k+1} a_{ij}$};
  \node[left] at (0,9.5) {$\underset{i\in [m]}{\sum}\lambda_i^{Z_j + L_j} a_{ij}$};
  \node[left] at (0,11.5) {$\underset{i\in [m]}{\sum}\lambda_i^{Z_j + L_j + 1} a_{ij}$};
\end{tikzpicture}
\begin{center}
    (b)
\end{center}

\caption{An illustration of the construction of the functions $c_j^\prime$ when: (a)  $p_j^{K_j} <1$ and (b) $p_j^1 > 0$.} 
\label{fig:mdm-construction}
\end{figure}

The above construction gives the functions $(c_j^\prime: j \in [n])$ which are absolutely continuous and strictly increasing piecewise linear functions within $[0,1]$. There must exist $(c_j: j \in [n])$ which are strictly convex and continuously differentiable quadratic functions. 

We next show that the constructed $(c_j^\prime: j \in [n])$ yield the observed data $(p_j^k: j \in \cI^{(k)}),$ for any observation $k \in \cK,$ when they are used in the optimality conditions \eqref{s-RAM-optcon} together with the assignment $\nu_j=0,$ for $j \in [n].$ In other words, given $(p_{j}^k: j\in\cI^{(k)}, k \in \cK)$,  we next verify that  
\begin{align*}
    &  - c_j^\prime(p_j^k) + \sum_{i\in [m]} \lambda_i^k a_{ij} + \mu_j^k - w_j^k= 0,\quad \forall j\in \cI^{(k)},\  \forall k\in \cK, \\
    & \mu_j^k p_j^k = 0, \ w_j^k (1-p_j^k) = 0, \quad \forall  j \in \cI^{(k)}, \ \forall k\in \cK,\\
    & \mu_j^k, w_j^k \geq 0, \ \forall j\in \cI^{(k)},\  \forall k\in \cK,\\
    &  \sum_{j\in \cI^{(k)}} a_{ij} p_j^{(k)} + s_i - b_i   = 0, \ \alpha_i s_i = 0, \quad \forall i \in [m],\\
    & s_i \geq 0 , \alpha_i \geq 0, \quad \forall i \in [m].
\end{align*}

First, since $\bp^{(k)} \in \text{conv}(\mathcal{X}^{(k)})$, we have $\boldsymbol{A}\bp^{(k)} \leq \boldsymbol{b}$. Hence, taking $\bar{\boldsymbol{s}} = \boldsymbol{b} - \boldsymbol{A}\bp^{(k)} \geq \boldsymbol{0}$, and setting $\alpha = \boldsymbol{0}$ always satisfies $\sum_{j\in \cI^{(k)}} a_{ij} p_j^{(k)} + \bar{s}_i - b_i   = 0$, $\alpha_i \bar{s}_i = 0$. Next, we verify that under different possible values of $\bp^{(k)}$, other optimality conditions hold. 

For any $p_j^k$ with $k\in \cK, j\in \cI^{(k)}$, when  $0 < p_j^k < 1$, we have from the construction of $c_j^\prime$ that $c_j^\prime ( p_j^k ) = \sum_{i\in [m]}\lambda_i^k a_{ij} $. Then for such $p_j^k$, we see that the optimality condition $- c_j^\prime(p_j^k) + \sum_{i\in [m]} \lambda_i^k a_{ij} + \mu_j^k - w_j^k= 0 $  holds by taking $\mu_j^k = w_j^k =0$. 

For any $p_j^k$ with $k\in \cK, j\in \cI^{(k)}$, when $p_j^k =0$, we have from Steps (a) and (c) of the above construction that $\sum_{i \in [m]} \lambda_i^k a_{ij} \leq \sum_{i \in [m]} \lambda_i^{Z_j} a_{ij} = c_j^\prime(0)$. Then if we take $\mu_j^k = \sum_{i \in [m]} \lambda_i^{Z_j} a_{ij} - \sum_{i \in [m]} \lambda_i^k a_{ij}  \geq 0 $ and $w_j^k =0 $, we have $ -c_j^\prime(p_j^k) +  \sum_{i\in [m]}\lambda_i^k a_{ij} + \mu_j^k - w_j^k =0$ holds.

For any $p_j^k$ with $k\in \cK, j\in \cI^{(k)}$, when $p_j^k =1$, we have from Steps (a) and (c) of the above construction that $\sum_{i \in [m]} \lambda_i^k a_{ij} \geq \sum_{i\in [m]}\lambda_i^{Z_j+L_j+1} a_{ij} = c_j^\prime(1)$. Then if we take $\mu_j^k = 0$ and $w_j^k =\sum_{i \in [m]} \lambda_i^k a_{ij} -  \sum_{i\in [m]}\lambda_i^{Z_j+L_j+1} a_{ij} \geq 0 $, we have $ - c_j^\prime (p_j^k) +  \sum_{i\in [m]}\lambda_i^k a_{ij} + \mu_j^k - w_j^k =0$ holds.

Lastly, checking whether the conditions in \eqref{eq:mdm-feascon-order} are satisfied for given choice data $(p_{j}^k: j\in[n], k \in \cK)$ is equivalent to testing if there exists an assignment for variables $(\lambda_i^k \in \mathbb{R}: i\in [m], k \in \cK)$ and $\epsilon > 0$ such that,
\begin{align*}
    &\sum_{i\in [m]}\lambda_i^k a_{ij} + \epsilon   \leq \sum_{i\in [m]}\lambda_i^l a_{ij}  \quad \text{if} \quad p_{j}^k < p_{j}^l, \quad  \forall l,k \in \cK, l\neq k,\ \forall j\in \cI^{(k)}\cap \cI^{(l)}, \\
    &\sum_{i\in [m]}\lambda_i^k a_{ij}   = \sum_{i\in [m]}\lambda_i^l a_{ij}  \quad \text{if} \quad 0< p_{j}^k = p_{j}^l <1, \quad \forall l,k \in \cK, l\neq k,\ \forall j\in \cI^{(k)}\cap \cI^{(l)}.
\end{align*}
This is possible in polynomial time by solving a linear program where the above conditions are formulated as constraints and maximizing $\epsilon.$  This linear program involves $mK$ variables for $(\lambda_i^k \in \mathbb{R}: i\in [m], k \in \cK)$ and one variable for $\epsilon,$ and $\mathcal{O}((n+m)K)$ constraints. 
\end{proof}

%\subsection{Proof of Proposition \ref{prop:wc-obj-reform}}

\begin{proof}{Proposition \ref{prop:wc-obj-reform}}
Here, as before, $\ps-RAM(\cK^\prime)$ denotes the collection of all s-RAM representable probabilities over the collection $\cI^{\cK^\prime},$ where $\cI^{\cK^\prime} = \cI^\cK \cup \{\cI^\text{new}\}$.
Due to Theorem \ref{thm:mdm-general-feascon}, we have  $\ps-RAM(\cK^\prime) = \text{Proj}_{\bx}(\Pi_{\cK^\prime}),$ where the lifted set $\Pi _{\cK^\prime}$ equals 
\begin{align*}
    \Pi _{\cK^\prime} = \Big\{  (\bx,\boldsymbol{\lambda}) :  \,  &\sum_{j\in \cI^{(k)}} a_{ij} x_j \leq b_i, \ \forall k \in \cK^\prime, \forall  i \in [m],  0 \leq x_j \leq 1, \ \forall k \in \cK^\prime,  \forall  j \in \cI^{(k)},   \\ 
    &\sum_{i\in [m]}\lambda_i^k a_{ij}   < \sum_{i\in [m]}\lambda_i^l a_{ij}  \; \text{if} \; x_{j}^k < x_{j}^l, \quad \forall k,l \in \cK^\prime, \forall j\in \cI^{(k)} \cap \cI^{(l)},  \\ 
    & \sum_{i\in [m]}\lambda_i^k a_{ij}   = \sum_{i\in [m]}\lambda_i^l a_{ij}  \; \text{if} \; 0< x_{j}^k = x_{j}^l <1, \quad \forall k,l \in \cK^\prime, \forall j\in \cI^{(k)} \cap \cI^{(l)}
    \Big\}.
\end{align*}
Since $\mathcal{U}_{\text{new}}:= \left\{ \bx: (\bx,\bx^\cK,\blam)  \in \Pi _{\cK^\prime}, \bx^\cK = \bp^\cK  \right\}$, the non-numbered constraints in the formulation in Proposition \ref{prop:wc-obj-reform} are obtained from the definition of $\mathcal{X}^{\text{new}}$ and by replacing $ \bx^\cK = \bp^\cK$ in the description of the lifted set $\Pi _{\cK^\prime}$. 

For deducing the  remaining constraints  \eqref{model:mdm_prediction-a}  - \eqref{model:mdm_prediction-b},  we proceed as follows: Consider any $k \in \cK$ and $ j \in \cI^{(k)}$. From the description of $\Pi_{\cK^\prime},$  observe that an assignment for $x_j, \{\lambda_i^{\text{new}}:i\in [m]\},  \{\lambda_i^k: i\in [m] \}$ in any $(\bx,\blambda) \in \Pi_\cK^\prime$  satisfying $\bx^\cK = \bp^\cK$ necessarily satisfies one of the following four cases: 

In \textit{Case 1,} we have $\sum_{i\in [m]}\lambda_i^\text{new} a_{ij}  < \sum_{i\in [m]}\lambda_i^k a_{ij}$ and  $x_{j} < p_{j}^k $: 
If $\sum_{i\in [m]}\lambda_i^\text{new} a_{ij}$, $\sum_{i\in [m]}\lambda_i^k a_{ij}$ are such that  $\sum_{i\in [m]}\lambda_i^\text{new} a_{ij}  < \sum_{i\in [m]}\lambda_i^k a_{ij}$,  this informs the restriction $\{x_{j}: x_{j} < p_{j}^k \}$ on the values $x_{j}$ can take. The closure of this  restricted collection $\{x_{j}: x_{j} < p_{j}^k\}$ equals $\{x_{j}: x_{j} \leq p_{j}^k\}.$ 

In \textit{Case 2,} we have $\sum_{i\in [m]}\lambda_i^\text{new} a_{ij}  > \sum_{i\in [m]}\lambda_i^k a_{ij}$  and  $x_{j} > p_{j}^k $: 
If $\sum_{i\in [m]}\lambda_i^\text{new} a_{ij}$, $\sum_{i\in [m]}\lambda_i^k a_{ij}$ are such that  $\sum_{i\in [m]}\lambda_i^\text{new} a_{ij}  > \sum_{i\in [m]}\lambda_i^k a_{ij}$, the closure of the corresponding restriction $\{x_{j}: x_{j} > p_{j}^k \}$ equals $\{x_{j}: x_{j} \geq p_{j}^k\}.$ 

In \textit{Case 3,} we have $\sum_{i\in [m]}\lambda_i^\text{new} a_{ij}  = \sum_{i\in [m]}\lambda_i^k a_{ij}$  and  $x_{j} = p_{j}^k \neq 0 $:
When $\sum_{i\in [m]}\lambda_i^\text{new} a_{ij}$, $\sum_{i\in [m]}\lambda_i^k a_{ij}$ are such that  $\sum_{i\in [m]}\lambda_i^\text{new} a_{ij}  = \sum_{i\in [m]}\lambda_i^k a_{ij}$ and $p_{j}^k \neq 0,$ the corresponding restriction on the values of $x_{j}$ is given by the closed set  $\{x_{j}: x_{j} = p_{j}^k\}.$ 

Finally, in \textit{Case 4,} we have $\sum_{i\in [m]}\lambda_i^\text{new} a_{ij} - \sum_{i\in [m]}\lambda_i^k a_{ij}$ unconstrained and  $x_{j} = p_{j}^k = 0$: Like in Case 3, the  restriction on the values of $x_{j}$ corresponding to this case equals  $\{x_{j}: x_{j} = 0\}.$ The relationship between $x_{j}$, $ p_{j}^k$, $\sum_{i\in [m]}\lambda_i^\text{new} a_{ij}$, $\sum_{i\in [m]}\lambda_i^k a_{ij}$ in this case is any one of the following sub-cases: Case (4a) $\sum_{i\in [m]}\lambda_i^\text{new} a_{ij}  > \sum_{i\in [m]}\lambda_i^k a_{ij}$ and $0 = x_{j} \leq p_{j}^k = 0,$ or Case (4b) $\sum_{i\in [m]}\lambda_i^\text{new} a_{ij}  < \sum_{i\in [m]}\lambda_i^k a_{ij}$ and $0 = x_{j} \geq p_{j}^k = 0,$ or Case (4c) $\sum_{i\in [m]}\lambda_i^\text{new} a_{ij}  = \sum_{i\in [m]}\lambda_i^k a_{ij}$ and $0 = x_{j} = p_{j}^k = 0.$ Similar arguments can be naturally extended to the case that we have $\sum_{i\in [m]}\lambda_i^\text{new} a_{ij} - \sum_{i\in [m]}\lambda_i^k a_{ij}$ unconstrained and  $x_{j} = p_{j}^k = 1.$

Combining the observations in the cases (1) \& (4a), (2) \& (4b), and (3) \& (4c), we obtain that the closure of $\mathcal{U}_{\text{new}}$ equals the collection of probability vectors $\bx = (x_j: j \in \cI^{\text{new}})$ for which there exists a function $\blam$ such that 
\begin{align*}
   & x_{j}  \  \leq \   p_{j}^k \quad \text{if} \quad  \sum_{i\in [m]}\lambda_i^\text{new} a_{ij}  <  \sum_{i\in [m]}\lambda_i^k a_{ij} , \quad \forall k \in \cK, \forall j\in \cI^\text{new} \cap \cI^{(k)},\\     
    & x_{j}  \  \geq \   p_{j}^k \quad \text{if} \quad  \sum_{i\in [m]}\lambda_i^\text{new} a_{ij}  > \sum_{i\in [m]}\lambda_i^k a_{ij}, \quad  \forall k \in \cK, \forall j\in \cI^\text{new} \cap \cI^{(k)}.
\end{align*}
%Additionally, $\blam$ satisfies %
% in addition to satisfying 
% $ \sum_{i\in [m]}\lambda_i^k a_{ij}   > \sum_{i\in [m]}\lambda_i^l a_{ij}$ if $p_{j}^k < p_{j}^l$,
% $ \sum_{i\in [m]}\lambda_i^k a_{ij}   = \sum_{i\in [m]}\lambda_i^l a_{ij} $ if $ p_{j}^k = p_{j}^l \neq 0$, for all $ k ,l \in \cK,$ for all $ j\in \cI^{(k)} \cap \cI^{(l)},$ and  $\lambda_i^k = 0 $ if $ b_i - \sum_{j\in \cI^{(k)}} a_{ij} p_j^k  > 0$, $\lambda_i^k \geq 0.$ for all $ k \in \cK$, for all $ i \in [m].$
% The constraints in the formulation in Proposition \ref{prop:wc-obj-reform} exactly specify these conditions describing the closure of $\mathcal{U}_{\text{new}}.$ 

Observe that in the objective $f(\bx)$ is continuous as a function of $\bx,$ and therefore, $\inf \{ f(\bx): \bx \in \mathcal{U}_{\text{new}} \} = \min\{ f(\bx) : \bx \in \mathcal{U}_{\text{new}}\}.$
\end{proof}

%\subsection{Proof of Proposition \ref{prop:wc-obj-milp} }

\begin{proof}{Proposition \ref{prop:wc-obj-milp}}
Recall the notation $\cI^{\cK^\prime} = \cI^\cK \cup \{\cI^\text{new}\}$.  
Observe that the variables $(\lambda_j^k : j\in \cI^{(k)}, k \in \cK^\prime)$ influence the value of the formulation in Proposition \ref{prop:wc-obj-reform} only via the sign of $\sum_{i\in [m]}\lambda_i^k a_{ij} - \sum_{i\in [m]}\lambda_i^l a_{ij},$ for any pair of variables $\sum_{i\in [m]}\lambda_i^k a_{ij}$,$\sum_{i\in [m]}\lambda_i^l a_{ij}$. 
For example, the optimal value of this optimization formulation is not affected by the presence of the following additional constraints: 
 $-1 \leq\sum_{i\in [m]}\lambda_i^k a_{ij} \leq 1$  for all $k \in \cK^\prime,$ for all $j\in \cI^{(k)}$ and 
\begin{align*}
    \sum_{i\in [m]}\lambda_i^k a_{ij} - \sum_{i\in [m]}\lambda_i^l a_{ij}  \geq \epsilon \quad \text{ if } \quad \sum_{i\in [m]}\lambda_i^k a_{ij} > \sum_{i\in [m]}\lambda_i^l a_{ij}, \quad  \forall k ,l \in \cK^\prime, \forall j\in \cI^{(k)} \cap \cI^{(l)},
\end{align*} 
for some suitably small value of $\epsilon > 0.$ Indeed, this is because the signs of the differences $\{ \sum_{i\in [m]}\lambda_i^k a_{ij} - \sum_{i\in [m]}\lambda_i^l a_{ij} : \forall k ,l \in \cK^\prime, \forall j\in \cI^{(k)} \cap \cI^{(l)}\}$ are not affected by these additional constraints. Taking  $\epsilon$ to be smaller than $1/(2nK),$ for example, ensures that there is a feasible assignment for  $(\sum_{i\in [m]} \lambda_i^k a_{ij}: j \in \cI^{(k)}, k \in  \cK^\prime)$ within the interval $[-1,1]$ even if all these variables take distinct values. Therefore, taking  $\epsilon$ to be smaller than $1/(2nK)$ ensures the problem is feasible. 

Let $F$ denote the feasible values for the variables $(\lambda_j^k : j\in \cI^{(k)}, k \in \cK^\prime)$ and $(x_{j}: j \in \cI^{\text{new}})$ satisfying the constraints introduced in the above paragraph besides those in the formulation in Proposition \ref{prop:wc-obj-reform}. 
 %Therefore, if we let $F$ denote the feasible values for the variables  $(\lambda_S: k \in \cK^\prime), (x_{i,A}: i \in A)$ satisfying these constraints in addition to those in the formulation in Proposition \ref{prop:wc-obj-reform}, the value of the resulting minimization is unaffected.  
Equipped with this feasible region $F$, we have the following deductions from \eqref{mdm:milp-a} - \eqref{mdm:milp-c} for $(\lambda_j^k : j\in \cI^{(k)}, k \in \cK^\prime)$ and $(x_{j}: j \in \cI^{\text{new}})$ in $F$: For every $j \in \cI^{\text{new}}$ and any $k \in \cK$ such that $\cI^{(k)}$ containing $j,$   
\begin{itemize}[leftmargin=*]
    \item[(i)] we have $ \sum_{i\in [m]}\lambda_i^\text{new} a_{ij} < \sum_{i\in [m]}\lambda_i^k a_{ij}$ if and only if  $\delta_j^{\text{new},k} = 1$ and $\delta_j^{k,\text{new}} = 0,$ due to the constraints \eqref{mdm:milp-a} and \eqref{mdm:milp-b}; in this case, we have from \eqref{mdm:milp-c} that $0 \leq x_{j} \leq  p_{j}^k ;$
     \item[(ii)] likewise, we have $ \sum_{i\in [m]}\lambda_i^\text{new} a_{ij} > \sum_{i\in [m]}\lambda_i^k a_{ij}$ if and only if  $\delta_j^{\text{new},k} = 0$ and $\delta_j^{k,\text{new}} = 1,$ due to the constraints \eqref{mdm:milp-a} and \eqref{mdm:milp-b}; in this case, we have from \eqref{mdm:milp-c} that $ p_{j}^k \leq  x_{j} \leq 1.$
     \item[(iii)] finally, $ \sum_{i\in [m]}\lambda_i^\text{new} a_{ij} = \sum_{i\in [m]}\lambda_i^k a_{ij}$ if and only if  $\delta_j^{\text{new},k} = 0$ and $\delta_j^{k,\text{new}} = 0;$ here we have from \eqref{mdm:milp-d} that $x_{j} = p_{j}^k.$
\end{itemize} 
Thus the binary variables $\{\delta_j^{\text{new},k},\delta_j^{k,\text{new}}: k \in \cK,\forall j \in \cI^\text{new} \cap \cI^{(k)} \}$ suitably model the  constraints collection \eqref{model:mdm_prediction-a} - \eqref{model:mdm_prediction-b} and provide an equivalent reformulation in terms of the constraints \eqref{mdm:milp-a} - \eqref{mdm:milp-d}. 

Therefore, the optimal value of the formulations in Propositions \ref{prop:wc-obj-reform} and \ref{prop:wc-obj-milp} are identical.  
\end{proof}

%\subsection{Proof of Proposition \ref{prop:limit-mdm}}
\begin{proof}{Proposition \ref{prop:limit-mdm}}
Due to Theorem \ref{thm:mdm-general-feascon}, we have  $\ps-RAM(\cK) = \text{Proj}_{\bx}(\Pi_\cK),$ following the definition in  \eqref{eq:lifted-set}. One can argue the closure of $\ps-RAM(\cK)$ by the similar arguments of the proof of Proposition \ref{prop:wc-obj-reform}. 

Consider any $k \in \cK$ and $ j \in \cI^{(k)}$. From the description of $\Pi_{\cK},$  observe that an assignment for $x_j^{k}, x_j^l$ and $\sum_{i\in [m]}\lambda_i^k a_{ij}$,$\sum_{i\in [m]}\lambda_i^l a_{ij}$ in any $(\bx,\blambda) \in \Pi_\cK$ necessarily satisfies one of the following four cases: 

In \textit{Case 1,} we have $\sum_{i\in [m]}\lambda_i^k a_{ij}  < \sum_{i\in [m]}\lambda_i^l a_{ij}$ and  $x_{j}^k < x_{j}^l$: If $\sum_{i\in [m]}\lambda_i^k a_{ij}$,$\sum_{i\in [m]}\lambda_i^l a_{ij}$ is such that  $\sum_{i\in [m]}\lambda_i^k a_{ij}  < \sum_{i\in [m]}\lambda_i^l a_{ij}$, the closure of the corresponding restriction $\{x_{j}^k , x_{j}^l: x_{j}^k < x_{j}^l\}$ equals $\{x_{j}^k , x_{j}^l: x_{j}^k \leq x_{j}^l\}.$
 
In \textit{Case 2,} we have $\sum_{i\in [m]}\lambda_i^k a_{ij}  > \sum_{i\in [m]}\lambda_i^l a_{ij}$ and $x_{j}^k > x_{j}^l$: If $\sum_{i\in [m]}\lambda_i^k a_{ij}$,$\sum_{i\in [m]}\lambda_i^l a_{ij}$ is such that $\sum_{i\in [m]}\lambda_i^k a_{ij}  > \sum_{i\in [m]}\lambda_i^l a_{ij}$, the closure of the corresponding restriction $\{x_{j}^k , x_{j}^l: x_{j}^k > x_{j}^l\}$ equals $\{x_{j}^k , x_{j}^l: x_{j}^k \geq x_{j}^l\}.$

In \textit{Case 3,} we have $\sum_{i\in [m]}\lambda_i^k a_{ij}  = \sum_{i\in [m]}\lambda_i^l a_{ij}$ and $x_{j}^k = x_{j}^l >0$: When $\sum_{i\in [m]}\lambda_i^k a_{ij}$,$\sum_{i\in [m]}\lambda_i^l a_{ij}$ is such that  $\sum_{i\in [m]}\lambda_i^k a_{ij}  = \sum_{i\in [m]}\lambda_i^l a_{ij}$, the corresponding restriction on the values of $x_{j}^k , x_{j}^l$ is given by the closed set $\{x_{j}^k , x_{j}^l: 0 < x_{j}^k = x_{j}^l < 1\}.$

Finally, in \textit{Case 4,} we have $\sum_{i\in [m]}\lambda_i^k a_{ij}$,$\sum_{i\in [m]}\lambda_i^l a_{ij}$ unconstrained and $x_{j}^k = x_{j}^l = 0$: in Case 3, the  restriction on the values of $x_{j}^k , x_{j}^l$ corresponding to this case equals $x_{j}^k = x_{j}^l = 0$. The relationship between $x_{j}^k , x_{j}^l$ and $\sum_{i\in [m]}\lambda_i^k a_{ij}$,$\sum_{i\in [m]}\lambda_i^l a_{ij}$ in this case is any one of the following sub-cases: Case (4a) 
$\sum_{i\in [m]}\lambda_i^k a_{ij}  > \sum_{i\in [m]}\lambda_i^l a_{ij}$ and $0 = x_{j}^k \leq x_{j}^l = 0,$ 
or Case (4b) $\sum_{i\in [m]}\lambda_i^k a_{ij}  < \sum_{i\in [m]}\lambda_i^l a_{ij}$ and $0 = x_{j}^k \geq x_{j}^l = 0,$  or Case (4c) $\sum_{i\in [m]}\lambda_i^k a_{ij}  = \sum_{i\in [m]}\lambda_i^l a_{ij}$ and $0 = x_{j}^k = x_{j}^l = 0.$ Similar arguments can be naturally extended to the case that we have $\sum_{i\in [m]}\lambda_i^k a_{ij}$,$\sum_{i\in [m]}\lambda_i^l a_{ij}$ unconstrained and $x_{j}^k = x_{j}^l = 1.$

Combining the observations in the cases (1) \& (4a), (2) \& (4b), and (3) \& (4c), we obtain that the closure of $\ps-RAM(\cK)$ equals the collection of probability vectors $\bx$ for which there exists  $\blam = \{\lambda_i^k: i\in [m], k\in \cK\}$ such that 
\begin{align*}
   &x_{j}^k \ \leq \ x_{j}^l \quad\text{ if } \quad \sum_{i\in [m]}\lambda_i^k a_{ij}  < \sum_{i\in [m]}\lambda_i^l a_{ij}, \quad \forall k,l\in \cK, \forall j\in \cI^{(k)}\cap \cI^{(l)},\\
    &x_{j}^k \ \geq \ x_{j}^l \quad\text{ if } \quad \sum_{i\in [m]}\lambda_i^k a_{ij}  > \sum_{i\in [m]}\lambda_i^l a_{ij}, \quad \forall k,l\in \cK, \forall j\in \cI^{(k)}\cap \cI^{(l)},\\
    &x_{j}^k \ = \ x_{j}^l \quad\text{ if } \quad \sum_{i\in [m]}\lambda_i^k a_{ij}  = \sum_{i\in [m]}\lambda_i^l a_{ij}, \quad \forall k,l\in \cK, \forall j\in \cI^{(k)}\cap \cI^{(l)}.
\end{align*}
The constraints in the formulation in Proposition \ref{prop:limit-mdm} exactly specify these conditions describing the closure of $\ps-RAM(\cK).$ Observe that the objective in $\inf \{ \text{loss}(\bp^\cK,\bx^\cK)\,:\, \bx^{\cK} \in \ps-RAM(\cK) \}$ is continuous as a function of $\bx.$ Therefore, $\inf \{ \text{loss}(\bp^\cK,\bx^\cK)\,:\, \bx^{\cK} \in \ps-RAM(\cK) \} = \min \{ \text{loss}(\bp^\cK,\bx^\cK)\,:\, \bx^{\cK} \in \text{closure}(\ps-RAM(\cK))\}.$
\end{proof}

%\subsection{Proof of Proposition \ref{prop:micp}}
\begin{proof}{Proposition \ref{prop:micp}}
Observe that the variables $(\lambda_i^k: k \in \cK, i \in [m])$ influence the value of the formulation in Proposition \ref{prop:micp} only via the sign of $\sum_{i\in [m]}\lambda_i^k a_{ij} - \sum_{i\in [m]}\lambda_i^l a_{ij},$ for any pair of variables $\sum_{i\in [m]}\lambda_i^k a_{ij}$,$\sum_{i\in [m]}\lambda_i^l a_{ij}$ from the collection $(\lambda_i^k: k\in \cK, i \in [m])$.  Therefore the optimal value of this optimization formulation is not affected by the presence of the following additional constraints:  $-1 \leq \sum_{i\in [m]}\lambda_i^k a_{ij} - \sum_{i\in [m]}\lambda_i^l a_{ij} \leq 1$  for all $k \in \cK,$ for all $j\in \cI^{(k)}$. 
Indeed, this is because the signs of the differences $\{ \sum_{i\in [m]}\lambda_i^k a_{ij} - \sum_{i\in [m]}\lambda_i^l a_{ij}: k,l \in \cK, j\in \cI^{(k)} \cap \cI^{(l)} \}$ are not affected by these additional constraints. Taking $\epsilon$ to be smaller than $1/(2nK),$ for example, ensures that there is a feasible assignment for  $(\sum_{i\in [m]}\lambda_i^k a_{ij} - \sum_{i\in [m]}\lambda_i^l a_{ij} : k,l \in \cK, j\in \cI^{(k)}\cap \cI^{(l)} )$ within the interval $[-1,1]$ even if all these variables take distinct values. 

Let $F$ denote the feasible values for the variables $(\lambda_i^k: k \in \cK, i \in [m])$ and $(x_j^k: k \in \cK, j \in \cI^{(k)})$
satisfying the constraints introduced in the above paragraph besides those in the formulation in Proposition \ref{prop:limit-mdm}.  
Equipped with this feasible region $F$, we have the following deductions for $(\lambda_i^k: k \in \cK, i \in [m])$ and $(x_j^k: k \in \cK, j \in \cI^{(k)})$ in $F$: For any $k,l \in \cK$ such that $j\in \cI^{(k)} \cap \cI^{(l)}$,

\begin{itemize}[leftmargin=*]
    \item[(i)] we have $ \sum_{i\in [m]}\lambda_i^k a_{ij} < \sum_{i\in [m]}\lambda_i^l a_{ij}$ if and only if  $\delta_j^{k,l} = 1$ and $\delta_j^{l,k} = 0,$
    due to the first set of constraints of \eqref{model:mdm_micp}; in this case, from the second set of constraints of \eqref{model:mdm_micp} and the constraint that $0 \leq x_j^k \leq 1$, for all $k\in \cK$ and $j\in \cI^{(k)}$, we have
    that $0 \leq x_{j}^k \leq x_{j}^l  \leq 1;$
     \item[(ii)] likewise, we have $ \sum_{i\in [m]}\lambda_i^k a_{ij} > \sum_{i\in [m]}\lambda_i^l a_{ij}$
     if and only if  $\delta_j^{k,l} = 0$ and $\delta_j^{l,k} = 1,$ due to the first set of constraints; in this case, from the second set of constraints of \eqref{model:mdm_micp} and the constraint that $0 \leq x_j^k \leq 1$, for all $k\in \cK$ and $j\in \cI^{(k)}$, 
     we have that $0 \leq x_{j}^l \leq x_{j}^k  \leq 1;$
     \item[(iii)] finally, $ \sum_{i\in [m]}\lambda_i^k a_{ij} = \sum_{i\in [m]}\lambda_i^l a_{ij}$ if and only if  $\delta_j^{k,l} = 0$ and $\delta_j^{l,k} = 0;$ here we have from the third set of constraints of \eqref{model:mdm_micp} that $ x_{j}^k = x_{j}^l .$
\end{itemize} 
Thus, the binary variables $\{\delta_j^{k,l}: k,l \in \cK,  \forall j\in \cI^{(l)} \cap \cI^{(k)}\}$ suitably model the first set of constraints of \eqref{model:lom} and provide an equivalent reformulation in terms of the constraints.

% To deduce the second set of constraints of \eqref{model:lom}, for all $k \in \cK$ and for all $i\in [m]$, we have 
% \begin{itemize}[leftmargin=*]
%     \item[(i)] $b_i - \sum_{j\in \cI^{(k)} } a_{ij} x_j^k  > 0$ if and only if $z_i^k = 1$ due to the fifth constraint of \eqref{model:mdm_micp} ; in this case, we have from the sixth constriant of \eqref{model:mdm_micp} that $\lambda_i^k \leq 0$. We also have $ \lambda_i^k \geq 0 $. This implies $\lambda_i^k = 0.$
%     \item[(ii)] $b_i - \sum_{j\in \cI^{(k)} } a_{ij} x_j^k  = 0$ if and only if $z_i^k = 0$ due to the forth constraint of \eqref{model:mdm_micp}; in this case, we have from the fifth constraint that $b_i - \sum_{j\in \cI^{(k)} } a_{ij} x_j^k \leq M $ from the sixth constraint that $\lambda_i^k \leq 1.$ The feasibility of \eqref{model:mdm_micp} is not affected.
% \end{itemize}
% Thus, the binary variables $\{ z_i^k: k\in \cK, i \in [m]\}$ suitably model the second set of constraints in \eqref{model:lom}  and provide an equivalent reformulation in terms of the constraints.

Therefore, the optimal value of the formulations in Propositions \ref{prop:limit-mdm} and \ref{prop:micp} are identical.  
\end{proof}

\section{Missing Details on Numerical Experiments in Section \ref{sec:experiment}} \label{app:exp_examples}

For reproducibility, we report the computational environment used in all experiments. Computations were performed on a MacBook Pro with a 12-core Apple M4 Pro chip and 24 GB RAM. The code was executed in Python 3.13.5 with NumPy 2.1.3, SciPy 1.15.3, pandas 2.2.3, and CVXPY 1.7.2. The optimization problems were solved using  Gurobi 12.0.3 and MOSEK 11.1.2.

\subsection{Data Generation}\label{app:exp_data}
\textbf{Gaussian Data for representability and estimation test.} 
We provide a detailed description of the data generation procedure used in our experiments for evaluating the representational power of s-RAM across four combinatorial problems: longest path, shortest path, assignment, and constrained subset selection.

For the path problems, we generate directed acyclic graphs (DAGs) with 10 nodes denoted as $\bV$ with $|\bV| = 10$. For the assignment problem, we generate bipartite graphs with 4 nodes on the left denoted as $\bV_L$, and 6 nodes on the right denoted as $\bV_R$, where $\bV = \bV_L \cup  \bV_R$ is the set of nodes of the whole set of nodes and $|\bV| =10$. For each of \( K \in \{5, 10, 15, 20, 25, 30\} \), we construct \(K\) historical graphs $\bG^1 = (\bV^1, \bE^1), \bG^2 =(\bV^2, \bE^2),\cdots, \bG^K = (\bV^K, \bE^K)$. For $k\in \cK$, we generate $\bE^k$ as follows. Edges in each graph are included independently with a probability drawn uniformly from \([0.4, 0.6]\). To ensure source-to-sink connectivity in the DAGs, we first randomly generate a path from the source node to the sink node and include all edges along this path. Additional edges are then added independently according to the sampled probability. All nodes being connected are considered as $\bV^k$ in the path problems. For the assignment problem, we ensure that each node has at least one incident edge. Therefore, $\bV^k = \bV$ in the assignment problem.
On average, the DAGs contain around 20 edges, and the bipartite graphs contain around 12 edges.

To generate data outside the s-RAM family, for each $e\in \bE$, we assign edge costs of the form $v_e + \epsilon_e$, where \(v_e\) is drawn independently from a uniform distribution over \([1, 7]\), and \(\epsilon_e\) is a Gaussian noise term. We consider three types of noise structures. In the independent setting, each \(\epsilon_e\) is drawn independently from a Gaussian distribution with zero mean and variance sampled uniformly from \([9, 25]\). In the negatively and positively correlated settings, we construct a covariance matrix with correlation coefficients sampled uniformly from \([-0.9, -0.5]\) or \([0.5, 0.9]\), respectively, and variances sampled uniformly from \([9, 25]\). The noise vector is then sampled from the corresponding multivariate Gaussian distribution.

For every collection of graphs, we generate an instance $\bp^\cK$ as follows. Given a graph $\bG^k$,  we fix the deterministic cost vector \(v_e\) and generate 100 realizations of the noise vector \(\epsilon_e\). For each realization, we solve the deterministic version of the corresponding combinatorial problem using the perturbed edge costs \(v_e + \epsilon_e\). The resulting 100 solutions are averaged to form the observed choice data for that graph. This procedure is repeated across all \(K\) historical graphs, resulting in \(K\) sets of observed choice data for each instance.  For each combination of \(K\) and correlation structure, we generate 20 such instances. To ensure consistency across noise structures, we reuse the same base graphs.

% In the constrained subset selection problem, for every collection of subsets $\{\cI^{(k):k\in \cK}\}$, we generate an instance $\bp^\cK$ as follows.  The generation of $\bp^{(k)}$ for the constrained subset selection problem follows the same procedure as in the longest path experiments, with item costs given by $v_i + \epsilon_i$ under the same noise structures and parameter settings. The only difference is that the deterministic constrained problem is solved to obtain the observed data.  

In the constrained subset selection problem, for each collection of available subsets $\{\cI^{(k):k\in \cK}\}$, we generate an instance $\bp^\cK$ as follows.
For each $k\in \cK$, the marginal allocation weights $bp^{(k)}$  are generated using the same procedure as in the longest path experiments. Specifically, item costs take the form $v_i + \epsilon_i$, where the deterministic components $v_i$ and noise terms $\epsilon_i$  follow the same distributions, noise structures, and parameter settings as before.
The only difference is that, for each realization, we solve the corresponding deterministic constrained subset selection problem to obtain the observed marginal allocation weights.

For an instance $\bp^\cK$, we verify whether the observed data $\bp^\cK$ is representable by s-RAM by checking the conditions in equation~\eqref{eq:mdm-feascon-order}. If the condition is violated, we compute the average absolute deviation loss using equation~\eqref{model:mdm_micp}, where the loss function is defined as $\frac{1}{K}\sum_{k \in \cK} \sum_{e \in \bE^{k}} |x_e^k - p_e^k|.$ Here, $p_e^k$ is the observed flow on edge $e$ in observation $k$, and $x_e^k$  is the probability estimated under the s-RAM model or for the constrained subset selection problem, $\frac{1}{K}\sum_{k \in \cK} \sum_{i \in \cI^{(k)}} |x_i^k - p_i^k|.$ Here, \(p_i^k\) is the observed marginal weight of item  $i$ in observation $k$, and $x_i^k$ is the probability estimated under the s-RAM model . 
Because the data are generated randomly, enforcing exact representability can be overly restrictive: instances with vanishingly small loss may be classified as non-representable due to numerical or sampling noise. 

Accordingly, we treat an instance as s-RAM representable if its average absolute deviation loss is below $10^{-6}$, and report the representable fraction under this criterion.
Finally, we report, for each \(K\) and noise setting, the proportion of representable instances and the average absolute deviation loss across the 20 generated instances.

\textbf{Parametric s-RAM data for prediction.}
We describe the data generation process used to evaluate the predictive accuracy of the nonparametric s-RAM framework on counterfactual tasks. In these experiments, observed data is generated from parametric s-RAM models, using their equivalent MDM formulations, and the goal is to assess whether the nonparametric s-RAM can accurately predict outcomes under new graph structures across different combinatorial problems, different underlying s-RAMs, and different counterfactual prediction tasks.

To retain interpretability grounded in utility theory, and given the equivalence between s-RAM and MDM, we generate data from parametric MDMs with structured assumptions on the marginal distributions. We still consider longest path problems, shortest path problems, and assignment problems. In each case, the evaluation focuses on whether s-RAM can accurately predict task-specific counterfactual quantities.

\textbf{Longest Path Problem.} Recall that for the longest path problem, we assume the noise marginals are nonidentical exponential distributions. Given a collection historical graphs $\bG^1 = (\bV^1, \bE^1), \bG^2 =(\bV^2, \bE^2),\cdots, \bG^K = (\bV^K, \bE^K)$, for a graph $\bG^k$ with $k\in \cK$, we solve a convex program defined by the corresponding parametric MDM model to obtain the choice probabilities. Let $\bv$ denote the deterministic weight matrix. Let $\lambda_{ij}$ be the rate parameters for each edge $(i,j)$ in the set of all edges $\bE$. From \citep{natarajan2009persistency}, we have that, for a graph $\bG^k$, the convex program with  nonidentical exponential marginal distributions for the longest path problem is as follows:
\begin{align}
\begin{aligned}
    \max_{\bx} & \sum_{(i,j) \in \bE^k}(\nu_{ij}+\frac{1}{\lambda_{ij}})x_{ij} - \sum_{(i,j) \in \bE^k} \frac{1}{\lambda_{ij}} x_{ij}\ln{x_{ij}}  \\ 
   \text{s.t.}\; & \sum_{j:(i,j)\in \bE^k} x_{ij} - \sum_{j:(j,i)\in \bE^k} x_{ji} = \begin{cases}
         1,\quad &i=s, \\
         -1,\quad &i=t, \\
         0,\quad &i\in \bV^k,
     \end{cases}\\
    &   x_{ij}\geq 0 \quad \forall \ (i,j)\in \bE^k.    
\end{aligned}\label{data:longestpath}
\end{align}
These choice probabilities are used as the observed data $\bp^\cK$ to train the nonparametric s-RAM model. For prediction, we generate a new graph $\bG^{\text{new}}$ and compute the probabilities via \ref{data:longestpath} with $\bG^{\text{new}}$.

%evaluate whether the model can predict the expected flow on the pre-specified key edge, which remains influential in the new graph.

\textbf{Shortest Path Problem.} Recall that for the longest path problem, we assume the noise marginals are nonidentical uniform distributions. Let $\bv$ denote the deterministic weight matrix, and let $\mathbf{a}$ and $\mathbf{b}$ denote the matrices of lower and upper bounds of the uniform noise distributions, respectively.  From \citep{natarajan2009persistency}, we have that, for a graph $\bG^k$ with $k\in \cK$, the convex program with  nonidentical uniform marginal distributions for the shortest path problem is as follows:
\begin{align}
\begin{aligned}
    \min_{\bx} & \sum_{(i,j) \in \bE^k}(\nu_{ij}+ a_{ij})x_{ij} + \sum_{(i,j) \in \bE^k} \frac{b_{ij} - a_{ij}}{2} x_{ij}^2  \\ 
   \text{s.t.}\; & \sum_{j:(i,j)\in \bE^k} x_{ij} - \sum_{j:(j,i)\in \bE^k} x_{ji} = \begin{cases}
         1,\quad &i=s, \\
         -1,\quad &i=t, \\
         0,\quad &i\in \bV^k,
     \end{cases}\\
    &   x_{ij}\geq 0 \quad \forall \ (i,j)\in \bE^k.    
\end{aligned}\label{data:shortestpath}
\end{align}
These choice probabilities are used as the observed data $\bp^\cK$ to train the nonparametric s-RAM model. For prediction, we generate a new graph $\bG^{\text{new}}$ and compute the probabilities via \ref{data:shortestpath} with $\bG^{\text{new}}$. 

\textbf{Assignment Problem.} Recall that for the assignment problem, we assume the noise marginals are nonidentical uniform distributions. Using the parameter notations in the shortest path problem above,  from \citep{natarajan2009persistency}, we have that, for a graph $\bG^k$ with $k\in \cK$, the convex program with  nonidentical uniform marginal distributions for the assignment problem is as follows: 
\begin{align}
\begin{aligned}
    \max_{\bx} & \sum_{i\in \bV_L^k }\sum_{j \in \bV_R^k } (\nu_{ij}+ b_{ij})x_{ij} - \sum_{i\in \bV_L^k }\sum_{j \in \bV_R^k }  \frac{b_{ij} - a_{ij}}{2} x_{ij}^2  \\ 
   \text{s.t.}\;& \sum_{j \in \bV_R^k} x_{ij} = 1, \quad \forall i \in \bV_L^k,\\
     & \sum_{i \in \bV_L^k} x_{ij} = 1, \quad \forall j \in \bV_R^k,\\
    & x_{ij}\geq 0 \quad \forall \ i\in \bV_L^k, j \in \bV_R^k. 
\end{aligned}\label{data:assignment}
\end{align}
These are used as observed data $\bp^\cK$ for training the nonparametric s-RAM model. For prediction, we generate a new graph $\bG^{\text{new}}$ and compute the total expected system welfare via \ref{data:assignment} with $\bG^{\text{new}}$, which is the optimal objective value under the new graph.

%we evaluate whether the nonparametric s-RAM can accurately predict the total expected system welfare, defined as the same objective above with the new graph $\bG^\text{new}$ in the selected matching.

\textbf{Constrained Subset Selection Problem.} Recall that for the constrained subset selection problem, we assume the noise marginals are nonidentical exponential distributions. Given a collection historical subsets $\cI^{(1)}, \cI^{(2)},\cdots, \cI^{(K)}$, for a subset $\cI^{(k)}$ with $k\in \cK$, we solve a convex program defined by the corresponding parametric MDM model to obtain the choice probabilities. Let $\bv$ denote the deterministic weight matrix. Let $\lambda_{i}$ be the rate parameters for each item $i\in \cN$. From \citep{natarajan2009persistency}, we have that, for a subset $\cI^{(k)}$, the convex program with  nonidentical exponential marginal distributions for the constrained subset selection problem is as follows:
\begin{align}
\begin{aligned}
    \max_{\bx} & \sum_{i \in \cI^{(k)}}(\nu_{i}+\frac{1}{\lambda_{i}})x_{i} - \sum_{i \in \cI^{(k)}} \frac{1}{\lambda_{i}} x_{i}\ln{x_{i}}  \\ 
   \text{s.t.}\; & \sum_{i\in \cI^{(k)} } x_i \leq C,\\
    &  0 \leq x_{ij} \leq 1 \quad \forall \ i\in \cI^{(k)}.    
\end{aligned}\label{data:subsetselection}
\end{align}
These are used as observed data $\bp^\cK$ for training the nonparametric s-RAM model. For prediction, we generate a new graph $\bG^{\text{new}}$ and compute the marginal allocation weights via \ref{data:subsetselection} with $\bG^{\text{new}}$, which is the probability of the chosen item.

\subsection{Representability, Estimation, Prediction Formulations}
\label{app:example_theory}

\subsubsection{Shortest Path Problem} \label{app:shortestpath}

Suppose that $\bG=(s\cup t \cup \bV, \bE)$ is the complete graph that describes all possible nodes and edges of the network, and we have flow data for a collection $K$ graphs $\bG^1,\bG^2, \cdots, \bG^K,$ with $\bG^k = (s\cup t \cup \bV^k, \bE^k)$, for every $k\in \cK$. For a graph $\bG^k$, let $p_{ij}^k \in [0,1]$ denote the flow of edge $(i,j)$ in  $\bG^k$. Given observed flow data $\bp^\cK = (p_{ij}^k: (i,j) \in \bE^k, k \in \cK)$, by Theorem \ref{thm:mdm-general-feascon}, the s-RAM characterization for shortest path problems is given below. 

\textbf{s-RAM characterization for Shortest Path Problems.} Under Assumptions \ref{asp:convex_perturbation} and \ref{asp:convexhull}, a data collection $\bp^\cK$ is representable by an s-RAM  if and only if there exists $(\lambda_i^k \in \mathbb{R}: i\in s\cup t\cup \bV^k, k \in \mathcal{K})$,
    %$\blam \in \mathbb{R}^m$ 
    such that for any two distinct graph $k,l \in \cK$ containing a common edge $(i,j) \in \bE^k\cap \bE^l,$ 
    %for all $(i,S),(i,T)\in \cI$:
\begin{align}
\begin{aligned}
     &\lambda_j^k- \lambda_i^k < \lambda_j^l - \lambda_i^l    \quad \text{if} \quad p_{ij}^k < p_{ij}^l,\\
     &\lambda_j^k- \lambda_i^k = \lambda_j^l - \lambda_i^l     \quad \text{if} \quad 0 < p_{ij}^k = p_{ij}^l <1.   \label{sp-feascon}
\end{aligned}
\end{align}
As a result, checking whether the given flow data $\bp^\cK$ satisfies the MDM hypothesis can be accomplished by solving a linear program %whose size is polynomial in the number of products and assortments. 
with $\mathcal{O}(\vert \bV \vert K)$ continuous variables and $\mathcal{O}(\vert \bE \vert K) $ constraints as follows, which is equivalent to check whether the objective value of \eqref{consistencyLP_shorestpath} below is strictly positive.
\begin{align}
\begin{aligned}\label{consistencyLP_shorestpath}
    \max_{\epsilon,\blam}\quad & \epsilon\\
  \text{s.t.}\quad &\lambda_j^k- \lambda_i^k + \epsilon \leq \lambda_j^l - \lambda_i^l  \quad \text{if} \quad p_{ij}^k < p_{ij}^l,\quad \forall k,l \in \cK, \ \forall (i,j) \in \bE^k\cap \bE^l,\\
     & \lambda_j^k- \lambda_i^k  = \lambda_j^l - \lambda_i^l  \quad \text{if} \quad 0< p_{ij}^k = p_{ij}^l <1 ,\quad \forall k,l \in \cK, \ \forall (i,j) \in \bE^k\cap \bE^l.  
\end{aligned}
\end{align}

\textbf{Predictions for Shortest Path Problems.} Suppose that $\bp^\cK$ is representable by s-RAM, given a new graph $\bG^{\text{new}} = (s\cup t \cup \bV^{\text{new}}, \bE^{\text{new}})$, by Proposition \ref{prop:wc-obj-milp}, the prediction formulation to compute $\underline{f}_\text{new}$ for the shortest path problem by taking $0 < \epsilon < 1 /(2|\bV|K),$ is as follows. 
%Suppose that $\bp^\cK$ is representable by s-RAM and $A$ is a graph such that $A\notin \cK$. Then for any $0 < \epsilon < 1 /(2|\bV||\cK|),$ the worst-case objective value $\underline{f}_\text{new}$ equals the value of the following mixed integer linear program:   
%Suppose that $\bp^\cK$ is representable by MDM and $A$ is a graph such that $A\notin \cK$. Then  for any $0 < \epsilon < 1 /(2n|\cK|),$ the best-case objective value $\underline{r}_A$ equals the value of the following mixed integer linear program:   
\begin{align}
    \min_{\bx,\blambda,\bdelta\!\!\,}\quad&  \quad   f(\bx) \label{model:sp_milp} \\
\text{s.t.}\quad \ -&\delta_{i,j}^{\text{new},k} \ \leq\    \lambda_j^{\text{new}} - \lambda_i^{\text{new}} -  \lambda_j^k+ \lambda_i^k \ \leq \  1 -  (1+\epsilon)\delta_{i,j}^{\text{new},k},  \ \forall \,  (i,j) \in \bE^\text{new}, \, \forall k\in \cK, (i,j) \in \bE^k, \nonumber \\
 -&\delta_{i,j}^{k, \text{new}}  \ \leq\    \lambda_j^k- \lambda_i^k  -  \lambda_j^{\text{new}} + \lambda_i^{\text{new}}\ \leq \  1 -  (1+\epsilon)\delta_{i,j}^{k, \text{new}},  \ \forall \,  (i,j) \in \bE^\text{new}, \, \forall k\in \cK, (i,j) \in \bE^k, \nonumber \\
&  \delta_{i,j}^{\text{new},k} - 1 \ \leq \ p_{ij}^k - x_{ij}^{\text{new}}   \ \leq \ 1 - \delta_{i,j}^{k, \text{new}}, \quad\qquad\quad\quad\;\;  \forall \,  (i,j) \in \bE^\text{new}, \, \forall k\in \cK, (i,j) \in \bE^k, \nonumber  \\
 -&(\delta_{i,j}^{\text{new},k}+  \delta_{i,j}^{k, \text{new}}) \ \leq \ p_{ij}^k- x_{ij}^{\text{new}}  \ \leq \ \delta_{i,j}^{\text{new},k} + \delta_{i,j}^{k, \text{new}},\ \quad     \forall \,  (i,j) \in \bE^\text{new}, \, \forall k\in \cK, (i,j) \in \bE^k, \nonumber  \\ 
& \lambda_j^k- \lambda_i^k - \lambda_j^l + \lambda_i^l \ \geq \ \epsilon, \qquad\! \forall \, k,l \in \cK, \  \forall (i,j) \in \bE^k \cap \bE^l \ \textnormal{ s.t. }\  p_{ij}^k > p_{ij}^l, \nonumber  \\
&\lambda_j^k- \lambda_i^k - \lambda_j^l + \lambda_i^l \ = \ 0, \quad\ \ \, \forall \, k,l \in \cK, \  \forall (i,j) \in \bE^k \cap \bE^l \ \textnormal{ s.t. }\  0< p_{ij}^k = p_{ij}^l <1 , \nonumber  \\
&  \sum_{j:(i,j)\in \bE^\text{new}} x_{ij} - \sum_{j:(j,i)\in \bE^\text{new}} x_{ji} = \begin{cases}
         1,\quad &i=s, \\
         -1,\quad &i=t, \\
         0,\quad &i\in \bV^\text{new},
     \end{cases} \nonumber  \\
% & 0 \leq \lambda_j^{\text{new}} - \lambda_i^{\text{new}}\leq 1,\;\  x_{ij}^{\text{new}} \ \geq\  0,  \ \forall \,  (i,j) \in \bE^\text{new}, \\
% & 0 \, \leq\,  \lambda_j^k- \lambda_i^k \leq\,  1, \quad \forall \, k,l \in \cK, \  \forall (i,j) \in \bE^k \cap \bE^l ,\\ 
& x_{ij}\geq 0, \  \forall k\in \cK, (i,j) \in \bE^k,\nonumber \\
& \delta_{i,j}^{\text{new},k}, \delta_{i,j}^{k, \text{new}} \, \in \, \{0,1\}, \ \forall  (i,j) \in \bE^\text{new}, \, \forall k\in \cK, (i,j) \in \bE^k. \nonumber 
\end{align}
Likewise, the best-case estimate for the objective value $\bar{f}_\text{new}$ equals the optimal value obtained by maximizing over the constraints in the above mixed integer program. In our experiment, we take $f(\bx) = x_{\hat{e}}$ where $\hat{e}$ is the pre-specified edge. 

\textbf{Estimation for Shortest Path Problems}
In the case of shortest path problems, by taking $0 < \epsilon < 1 /(2|\bV|K),$ formulation \eqref{model:mdm_micp} can be written as: 
\begin{align}
\begin{aligned}
    \min_{\bx,\blambda,\bdelta\!\!\,}\quad&  \quad   \textnormal{loss}(\bp_\cK,\bx_\cK) \\
    \text{s.t.}\quad \ -&\delta^{k,l}_{i,j} \ \leq\    \lambda_j^k- \lambda_i^k  -  \lambda_j^l + \lambda_i^l \ \leq \  1 -  (1+\epsilon)\delta^{k,l}_{i,j}, \quad\  \ \forall \,   k,l \in \cK,\ \forall (i,j) \in \bE^{k}\cap \bE^{l}, \\
&  \delta^{k,l}_{i,j} - 1 \ \leq \  x_{ij}^{l}  -x_{ij}^{k}  \ \leq \ 1 - \delta^{l,k}_{i,j}, \quad\qquad\quad\,   \forall \,   k,l \in \cK,\ \forall (i,j) \in \bE^{k}\cap \bE^{l},\\
 -&(\delta^{k,l}_{i,j} +  \delta^{l,k}_{i,j} ) \ \leq \ x_{ij}^{l}  -x_{ij}^{k}    \ \leq \ \delta^{k,l}_{i,j} + \delta^{l,k}_{i,j} ,\ \quad     \forall \,   k,l \in \cK,\ \forall (i,j) \in \bE^{k}\cap \bE^{l},\\ 
&  \sum_{j:(i,j)\in \bE^k} x_{ij} - \sum_{j:(j,i)\in \bE^k} x_{ji} = \begin{cases}
         1,\quad &i=s, \\
         -1,\quad &i=t, \\
         0,\quad &i\in \bV^k,
     \end{cases} \forall k \in \cK,\\
%& 0 \, \leq\,  \lambda_j^k- \lambda_i^k \leq\,  1, \quad \forall \, k,l \in \cK, \  \forall (i,j) \in \bE^k \cap \bE^l ,\\
& x_{ij}\geq 0, \  \forall \,  k \in \cK,\  \forall (i,j) \in \bE^k,\\
&\delta^{k,l}_{i,j} \, \in \, \{0,1\}, \  \forall \,  k,l \in \cK,\  \forall (i,j) \in \bE^k \cap \bE^l. \label{model:sp_micp}
\end{aligned}
\end{align}

\subsubsection{Assignment Problem} \label{app:assignment}

Suppose that $\bG=(\bV_L, \bV_R, \bE)$ is the complete bipartite graph that describes all the possible nodes and edges of the network, and we have flow data for a collection $K$ graphs $\bG^1,\bG^2, \cdots, \bG^K,$ with $\bG^k = (\bV_L^k, \bV_R^k, \bE^k)$, for every $k\in \cK$. For a graph $\bG^k$, let $p_{ij}^k \in [0,1]$ denote the flow of edge $(i,j)$ in  $\bG^k$. Given observed flow data $\bp^\cK = (p_{ij}^k: (i,j) \in \bE^k, k \in \cK)$, by Theorem \ref{thm:mdm-general-feascon}, the s-RAM characterization for assignment problems is given below. 
% More formally, suppose that $\bG=(\bV_L, \bV_R, \bE)$ is the complete bipartite graph that describes all the possible nodes and edges of the network, and we have flow data for a collection $\mathcal{S}$ of graphs of $\bG.$ For each graph $k \in \mathcal{K},$ let $p_{ij}^k \in [0,1]$ denote the fraction of choosing edge $(i,j)$ in graph $k$. Our goal is to identify necessary and sufficient conditions for the data $\bp^\cK = (p_{ij}^k: (i,j) \in \bE^k, k \in \cK)$  to be representable by an s-RAM instance.  We are interested in identifying necessary and sufficient conditions on the observable flow data  $\bp^\cK$ which make it consistent with the MDM hypothesis. Applying the result in Theorem \ref{thm:mdm-general-feascon}, 

\textbf{s-RAM characterization for Assignment Problems.} Under Assumptions \ref{asp:convex_perturbation} and \ref{asp:convexhull}, an observed flow collection $\bp^\cK$ is representable by an s-RAM  if and only if there exists $(\lambda_{i,L}^k, \lambda_{j,R}^k \in \mathbb{R}: i\in \bV^k_L, j\in \bV^k_R, k \in \mathcal{K})$, such that for any two distinct graph $k,l \in \cK$ containing a common edge $(i,j) \in \bE^k \cap \bE^l,$ 
\begin{align}
\begin{aligned}
     &\lambda_{i,L}^k + \lambda_{j,R}^k  > \lambda_{i,L}^l + \lambda_{j,R}^l \quad \text{if} \quad p_{ij}^k > p_{ij}^l,\\
     &\lambda_{i,L}^k + \lambda_{j,R}^k = \lambda_{i,L}^l + \lambda_{j,R}^l \quad \text{if} \quad 0< p_{ij}^k = p_{ij}^l <1.   \label{eq:mdm-feascon-assignment}
\end{aligned}
\end{align}
As a result, checking whether the given flow data $\bp^\cK$ satisfies the MDM hypothesis can be accomplished by solving a linear program with $\mathcal{O}(\vert \bV \vert K)$ continuous variables and $\mathcal{O}(\vert \bE \vert K) $ constraints as follows, which is equivalent to check whether the optimal objective value of \eqref{consistencyLP_assignment} below is strictly positive:
\begin{align}
\begin{aligned}\label{consistencyLP_assignment}
    \max_{\epsilon,\blam}\quad & \epsilon\\
  \text{s.t.}\quad &\lambda_{i,L}^k + \lambda_{j,R}^k   \geq \lambda_{i,L}^l + \lambda_{j,R}^l + \epsilon \quad \text{if} \quad p_{ij}^k > p_{ij}^l,\ \forall k,l \in \cK, \ \forall (i,j) \in \bE^k\cap \bE^l,\\
     & \lambda_{i,L}^k + \lambda_{j,R}^k = \lambda_{i,L}^l + \lambda_{j,R}^l  \quad \text{if} \quad 0 < p_{ij}^k = p_{ij}^l < 1,\ \forall k,l \in \cK, \ \forall (i,j) \in \bE^k\cap \bE^l.  
\end{aligned}
\end{align}

\textbf{Predictions for Assignment Problems.} Suppose that $\bp^\cK$ is representable by s-RAM, given a new graph $\bG^{\text{new}} = (\bV_L^{\text{new}}, \bV_R^{\text{new}}, \bE^{\text{new}})$, by Proposition \ref{prop:wc-obj-milp}, the prediction formulation to compute $\underline{f}_\text{new}$ for the assignment problem by taking $0 < \epsilon < 1 /(2|\bV|K),$ is as follows.    
\begin{align}
    \min_{\bx,\blambda,\bdelta\!\!\,}\quad&  \quad  f(\bx) \label{model:assignment_milp} \\
\text{s.t.}\quad \ -&\delta_{i,j}^{\text{new},k} \ \leq\    \lambda_{i,L}^{\text{new}} + \lambda_{j,R}^{\text{new}}   - \lambda_{i,L}^k -  \lambda_{j,R}^k  \ \leq \  1 -  (1+\epsilon)\delta_{i,j}^{\text{new},k}, \ \forall   (i,j) \in \bE^\text{new},  \forall k\in \cK, (i,j) \in \bE^k, \nonumber  \\
 -&\delta_{i,j}^{k, \text{new}}  \ \leq\    \lambda_{i,L}^k + \lambda_{j,R}^k  - \lambda_{i,L}^{\text{new}}  -  \lambda_{j,R}^{\text{new}} \ \leq \  1 -  (1+\epsilon)\delta_{i,j}^{k, \text{new}}, \   \forall  (i,j) \in \bE^\text{new},  \forall k\in \cK, (i,j) \in \bE^k,\nonumber  \\
&  \delta_{i,j}^{\text{new},k} - 1 \ \leq \ p_{ij}^k  -  x_{ij}^{\text{new}} \ \leq \ 1 - \delta_{i,j}^{k, \text{new}}, \quad\qquad\qquad\;\;   \forall \,  (i,j) \in \bE^\text{new}, \, \forall k\in \cK, (i,j) \in \bE^k, \nonumber  \\
 -&(\delta_{i,j}^{\text{new},k}+  \delta_{i,j}^{k, \text{new}}) \ \leq \ p_{ij}^k  - x_{ij}^{\text{new}} \ \leq \ \delta_{i,j}^{\text{new},k} + \delta_{i,j}^{k, \text{new}},\ \quad     \forall \,  (i,j) \in \bE^\text{new}, \, \forall k\in \cK, (i,j) \in \bE^k, \nonumber  \\ 
& \lambda_{i,L}^k + \lambda_{j,R}^k  - \lambda_{i,L}^l - \lambda_{j,R}^l \ \geq \ \epsilon, \qquad\! \forall \, k,l \in \cK, \  \forall (i,j) \in \bE^k \cap \bE^l \ \textnormal{ s.t. }\  p_{ij}^k > p_{ij}^l,  \nonumber \\
&\lambda_{i,L}^k  + \lambda_{j,R}^k - \lambda_{i,L}^l - \lambda_{j,R}^l \ = \ 0, \quad\ \ \, \forall \, k,l \in \cK, \  \forall (i,j) \in \bE^k \cap \bE^l \ \textnormal{ s.t. }\  0< p_{ij}^k = p_{ij}^l <1, \nonumber  \\
& \sum_{j \in \bV_R^{\text{new}}} x_{ij} = 1, \quad \forall i \in \bV_L^{\text{new}},  \nonumber  \\
& \sum_{j \in \bV_L^{\text{new}}} x_{ij} = 1, \quad \forall i \in \bV_R^{\text{new}},  \nonumber  \\
& x_{ij} \geq 0 ,  \   \ \forall \,  (i,j) \in \bE^\text{new} ,\\
& \delta_{i,j}^{\text{new},k}, \delta_{i,j}^{k, \text{new}} \, \in \, \{0,1\}, \   \ \forall \,  (i,j) \in \bE^\text{new} , \, \forall k\in \cK, (i,j) \in \bE^k. \nonumber
\end{align}
Likewise, the best-case estimate for the objective value $\bar{f}_\text{new}$ equals the optimal value obtained by maximizing over the constraints in the above mixed integer linear program. In our experiment, to evaluate the total welfare of the system under the equivalent parametric MDM with marginals being nonidentical uniform distributions, we take $f(\bx) = \sum_{i\in \bV_L^\text{new} }\sum_{j \in \bV_R^\text{new}  } (\nu_{ij}+ b_{ij})x_{ij} + \sum_{i\in \bV_L^\text{new}  }\sum_{j \in \bV_R^\text{new}} \frac{b_{ij} - a_{ij}}{2} x_{ij}^2$ where the uniform distributions are in the form of $\text{Uniform}(a_{ij}, b_{ij})$ for $(i,j)\in \bE^\text{new}$. 

\textbf{Estimation for Assignment Problems.} In the case of assignment problems, by taking $0 < \epsilon < 1 /(2|\bV|K),$ formulation \eqref{model:mdm_micp} can be written as: 
\begin{align}
    \min_{\bx,\blambda,\bdelta\!\!\,}\quad&  \quad   \textnormal{loss}(\bp_\cK,\bx_\cK) \nonumber \\
\text{s.t.}\quad \ 
 -&\delta_{i,j}^{k,l}  \ \leq\    \lambda_{i,L}^k + \lambda_{j,R}^k  - \lambda_{i,L}^l  -  \lambda_{j,R}^l\ \leq \  1 -  (1+\epsilon)\delta_{i,j}^{k,l},   \ \forall \, k,l \in \cK,\ \forall (i,j) \in \bE^k \cap \bE^l, \nonumber \\
&  \delta_{i,j}^{k,l} - 1 \ \leq \  x_{ij}^{l} - x_{ij}^{k}    \ \leq \ 1 - \delta_{i,j}^{l,k}, \quad\qquad\quad  \forall \,   k,l \in \cK,\  \forall (i,j) \in \bE^k \cap \bE^l, \nonumber  \\
 -&(\delta_{i,j}^{k,l} +  \delta_{i,j}^{l,k}) \ \leq \ x_{ij}^{l} - x_{ij}^{k}   \ \leq \ \delta_{i,j}^{k,l} + \delta_{i,j}^{l,k},\ \quad     \forall \,   k,l \in \cK,\ \forall (i,j) \in \bE^{k}\cap \bE^{l},  \nonumber  \\ 
& \sum_{j \in \bV_R^k} x_{ij} = 1, \quad \forall k \in \cK, \forall i \in \bV_L^k,  \label{model:assignment_micp} \\
& \sum_{j \in \bV_L^k} x_{ij} = 1, \quad \forall k \in \cK, \forall i \in \bV_R^k,  \nonumber \\
%& 0 \, \leq\,  \lambda_{i,L}^k + \lambda_{j,R}^k  \leq\,  1, \quad \forall \, k,l \in \cK, \  \forall (i,j) \in \bE^k \cap \bE^l ,\\
& x_{ij} \geq 0 ,\   \ \forall \,  k \in \cK,\ \forall (i,j) \in \bE^{k},\\
& \delta_{i,j}^{k,l} \, \in \, \{0,1\}, \   \ \forall \,  k,l \in \cK,\ \forall (i,j) \in \bE^{k}\cap \bE^{l}. \nonumber
\end{align}

\subsubsection{Constrained Subset Selection Problem} \label{app:subsetselection}

Suppose $\cN$ is the complete set of available items. We focus on the cardinality-constrained case with $\mathcal{X} = \{\bx: \sum_{i\in \cN} x_i \leq C, \ x_{i}\in\{0,1\},\forall i \in \cN \}$,  where $C$ is a positive integer such that $2 \leq C \leq |\cN| - 1$. We have marginal allocation weight data for a collection of $K$ observations with available subsets $\cI^{(1)}, \cI^{(2)}, \ldots, \cI^{(K)}$ such that $\cI^{(k)} \subseteq \cN$, for $k \in \cK$. For a set $\cI^{(k)}$, let $p_i^k \in [0,1]$ denote the marginal allocation weight for item $i$ in $\cI^{(k)}$. Given observed marginal allocation weight data $\bp^\cK = (p_i^k:i\in \cI^{(k)},k\in \cK)$, by Theorem~\ref{thm:mdm-general-feascon}, the s-RAM characterization for the constrained subset selection problem is given below.

% Suppose $\cN$ is the complete set of available items. We focus on the cardinality-constrained case with $\mathcal{X} = \{\bx: \sum_{i\in \cN} x_i \leq C, \ x_{i}\in\{0,1\},\forall i \in \cN \}$, where $C$ is a postive integer such that $2\leq C\leq |\cN| - 1$.
% We have marginal allocation weight data for a collection $K$ observations with avaiable subsets $\cI^{(1)}, \cI^{(2)},..., \cI^{(K)}$ such that $\cI^{(k)}\subseteq \cN$, for $k\in\cK.$ For a set  $\cI^{(k)}$, let $\bp_i^k\in [0,1]$ denote the mariginal allocation weight for item $i$ in $\cI^{(k)}$. Given observed marginal allocation weight data $\bp^\cK = (p_i^k:i\in \cI^{(k)},k\in \cK)$, by Theorem \ref{thm:mdm-general-feascon}, the s-RAM characterization for the constrained subset selection problem is given below. 

\textbf{s-RAM characterization for Constrained Subset Selection Problems with Cardinality Constraint.} Under Assumptions \ref{asp:convex_perturbation} and \ref{asp:convexhull}, a data  collection $\bp^\cK$ is representable by an s-RAM if and only if there exists $(\lambda^k\in \mathbb{R}: k\in \cK)$ such that for any two observation $k,l \in \cK$ containing a common item $i \in \cN,$ 
\begin{align}
\begin{aligned}
     &\lambda_k  > \lambda_l \quad \text{if} \quad p_i^k < p_i^l,\\
     & \lambda_k  > \lambda_l   \quad \text{if}\quad  p_i^k = p_i^l \neq 0.   \label{eq:subsetselection-feascon}
\end{aligned}
\end{align}
As a result, checking whether the given choice data $\bp^\cK$ satisfies the s-RAM hypothesis can be accomplished by solving a linear program with $\mathcal{O}(\vert \cK \vert)$ continuous variables and $\mathcal{O}(n \vert \cK \vert) $ constraints as follows, which is equivalent to checking whether the optimal objective value of \eqref{model:subsetselection_feas_lp} below is strictly positive.
\begin{align}
\begin{aligned}
\max_{\epsilon,(\lambda_k:k\in\cK)}  \quad & \eps &&   &&& \\
 \text{s.t.}\qquad &\lambda_k  \geq \lambda_l  + \epsilon &&\text{if }\quad p_{i}^k < p_{i}^l,  &&&\forall k,l \in \cK, i\in \cI^{(k)}\cap \cI^{(l)} ,\\
  &\lambda_k  = \lambda_l  &&\text{if } \quad p_{i}^k = p_{i}^l \neq 0,  &&&\forall k,l \in \cK, i\in \cI^{(k)}\cap \cI^{(l)}.
\end{aligned}
\label{model:subsetselection_feas_lp}
\end{align}

\textbf{Prediction for Constrained Subset Selection Problems.} Suppose $\bp^\cK$ is representable by s-RAM, given a new set $\cI^{\text{new}}$, by Proposition \ref{prop:wc-obj-milp}, the prediction formulation to compute $\underline{f}_{\text{new}}$ for the constrained subset selection problem by taking  any $0 < \epsilon < 1 /(2\vert \cK \vert)$ is as follows:   
\begin{align}
\begin{aligned}
\min_{\bx,\blambda,\bdelta\!\!\,}\quad&  \quad  f(\bx) \\
 \text{s.t.}\quad \ -&\delta_{\text{new},k} \ \leq\    \lambda_{\text{new}}  -  \lambda_k \ \leq \  1 -  (1+\epsilon)\delta_{\text{new},k}, \qquad\quad\  \ \forall k\in \cK, \forall \,  i \in \cI^{\text{new}} \cap \cI^{(k)},  \\
 -&\delta_{k,\text{new}}  \ \leq\    \lambda_k  -  \lambda_{\text{new}} \ \leq \  1 -  (1+\epsilon)\delta_{k,\text{new}}, \qquad\quad\  \ \forall k\in \cK, \forall \,  i \in \cI^{\text{new}} \cap \cI^{(k)}, \\
&  \delta_{\text{new},k} - 1 \ \leq \  x_{i} - p_{i}^k \ \leq \ 1 - \delta_{k,\text{new}}, \quad\qquad\quad\,   \forall k\in \cK, \forall \,  i \in \cI^{\text{new}} \cap \cI^{(k)}, \\
 -&(\delta_{\text{new},k}  +  \delta_{k,\text{new}}) \ \leq \ x_{i,A} - p_{i,S} \ \leq \ \delta_{\text{new},k}  + \delta_{k,\text{new}},\ \quad     \forall k\in \cK, \forall \,  i \in \cI^{\text{new}} \cap \cI^{(k)},\\ 
 & \lambda_k - \lambda_l \ \geq \ \epsilon, \qquad\qquad\qquad\! \forall \, k,l \in \cK, i\in \cI^{(k)}\cap \cI^{(l)}, \ \textnormal{ s.t. }\  p_{i}^k < p_{i}^l ,\\
&\lambda_k - \lambda_l \ = \ 0, \qquad\qquad\quad\ \ \, \forall \, k,l \in \cK, i\in \cI^{(k)}\cap \cI^{(l)}, \textnormal{ s.t. }\  p_{i}^k  = p_{i}^l \neq 0 ,\\
\
& \sum_{i \in \cI^{\text{new}}} x_{i} \ = \ C,   \\
& 0 \leq \lambda_A \leq 1,\quad\  x_{i} \ \geq\  0,  \ \forall \,  i \in \cI^{\text{new}} ,\quad\ 0 \, \leq\,  \lambda_k \leq\,  1, \quad  \delta_{\text{new},k}, \delta_{k,\text{new}} \, \in \, \{0,1\}, \  \forall \, k \in \cK. 
\end{aligned}\label{model:subsetselection_predition}
\end{align}
Likewise, the best-case estimate for the objective value $\bar{f}_{\text{new}}$ equals the optimal objective value obtained by maximizing over the constraints in the above mixed integer program.

\textbf{Estimation for Constrained Subset Selection Problems.} Suppose that Assumptions \ref{asp:convex_perturbation} and \ref{asp:convexhull} is satisfied. Then for any $0 < \epsilon < 1 /(2\vert \cK \vert),$ formulation \ref{model:mdm_micp} can be written as:
\begin{align}
\begin{aligned}
& \min_{\bx,\blam,\boldsymbol{\delta}}\quad   \textnormal{loss}(\bp^\cK,\bx^\cK)  \\
& \text{s.t.}\quad   -\delta_{k,l}  \leq   \lambda_{k} - \lambda_{l} \leq   1 -  (1+\epsilon)\delta_{k,l}, \quad\quad\quad\;\;\; \forall k,l \in \cK, i\in \cI^{(k)} \cap \cI^{(l)}, \\
& \quad\quad\quad \delta_{k,l} - 1  \leq  x_{i}^k-x_{i}^l  \leq 1 - \delta_{l,k}, \quad\quad\quad\quad\    \forall k,l \in \cK, i\in \cI^{(k)} \cap \cI^{(l)}, \\
&\quad\quad-(\delta_{k,l} + \delta_{l,k})  \leq x_{i}^k - x_{i}^l  \leq  \delta_{k,l} + \delta_{l,k}, \quad\, \forall k,l \in \cK, i\in \cI^{(k)} \cap \cI^{(l)}, \\
 &\quad\quad\quad \sum_{i\in \cI^{(k)}} x_{i}^k \leq C, \ \forall \,k\in \cK, \;  0\leq x_{i}^k \leq 1,\  \forall\, k\in \cK, i\in \cI^{(k)}, \\
 &\quad\quad\quad0 \leq \lambda_{k} \leq 1,\ \forall\, k\in \cK, \;\\
 &\quad\quad\quad\delta_{k,l} \in \{0,1\},\   \forall\, k,l \in \cK.   
\end{aligned}\label{model:subsetselection_estimation}
\end{align}

\end{document}